\documentclass[11pt,reqno]{amsproc}

\title[Scattering for focusing supercritical wave equations in four dimensions]{Scattering for radial bounded solutions of focusing supercritical wave equations in four dimensions}

\author[G.~Camliyurt]{Guher Camliyurt}
\address{Department of Mathematics, Virginia Tech, Blacksburg, VA 24061}
\email{gcamliyurt@vt.edu}

\author[C. E.~Kenig]{Carlos E. Kenig}
\address{Department of Mathematics, University of Chicago, Chicago, IL 60637}
\email{cek@math.uchicago.edu}

\usepackage[margin=1in]{geometry}
\usepackage{amsmath, amsthm, amssymb}
\usepackage{mathrsfs}
\usepackage{graphicx}
\usepackage{times}
\allowdisplaybreaks

\usepackage[usenames,dvipsnames]{color}
\usepackage[colorlinks=true, pdfstartview=FitV, linkcolor=blue, citecolor=blue, urlcolor=blue]{hyperref}
\usepackage{amsmath}
\usepackage{amsfonts}
\usepackage{amssymb}

\usepackage{amsmath, amsthm, amssymb, enumerate}

\usepackage{fancyhdr}
\usepackage{mathtools}
\usepackage{bbm}

\DeclarePairedDelimiter\evaluat{.}{\rvert}

\numberwithin{equation}{section}

\definecolor{colorgggg}{rgb}{0.5,0.0,0.4}
\def\cole{}

\def\colb{\color{black}}
\def\inon#1{\,\,\,\,\,\,\hbox{#1}}                

\def\epsi{\epsilon}    
\def\norm#1{\left\Vert #1\right\Vert} 

\def\llabel#1{\nonumber}

\newtheorem{Theorem}{Theorem}[section]
\newtheorem{Corollary}[Theorem]{Corollary}
\newtheorem{Proposition}[Theorem]{Proposition}
\newtheorem{Lemma}[Theorem]{Lemma}
\newtheorem{Remark}[Theorem]{Remark}

\newtheorem{Notation}[Theorem]{Notation}
\newtheorem{Claim}[Theorem]{Claim}
\newtheorem*{Claim*}{Claim}
\newtheorem*{Proposition*}{Proposition}
\newtheorem*{Lemma*}{Lemma}
\newtheorem{definition}[Theorem]{Definition}

\def\comma{ {\rm ,\qquad{}} }            
\def\divv{\mathop{\rm div}\nolimits}    
\def\supp{\mathop{\rm supp}\nolimits}    
\def\les{\lesssim}

\def\indeq{\quad{}}                     
\def\indeqtimes{\indeq\indeq\times}
\def\period{.}                           

\def\Imax{\textmd{I}_{\scriptsize{\mbox{max}}}}
\def\bigO{\mathcal{O}}

\begin{document}

\nocite{*} 

\begin{abstract}
We consider the focusing wave equation with energy supercritical nonlinearity in dimension four. 
We prove that any radial solution that remains bounded in the critical Sobolev space is global and 
scatters to free waves as $t \to \pm \infty$. 
\end{abstract}



\thanks{{\em Date}. \today}

\maketitle

\tableofcontents

\section{Introduction}
\label{sec:intro}
We consider the semilinear focusing wave equation 
\begin{equation}
     \partial^2_{t} u
    - \Delta u 
    -  |u|^{p-1} u
     = 0   \comma (t,x) \in  I \times \mathbb{R}^4 ,
   \label{wave1.1}
 \end{equation}
  with initial data
  \begin{equation}
  u (0, x) = u_0 (x) \comma \partial_t u (0, x) = u_1 (x).
      \label{wave1.2}
  \end{equation}
The set $I$ in \eqref{wave1.1} is an interval with $0 \in I$. Equation~\eqref{wave1.1} has the following scaling invariance: if $u (t,x)$ is a solution of \eqref{wave1.1}, then 
\begin{equation}
u_{\lambda} \left( \lambda t, x \right) = \lambda^{\frac{2}{p-1}} u (\lambda t, \lambda x)
    \label{wave1.3}
\end{equation}
for $\lambda >0$ is also a solution. The scaling in \eqref{wave1.3}
determines the critical regularity for \eqref{wave1.1}--\eqref{wave1.2} as
\begin{equation*}
    s_p = \frac{d}{2} - \frac{2}{p-1}
\end{equation*}
due to the isometry of the norms
\begin{equation*}
    \norm{\left( 
    u_{\lambda} (\lambda t, \cdot), \partial_t u_{\lambda} (\lambda t, \cdot)
    \right)}_{\dot{H}^{s_p} \times \dot{H}^{s_p -1}}
    = 
    \norm{ \left(
    u (\lambda t, \cdot), \partial_t
    u (\lambda t, \cdot)
    \right)
    }_{\dot{H}^{s_p} \times \dot{H}^{s_p -1}} .
\end{equation*}
Here, $\dot{H}^{s_p}$ 
denotes the closure of the Schwartz functions under the 
 \textit{homogeneous Sobolev norm }
 \begin{equation*}
 \norm{f}_{\dot{H}^{s_p}  (\mathbb{R}^4)} = \norm{D^{s_p} f}_{L^2 (\mathbb{R}^4)} .
 \end{equation*}
In this paper we study radial solutions of the initial value problem \eqref{wave1.1}--\eqref{wave1.2} with energy supercritical nonlinearities $p>3$, i.e., $s_p >1$. We prove that any radial solution whose norm stays bounded in the scale invariant space  $\dot{H}^{s_p} \times \dot{H}^{s_p -1}$ throughout its maximal interval of existence must be global and must scatter to a linear solution. The precise statement of our main result is given below. 
\begin{Theorem}
\label{main thm_1}
Let $p \geq 5$ be an odd integer and let 
$\vec {u} (t)= (u(t), \partial_t u(t)) $  be a radial solution to the Cauchy problem
\begin{align}
   \begin{split}
    & \partial^2_{t} u
    - \Delta u 
    -  |u|^{p-1} u
     = 0   \comma \mbox{in } I \times \mathbb{R}^4 ,
    \\
   &  \vec{u} (0)
    = (u_0, u_1) \in \dot{H}^{s_p} \times \dot{H}^{s_p -1} (\mathbb{R}^4)
    \end{split}
\label{wave equations}
\end{align}
with maximal interval of existence $\Imax (\vec{u}) = (T_{-} (\vec{u}), T_{+} (\vec{u}))$ such that 
\begin{align}
\sup_{t \in (0, T_+ (\vec{u}))}   \Vert (u(t), \partial_t u(t)) \Vert_{\dot{H}^{s_p} \times \dot{H}^{s_p -1} (\mathbb{R}^4)} 
< \infty
.
\label{wave1.4}
\end{align}
Then, $\Imax ( \vec{u}) \cap (0, \infty)  = (0, \infty)$  and $\vec{u} (t)$ scatters to a free wave as $t \to \infty$. 
\end{Theorem}
The analogous result was first established in three dimensions by Duyckaerts, Kenig, and Merle 
\cite{DKM14}, extended to five dimensions by 
Dodson and Lawrie \cite{DL15}, and further generalized to higher odd dimensions (greater than or equal to seven) by the authors of the present paper \cite{CK23}. We remark that  
Theorem~\ref{main thm_1} and analogous versions in \cite{DKM14, DL15, CK23} directly imply that any finite time blow-up solution $\vec{u} (t)$ must be accompanied with a sequence of times $\{ t_n \}$ approaching $T_+ (\vec{u})$ so that
\begin{equation*}
    \limsup_{t_n \to T_+ (\vec{u})} \Vert (u(t), \partial_t u(t)) \Vert_{\dot{H}^{s_p} \times \dot{H}^{s_p -1}} = \infty . 
\end{equation*}
In dimension three, these issues were studied in generalized scale-invariant spaces that are non-Hilbertian by Duyckaerts and Roy in \cite{DR17} and by Duyckaerts and Yang in \cite{DY18}. We remark that similar conditional scattering results, such as Theorem~\ref{main thm_1}, were first obtained in the case of defocusing nonlinearities: see \cite{KM11, KV11, B14} for analogous results;  see also  \cite{B12, B15, KV11.2} for results in the non-radial setting.

The boundedness condition \eqref{wave1.4} provides a priori control over a scaling-critical norm of the solution, which in turn serves as a crucial tool for analyzing its qualitative behavior. This criterion helps compensate for the absence of conservation laws and monotonicity formulae in energy-supercritical regimes.
The condition \eqref{wave1.4} is also closely related to the definition of type-II solutions. Namely, 
a solution $\vec{u} (t)$ is classified as type-II if its scaling-critical norm remains bounded over $(T_- (\vec{u}), T_+ (\vec{u}))$. This definition is motivated by type-I blow-up phenomena. To give an example, we observe that
\begin{equation*}
    \varphi_T (t) = \frac{C(p)}{(T-t)^{\frac{2}{p-1}}}
    \comma C (p) = \left( \frac{2 (p+1)}{(p-1)^2} \right)^{\frac{1}{p-1}}
\end{equation*}
solves the ODE $\partial_{t}^2 \varphi = \varphi^p$ and blows up at time $t=T$. By setting a truncated version of $\varphi_T$ as initial data, we may also obtain a compactly supported solution $u_T$ to \eqref{wave1.1}. 
By the finite speed of propagation, $u_T (t,x)$ admits an unbounded critical Sobolev norm, i.e.,  
\begin{equation*}
    \lim_{t \to T} \norm{ u_T (t, x) }_{\Vert_{\dot{H}^{s_p} \times \dot{H}^{s_p -1}} } = \infty . 
\end{equation*}

For the focusing energy-critical equation in dimensions $d \geq 3$
\begin{align}
 \begin{split}
    & \partial^2_{t} u
    - \Delta u 
    =   |u|^{\frac{4}{d-2}} u
       \comma \mbox{in } I \times \mathbb{R}^d ,
    \\
   &  \vec{u} (0)
    = (u_0, u_1) \in \dot{H}^{1} \times L^2 (\mathbb{R}^d)
    \end{split}
   \label{wave equations_ec}
\end{align}
the type-II category includes a wide range of solutions. 
In \cite{KST09}, Krieger, Schlag, and Tataru provided the first construction of a radial type-II 
 solution to \eqref{wave equations_ec} 
in dimension $3$ that blows up in finite time. Their finite time blow-up construction 
$u(t,x)$ has the form 
\begin{equation}
u (t,x) = \lambda (t)^{\frac{1}{2}} W (\lambda (t) x) + \eta (t, x)
\end{equation}
near the blow-up time $T_- =0$. Here, $W$ is the unique (up to sign and scaling) radial ground state solution for the underlying elliptic equation. For similar constructions, see also \cite{KS14, DHKS14}. 

In \cite{DKM13}, Duyckaerts, Kenig, and Merle  proved that any non-scattering radial type-II solution to \eqref{wave equations_ec} in dimension $3$ 
decouples asymptotically as a sum of solitary waves and a radiation term in the energy space. This type of behavior is typically referred to as \textit{soliton resolution}.  
The proof of soliton resolution was also established in odd dimensions $d\geq 5$ in \cite{DKM23}. The dimensions $d=4$ and $d=6$ were treated in \cite{DKMM22} and \cite{CDKM24}, respectively.  An essential ingredient of these results has come to be known as the \textit{channels of energy method}. In \cite{JL23}, following an independent approach, Jendrej and Lawrie gave a proof  soliton resolution for radial type-II solutions  in every dimension $d \geq 4$.


 \subsection{Outline of the proof of Theorem~\ref{main thm_1}} 

 As in \cite{KM11, DKM14, DL15, CK23}, 
the proof of Theorem~\ref{main thm_1} is built upon
 the concentration compactness technics and rigidity methods. In Sections~\ref{sec:A review of the Cauchy problem}--\ref{sec:Decay}, we recall the familiar aspects
of the proofs given for analogous results in odd dimensions \cite{KM11, DKM14, DL15, CK23}. In Section~\ref{sec:A review of the Cauchy problem}, we recall some preliminaries and introduce the function spaces for the local theory for the Cauchy problem \eqref{wave equations}.  We then review the standard and inhomogeneous Strichartz estimates and extend the local wellposedness results of \cite[Section~2]{CK23}. In Section~\ref{sec:reduction}, we outline the details behind the concentration compactness method following closely the discussion in \cite{KM10}. In Section~\ref{sec:Decay}, we adapt the strategy behind the double Duhamel method of \cite{CKSTT08}  as in \cite[Section~4]{CK23} to show that solutions with the compactness property exhibit better decay estimates.

Although Theorem~\ref{main thm_1} may appear to be a natural extension of global wellposedness and scattering results of \cite{DKM14, DL15, CK23}, proving these assertions in even dimensions require tools and technics beyond those used in odd dimensional problems. The second half of the proof adapts a road map introduced in \cite{DKMM22}, where the rigidity methods are tailored to address the challenges that are unique to four dimensions.

The main difference between the treatments in even and odd dimensions stems from the 
lack of the strong Huygens principle and its effect on
applicability of the channels of energy method, which is closely tied to the parity of the
space dimension $d$. We refer our readers to \cite[Subsection~1.2]{DKMM22} for a comprehensive analysis  of channels of energy arguments and their  use in rigidity methods. In this paper, we rely on the following channels of energy estimate from \cite[Lemma~3.9]{DKMM22} (see also \cite{LSW24}):
\begin{Lemma*}
    Let $v_0 \in \dot{H}^1$ and $f \in L^1 ((0, \infty), L^2 (\mathbb{R}^4))$ have radial symmetry. Let $\vec{v} (t,r)$ solve the linear inhomogeneous equation
    \begin{equation*}
        \partial_t^2 v - \Delta v = f \comma (t,x) \in \mathbb{R} \times \mathbb{R}^4
    \end{equation*}
    with initial data $(v_0, 0)$. Then we have     
\begin{equation}
  \Vert{ \pi_R^{\perp}v_0} \Vert_{\dot{H}^1 ({R})}
    \leq \sqrt{\frac{20}{3}}
    \left( \lim_{t \to + \infty} \norm{\nabla v (t)}_{L^2 (|x|> t +R)} 
    + \norm{ f}_{L^1 ((0, \infty), \, L^2 (|x|>t + R)  )} 
    \right)
    \label{eq5.13intro}
\end{equation}
for every $R>0$.
\end{Lemma*}
The lower bound on the left-hand side above measures the part of $v_0$ orthogonal to ${1}/{r^2}$ in $\dot{H}^1 (R)$.   

Another related  challenge is faced in even dimensions due to the failure of the strong Huygens principle. More precisely, let 
$v_L$ be a solution of the linear wave equation with data $(v_0, v_1)$. In odd dimensions, if 
$\supp ((v_0, v_1)) \subset B(x_0, a)$,  then we have
\begin{equation}
    \supp (v_L (t)) \subset \{x: t- a \leq |x-x_0| \leq t+a \}
    \label{eq0.10intro}
\end{equation}
for  $t>0$.
This allows us to extract additional decay
of  solutions with the compactness property by means of the double Duhamel method. 
This method
was introduced in \cite{CKSTT08} and has been employed in many results  (see for instance \cite{KV10, KV10.2,  KV11.2, KV13, B12, DL15, CK23}). 
As a result of the double Duhamel method, combined with an interpolation argument, we prove in \cite{CK23} that
the trajectories of such solutions belong to the energy space.  In particular, 
we obtain the following vanishing: For all $R >0$,
\begin{align}
\limsup_{t \to -\infty} \|{ \vec{u} (t)}\|_{\dot{H}^1 \times L^2 (r \geq R + |t|)} = 
\limsup_{t \to \infty} \|{ \vec{u} (t)}\|_{\dot{H}^1 \times L^2 (r \geq R + |t|)} = 0.
\label{eq0.11intro}
\end{align}
In \cite[Section~4]{CK23}, the zero-limits in \eqref{eq0.11intro} are achieved by analyzing the size of the pairing 
\begin{equation}
 \left\langle P_M \int_0^{\infty} \frac{e^{i \tau \sqrt{-\Delta}}}{\sqrt{-\Delta}}  |u|^{p-1} u (\tau)~d\tau , 
   P_M \int^0_{-\infty} \frac{e^{i \tau \sqrt{-\Delta}}}{\sqrt{-\Delta}}  |u|^{p-1} u (\tau)~d\tau \right\rangle_{\dot{H}^{s_p}} 
   \label{eq0.14}
\end{equation}
in M, where $P_M$ denotes the Littlewood-Paley projection corresponding to the dyadic number $2^M$. This term may be split into 
four main components 
\begin{align*}
\begin{split}
\langle A+B, \tilde{A}+ \tilde{B} \rangle_{\dot{H}^{s_p}} = 
\langle A, \tilde{A}+ \tilde{B} \rangle_{\dot{H}^{s_p}} + 
\langle A+B,  \tilde{A} \rangle_{\dot{H}^{s_p}} 
- \langle A, \tilde{A} \rangle_{\dot{H}^{s_p}} 
+ \langle B, \tilde{B} \rangle_{\dot{H}^{s_p}} 
\end{split}
\end{align*} 
 where 
\begin{align*}
 \begin{split}
 A := 
  \int_0^{\Lambda/M} \frac{e^{i \tau \sqrt{-\Delta}}}{\sqrt{-\Delta}} P_M  |u|^{p-1} u (\tau)~d\tau
  +  \int_{\Lambda/M}^{\infty} \frac{e^{i \tau \sqrt{-\Delta}}}{\sqrt{-\Delta}}  (1- \chi) \left(4x/{ \tau} \right) P_M  |u|^{p-1} u (\tau)~d\tau 
 \end{split}
 \end{align*}
 and, with a smooth radial bump function $\chi$, we denote by
\begin{equation*}
B := \int_{\Lambda/M}^{\infty} \frac{e^{i \tau \sqrt{-\Delta}}}{\sqrt{-\Delta}}  \chi \left(4x/{ \tau} \right) P_M  |u|^{p-1} u (\tau)~d\tau \period 
\end{equation*} 
The terms $\tilde{A}$ and $\tilde{B}$ may be defined analogously in the negative time direction. 
The lower bound on the support set of $v_L$ in \eqref{eq0.10intro} helps us estimate the real part of
$\langle B, \tilde{B} \rangle_{\dot{H}^{s_p}}$ as follows: we express 
$\mbox{Re} (\langle B, \tilde{B} \rangle)$
 as the sum of the main term
   \begin{align}
\int_{\Lambda/M}^{\infty}
 \int_{-\infty}^{- \Lambda /M} 
\left\langle
    \chi \left(4x/{t} \right) F_1 (u(t)), 
  {\cos (( t - \tau)  \sqrt{-\Delta}) }  \chi \left(4 x/{\tau} \right)  F_2 ( u (\tau))  \right\rangle_{L^2} d\tau dt 
  \label{eq0.13}
\end{align}
and commutator terms. Due to the strong  Huygens principle as in \eqref{eq0.10intro}, it follows that
\begin{align*}
\mbox{supp} \left( \cos ( (t - \tau)  \sqrt{- \Delta}) \chi (4x/t) F_2 (u (\tau)) \right) 
\cap \mbox{supp} \left(   \chi \left(4x/{t} \right) F_1 (u (t)) \right) = \emptyset
\end{align*}
and we deduce that the main term \eqref{eq0.13} is indeed zero. 

For even dimensions, the strong Huygens principle does not provide a lower bound on the support of linear solutions. We have 
 \begin{equation}
 \mbox{ supp } ( v_L (t)) \subset \{ x:  |x- x_0| \leq t+a \}
 \label{eq0.12}
 \end{equation}
provided that $\supp ((v_0, v_1)) \subset B(x_0, a)$. 
Therefore, the inner product in an analogous term yields non-zero contribution. In Section~\ref{sec:Decay} of the present paper, 
we instead 
 investigate the size 
 of $\langle B, \tilde{B} \rangle$ directly by expanding 
 \begin{equation*}
     \cos ( (t - \tau)  \sqrt{- \Delta}) \chi (4x/t) F_2 ( u (\tau))
 \end{equation*}
  in the physical variable $x$ and analyzing the double integrability of the resulting inner product in time. In the current implementation of the Double Duhamel trick, we also make an adjustment to the splitting of  \eqref{eq0.14} to avoid the commutator terms we collected in \cite[Section~4]{CK23}. 

 In Section~\ref{sec: rigidity}, we recall some preliminary results needed for the rigidity method. For the remainder of the paper, we are concerned with the Cauchy problem
\begin{align}
   \begin{split}
    & \partial^2_{t} u
    - \partial^2_{r} u - \frac{3}{r} \partial_r u
    = |u|^{p-1} u
      \comma r > R + |t|,
    \\
   &  \vec{u} (0)
    = (u_0, u_1) \in \mathcal{H} (R)  \comma r> R
    \end{split}
   \label{equation4.1.intro}
  \end{align}
to study radial solutions outside a channel of width $R>0$. The local wellposedness theory of these solutions is covered in Subsection~\ref{subsec: Cauchy problem}. Additionally, we include the essential ingredients for the rigidity method: we recall (i)  the 
channels of energy estimates from \cite[Lemma~3.9]{DKMM22}; (ii) the construction  from \cite{CK23} of a family of smooth   solutions to 
\begin{align*}
- \partial_{rr} \varphi - \frac{3}{r} \partial_r \varphi = | \varphi |^{p-1} \varphi \comma r>0
\end{align*}
which fail to belong in the critical space $\dot{H}^{s_p} \times \dot{H}^{s_p -1} (\mathbb{R}^4)$. 
We close Section~\ref{sec: rigidity} by stating the main rigidity  theorem for our paper:

\begin{Proposition*}[Proposition~\ref{rigidity proposition}]
Suppose that $u(t,r)$ is a radial solution of \eqref{wave equations}
with $\vec{u} (t) \in C_b (\mathbb{R}, \dot{H}^{s_p} \times \dot{H}^{s_{p}-1})$. Suppose that for all $R_0 >0$, $u(t,r)$ is a radial solution of  
\eqref{equation4.1.intro}  for $r> R_0 + |t|$, with initial data in $\mathcal{H}(R_0)$ and $\norm{\vec{u} (t)}_{\mathcal{S}(R_0)} < \infty$. 
Assume that
       \begin{equation}
        \sum_{\pm} \lim_{t \to \pm \infty}     \int_{r > R_0 +|t|} |\nabla_{t,r} u (t,r)|^2  r^3 \, dr =0.
\label{eq5.16intro}
    \end{equation}
Then, there exists $R > 0$ such that 
    \begin{equation}
        (u_0 (r), u_1(r)) = (0,0) \comma \mbox{ for all }~ r > R. 
    \end{equation}
\end{Proposition*} 
 As in \cite{DKMM22}, the proof of the rigidity theorem consists of two key parts. Section~\ref{sec: rigidity section 1} and Section~\ref{sec: rigidity section 2} are devoted to covering Part~I and Part~II of the rigidity proof, respectively.

 In Part I, we study solutions as given in Proposition~\ref{rigidity proposition} that also have a constant sign on $\{r> R_0 + |t|\}$. For such a solution $u(t)$, we analyze the dominant even part
 \begin{equation*}
     u_+ (t) = \frac{u(t) + u(-t)}{2}. 
 \end{equation*}
 Observing that $u_+ (t)$ solves an inhomogeneous wave equation with data $(u_0, 0)$, we apply the exterior energy estimates in Lemma~\ref{Lemma3.9} and, with further analyses,  
 obtain asymptotic estimates for $u(t)$ on outer channels $\{r> R+ |t| \}$ with $R> R_0$. These estimates are given in Lemmas~\ref{Lemma4.8}--\ref{Lemma4.10}. Proceeding as in \cite[Proposition~4.6]{DKMM22}, we have the following  control
 \begin{equation*}
  \frac{M (\rho/2)}{2} \leq   M(\rho)  = \sup_{\substack{t_0 \in \mathbb{R} \\\sigma \geq \rho + |t_0|}}
      \sigma^2 |u (t_0, r) - Z_{\ell} (r)|  \leq \frac{C}{\rho^2}
 \end{equation*}
on the difference term $M (\rho)$ for large $\rho$, which leads to showing that there exists $R > R_0$ so that
\begin{equation*}
(u_0 (r) , u_1 (r)) = (Z_{\ell} (r), 0) \comma r > R
\end{equation*}
where $Z_{\ell}$ is a radial stationary solution to \eqref{equation4.1.intro} as constructed in Subsection~\ref{sec:stationary solutions}. At the end of Section~\ref{sec: rigidity section 1}, we prove that the above equality produces a contradiction by appealing to an argument in \cite[Proposition~4.7]{DKM20}. The main idea behind the contradiction stems from the fact that compactly supported nontrivial perturbations of a stationary solution  cannot obey the zero-limit condition in \eqref{eq5.16intro}.

In part~II, we focus on solutions $u(t)$ as in Proposition~\ref{rigidity proposition} whose data is of the form $(0, u_1)$.  The theory in Section~\ref{sec: rigidity section 1} allows us to reduce our main rigidity result to the following proposition:

\begin{Proposition*}[Proposition~\ref{Prop5.1}]
   Let $R>0$ and $\vec{u} (t)$ be a radial solution of \eqref{equation4.1} on $\{ r> R+ |t| \}$ with data $(0, u_1) \in \mathcal{H} (R)$. Assume that
   
   \begin{equation}
       \sum_{\pm} \lim_{t \to \pm \infty}  \int_{|x| > R + |t|} |\nabla_{t,x} u (t, x)|^2 \, dx =0.
       \label{eq7.1}
   \end{equation}
   Then, we have for all $t \in \mathbb{R}$ 
   \begin{equation}
     u(t,r) =0
    \comma r \geq R+ |t| .
       \label{eq7.2intro}
   \end{equation}
\end{Proposition*}

To prove Proposition~\ref{Prop5.1}, we first seek to achieve \textit{gain of regularity in time} as in \cite{DKMM22}. More precisely, we show that if $\vec{u} (t)$ is as in Proposition~\ref{Prop5.1}, then $(\partial_t u (t), \partial_t^2 u (t))$ restricted to $\{r> R+ |t| \}$ belongs to $C(\mathbb{R}, \dot{H}^1 \times L^2)$. 

Next, we utilize the function 
\begin{equation}
    a(t,r) = \frac{t}{r^2 (\log r)^{1/2}}
\end{equation}
as an odd in time approximate solution of \eqref{equation4.1.intro}. We note that the associated initial data $(a(0,r), \partial_t a(0,r))$ does not belong in the energy space, but the approximate solution $a(t,r)$ satisfies the same space-time exterior estimates as our solution $\vec{u} (t)$. We close the proof by the way of contradiction: if $u_1$ is not identically zero for $r> R$, then $u (t,r)$ stays asymptotically close to $a(t,r)$. We obtain the desired contradiction since $a(t,r)$ does not have finite energy.


\section{A review of the local theory}
\label{sec:A review of the Cauchy problem}

In this section, we concentrate on the theory of the local Cauchy problem for the nonlinear wave equation 
\begin{align}
   \begin{split}
    & \partial^2_{t} u
    - \Delta u 
    -  |u|^{p-1} u
     = 0   \comma \mbox{in } I \times \mathbb{R}^4 ,
    \\
   &  \vec{u} (0)
    = (u_0, u_1) \in \dot{H}^{s_p} \times \dot{H}^{s_p -1} (\mathbb{R}^4)
    \end{split}
   \label{nonlinear wave equations}
\end{align}
 where the exponent $p$ obeys the condition
 \begin{equation}
p >  4  \qquad \mbox{ and } \qquad  s_p = 2 - 2 / (p-1) \in (4/3, 2)
\period
\label{p-condition}
\end{equation}
The local well-posedness theory for the nonlinear wave equation 
\begin{equation*}
    \partial_t^2 - \Delta u = \pm |u|^{p-1} u
\end{equation*}
is established in \cite{KM11, DL15, CK23} in dimensions three, five, and greater than or equal to seven. To address the local theory of the Cauchy problem \eqref{nonlinear wave equations} in four dimensions, we proceed as in \cite{CK23}. 
We first briefly review the preliminaries, such as wave-admissible function spaces \eqref{eq2.2} and the Strichartz estimates \eqref{eq2.3}. Then we verify linear embeddings \eqref{embd1}--\eqref{embd4}, an inhomogeneous Strichartz estimate \eqref{Inhom SE}, and nonlinear estimates \eqref{EQ_NE1}--\eqref{EQ_NE4} leading to the small data theory and a perturbation result. 

We begin with the linear wave equation in $I \times \mathbb{R}^4$
\begin{align}
   \begin{split}
    & \partial^2_{t} \omega
    - \Delta \omega 
     = h   \comma 
    \\
   &  \vec{\omega} (0)
    = (\omega_0, \omega_1) \in \dot{H}^{s} \times \dot{H}^{s -1} (\mathbb{R}^4) .
    \end{split}
   \label{linear wave equations}
  \end{align}
A solution to \eqref{linear wave equations} with $h=0$ is called a \textit{free wave} and is given by
\begin{align}
 S(t) (\omega_0, \omega_1) = \cos (t \sqrt{- \Delta}) ~\omega_0 + (- \Delta)^{-1/2} \sin (t \sqrt{- \Delta})~\omega_1 .
 \label{eq1.3}
\end{align}
The right-hand side is defined via the Fourier transform:
we have
\begin{equation}
 \mathcal{F} (\cos (t \sqrt{- \Delta}) ~ f ) (\xi )  = \cos ( t |\xi|) \mathcal{F}( f ) (\xi) 
 \label{eq1.4}
\end{equation}
and
\begin{equation}
\mathcal{F} ( (- \Delta)^{-1/2} \sin (t \sqrt{- \Delta})~ f ) (\xi ) = | \xi |^{-1} \sin (t |\xi| ) \mathcal{F}( f ) (\xi).
\label{eq1.5}
\end{equation}
Similarly, the fractional differentiation operators are defined by
\begin{equation}
\mathcal{F} (D^{\alpha} f) (\xi) = |\xi |^{\alpha} \mathcal{F} (f) (\xi). 
\label{eq1.6}
\end{equation}
We recall the agreement  $D^2 = - \Delta$ under the Fourier transform.  For $s \in \mathbb{R}$, we may define the homogeneous Sobolev norm 
\begin{equation*}
    \norm{f}_{\dot{H}^{s,p}} : = \norm{D^{s} f}_{L^p (\mathbb{R}^4)}. 
\end{equation*}
By Duhamel's principle, the solution operator to \eqref{linear wave equations} is formally given by
\begin{align}
\omega (t) = S (t) (\omega_0, \omega_1) + \int_0^{t} \frac{\sin ((t-s) \sqrt{- \Delta})}{\sqrt{- \Delta}} h(s)~ds.
\label{eq1.2}
\end{align} 

Next, we recall the definition of Littlewood-Paley projections and homogeneous Besov spaces. Given a radial bump function $\varphi (\xi)$ supported in the ball 
$\{ \xi \in \mathbb{R}^d: |\xi| \leq 2\}$
with 
$
\varphi (\xi) = 1$ on $\{\xi \in \mathbb{R}^d: |\xi| \leq 1\}$, we define 
\begin{align}
& \mathcal{F} (P_{\leq N} f) (\xi)  = \varphi (\xi/N) \mathcal{F} (f) (\xi) \\
& \mathcal{F} (P_{> N} f) (\xi)  = (1- \varphi (\xi/N)) \mathcal{F} (f) (\xi) \\
& \mathcal{F} (P_N f) (\xi)  = (\varphi (\xi/N) - \varphi (2 \xi /N)) \mathcal{F} (f) (\xi). \label{eq1.7p}
\end{align} 
The frequency size $N$ will usually represent a dyadic number $N = 2^k$ for $k \in \mathbb{Z}$.

Let $s \in \mathbb{R}$ and $1\leq r \leq \infty$, 
$1 \leq q \leq \infty$. We denote by $\dot{B}^s_{r, q}$ the homogeneous Besov space 
\begin{equation*}
\dot{B}^s_{r, q} (\mathbb{R}^d ) = \{ u \in \mathcal{S}' (\mathbb{R}^d): \norm{u}_{\dot{B}^s_{r, q}} < \infty \} 
\end{equation*}
with the norm 
\begin{equation}
\norm{u}_{\dot{B}^s_{r, q}} = \left( \sum_{N \in 2^{\mathbb{Z}}} (N^s \norm{ P_N u}_{L^r (\mathbb{R}^d)} )^q   \right)^{\frac{1}{q}}
.
\label{EQ_Besov1}
\end{equation}

We say that a triple $(q, r, \gamma)$ is admissible if 
\begin{align*}
q, r \geq 2 \comma \frac{1}{q}  \leq \frac{3}{2} \left( \frac{1}{2} - \frac{1}{r} \right)   \comma \frac{1}{q}+ \frac{4}{r} = 2 - \gamma .
\end{align*}  
Given admissible triples $(q, r, \gamma)$ and $(a,b, \rho)$, Strichartz estimates \cite{T10, GV95}  yield that any solution $\omega (t)$ to the Cauchy problem \eqref{linear wave equations} satisfies
\begin{align}
\| \vec{\omega} (t)  \|_{L^{q}_t \left( I; \dot{B}^{s -\gamma}_{r,2} \times  \dot{B}^{s -\gamma -1}_{r,2} \right) }
 \les \| \vec{\omega} (0) \|_{\dot{H}^{s} \times \dot{H}^{s -1}}
 + \| h \|_{L^{a'}_t \left( I; \dot{B}^{s -1 + \rho}_{b',2} \right)}
 \period
 \label{eq2.2}
\end{align}
On the right-hand side above 
 $a'$ and $b'$ denote the conjugate indices to $a$ and $b$, respectively, i.e., we have
\begin{align*}
\frac{1}{a} + \frac{1}{a'} = 1 \comma \frac{1}{b} + \frac{1}{b'} =1. 
\end{align*}
Let $I\in \mathbb{R}$ be a time interval. 
As in \cite{CK23}, we follow the functional setting introduced in \cite[Section~2]{BCLPZ13}, and we define
\begin{equation*}
    \norm{u}_{\dot{S} (I)} = \sup \{ \norm{u}_{L^q_t \dot{B}^{s_p - \gamma}_{r,2} (I \times \mathbb{R}^4)}, \norm{ \partial_t u}_{L^q_t \dot{B}^{s_p -1 -\gamma}_{r,2} (
I\times \mathbb{R}^4)}: (q,r, \gamma)~ \mbox{is an admissible triple } \}.
\end{equation*}
We now define the function spaces $S_p (I)$, $W(I)$, and $X(I)$ associated to the respective admissible triples 
$(2 (p-1), \frac{8}{3} (p-1), s_p)$, $(2,8,1)$, and $(2(p-1), \frac{8 (p-1)}{(3p -4)}, \frac{1}{2})$. Namely, we define
\begin{align}
\begin{split}
& \norm{u}_{S_p(I)} = \norm{u}_{L^{2(p-1)}_t L^{\frac{8}{3} (p-1)}_x  (I \times \mathbb{R}^4 )} \\
& \norm{u}_{W(I)} = \norm{u}_{L^2_t \dot{B}^{s_p -1}_{8,2} (I\times \mathbb{R}^4) } \\
& \norm{u}_{X(I)} = \norm{u}_{L^{q_x}_t \dot{B}^{s_p - \frac{1}{2}}_{r_x,2} (I\times \mathbb{R}^4)} \comma q_x = 2(p-1), r_x = \frac{8 (p-1)}{3 p -4} 
\period
\end{split}
\label{norms}
\end{align}
By construction, as in Lemma~2.4 in \cite{CK23} (see also \cite[Section~2]{BCLPZ13}), the following embeddings
 hold in dimension four. We have 
\begin{align}
& \norm{u}_{X(I)} \les \norm{u}_{\dot{S}(I)}  \label{embd1} \\
& \norm{u}_{W(I)} \les \norm{u}_{\dot{S}(I)} \label{embd2} \\
& \norm{u}_{S_p (I)} \les \norm{u}_{X(I)}   \label{embd3} \\
& \norm{u}_{X(I)} \les \norm{u}^{\theta}_{S_p (I)} \norm{u}_{\dot{S}(I)}^{1- \theta} \quad \quad \mbox{ for some } \theta \in (0,1). \label{embd4}
\end{align} 
We also introduce the conjugate spaces for $W(I)$ and $X(I)$ by
\begin{align}
\begin{split}
    & \norm{u}_{W' (I)} = \norm{u}_{L^1_t \dot{H}^{s_p -1} (I\times \mathbb{R}^4) } \\
& \norm{u}_{X' (I)} = \norm{u}_{L^{\tilde{q}'}_t \dot{B}^{s_p -\frac{1}{2}}_{\tilde{r}', 2} (I\times \mathbb{R}^4) } \comma \tilde{q} = \frac{2(p-1)}{(p-2)}, \tilde{r}= \frac{8(p-1)}{2p -1} 
\period
\end{split}
\label{norms2}
\end{align}

As done in \cite{CK23}, we apply \eqref{eq2.2} with 
$(a, b, \rho) = (\infty, 2, 0)$ and any other admissible triple to obtain 
 the following Strichartz estimate
 for solutions to the linear Cauchy problem \eqref{linear wave equations}:
\begin{align}
\begin{split}
\sup_{t \in \mathbb{R}} \Vert  \vec{\omega} (t) \Vert_{\dot{H}^{s_p} \times \dot{H}^{s_p -1}}
+ \Vert \omega \Vert_{ \dot{S} (I)} 
  \les
\Vert \vec{\omega} (0) \Vert_{\dot{H}^{s_p} \times \dot{H}^{s_p -1}}
+
\Vert  h \Vert_{W' (I)}.
\end{split}
\label{eq2.3}
\end{align}

Next, we state an inhomogeneous Strichartz estimate. The lemma below is an application of \cite[Corollary~8.7]{T10}. An analogous version is given in \cite{CK23} for dimensions $d \geq 7$. 
\begin{Lemma}[Inhomogeneous Strichartz Estimates]
Let  $I \in \mathbb{R}$ denote any time interval. Then, we have
\begin{align}
\norm{ \int_0^t \frac{\sin ( (t- \tau) \sqrt{- \Delta}) }{\sqrt{- \Delta}} f(\tau) ~d \tau }_{X(I)} \les \norm{f}_{X'(I)} . \label{Inhom SE}
\end{align}
\end{Lemma}

\begin{proof}
We obtain \eqref{Inhom SE} by directly verifying the conditions for the $X(I)$-$X'(I)$ norms in \cite[Corollary~8.7 ]{T10}.
As outlined in \cite{CK23}, we begin with exponents 
$(q_x, r_1)$ 
and $(\tilde{q}, \tilde{r}_1)$, where
\begin{align}
(q_x, r_1) = (2 (p-1), {8}/{3})
\comma
(\tilde{q}, \tilde{r}_1 ) = (2 (p-1)/ (p-2) , 24/7).
\end{align}
We have  $r_1, \tilde{r}_1 \in (2, \infty)$ satisfying
the conditions
\begin{align*}
\frac{1}{q_x}  < 3 \left( \frac{1}{2} - \frac{1}{r_1}  \right) \comma
\frac{1}{\tilde{q}}  < 3 \left(  \frac{1}{2} - \frac{1}{\tilde{r}_1} \right)
\end{align*}
as well as 
\begin{equation*}
\frac{1}{q_x} + \frac{1}{\tilde{q}} = \frac{3}{2} \left(  1- \frac{1}{r_1} - \frac{1}{\tilde{r}_1}  \right) 
\end{equation*}
and 
\begin{equation*}
\frac{1}{r_1} \leq \frac{3}{\tilde{r}_1} \comma \frac{1}{\tilde{r}_1} \leq \frac{3}{r_1}. 
\end{equation*}
We also have $r_x \geq r_1$ and $\tilde{r} \geq \tilde{r}_1$ as $p>4$. 
Finally, we check the scaling condition:
\begin{equation*}
s_p - \frac{1}{2} + 4 \left( \frac{1}{2} - \frac{1}{r_x} \right) - \frac{1}{q_x} = 1- \left( - s_p + \frac{1}{2} + 4 \left( \frac{1}{2}- \frac{1}{\tilde{r}} \right) - \frac{1}{\tilde{q}} \right) . 
\end{equation*}
Having checked all the conditions for Corollary~$8.7$, we obtain \eqref{Inhom SE}.
\end{proof}

As in  \cite{CK23}, the functional setting introduced in \eqref{norms}-\eqref{norms2} allows us to verify the nonlinear estimates below. 
Let $F(u) = |u|^{p-1} u$, where $p$ obeys the condition  \eqref{p-condition}. The proof of Lemma~2.8 in \cite{CK23} guarantees that we have 
\begin{align}
\norm{ F(u)}_{W' (I)} \les \norm{u}^{p-1}_{S_p (I)} \norm{u}_{W (I)}
\label{EQ_NE1}
\end{align}
and 
\begin{align}
\begin{split}
\norm{F(u) - F(v)}_{W' (I)} & \les \norm{u-v}_{W(I)} \left( \norm{u}^{p-1}_{S_p (I)} + \norm{v}^{p-1}_{S_p (I)} \right) \\
& \indeq + \norm{u-v}_{S_p (I)} \left( \norm{u}^{p-2}_{S_p (I)} + \norm{v}^{p-2}_{S_p (I)} \right) \\
& \indeq \indeq \indeqtimes \left( \norm{u}_{W(I)} + \norm{v}_{W(I)} \right) \period
\end{split}
\label{EQ_NE2}
\end{align}
Additionally, we  observe the nonlinear estimates with the $X(I)$-$X' (I)$ pairing by \cite[Lemma~2.9]{CK23}. 
\begin{align}
\norm{ F(u)}_{X' (I)} \les \norm{u}^{p-1}_{S_p (I)} \norm{u}_{X (I)}
\label{EQ_NE3}
\end{align}
and 
\begin{align}
\begin{split}
\norm{F(u) - F(v)}_{X' (I)} & \les \norm{u-v}_{X(I)} \left( \norm{u}^{p-1}_{S_p (I)} + \norm{v}^{p-1}_{S_p (I)} \right) \\
& \indeq + \norm{u-v}_{S_p (I)} \left( \norm{u}^{p-2}_{S_p (I)} + \norm{v}^{p-2}_{S_p (I)} \right) \\
& \indeq \indeq \indeqtimes \left( \norm{u}_{X (I)} + \norm{v}_{X (I)} \right) \period
\end{split}
\label{EQ_NE4}
\end{align}

\subsection{Local Well-Posedness: small data theory, scattering, and a long time perturbation result}
\label{subsec:small data theory}
Below, we state the main local well-posedness results for the Cauchy problem \eqref{nonlinear wave equations}. The details behind the main results below follow closely those of \cite[Section~2]{CK23} and \cite[Section~2]{BCLPZ13}. 
We begin by recalling the definition of strong solutions.

\begin{definition}\label{def1}
Let $I \subset \mathbb{R}$ and $t_0 \in I$. 
We say that $u$ is a \textit{strong solution} of  the Cauchy problem \eqref{nonlinear wave equations} on $I$  if
$(u, \partial_t u) \in C(I; \dot{H}^{s_p} \times \dot{H}^{s_p -1})$, $u \in {S_p}(I) \cap W(I)$, $(u (t_0), \partial_t u (t_0))  = (u_0, u_1)$
and Duhamel's formula
\begin{equation*}
u (t) = S(t- t_0) ((u_0, u_1)) - \int_{t_0}^{t} \frac{\sin ((t-s)\sqrt{- \Delta})}{\sqrt{- \Delta}} F(u(s))~ds \comma t \in I
\end{equation*}
holds with $F(u)= |u|^{p-1} u$. 
\end{definition}
The local well-posedness result below is the four dimensional version of Theorem~2.12 from \cite{CK23}. Although Theorem~2.12 is stated for dimensions $d \geq 7$, the estimates \eqref{EQ_NE1}--\eqref{EQ_NE4} and Lemma~\ref{Inhom SE} guarantee that the arguments used in the proof of Theorem~2.12 stay valid for the Cauchy problem \eqref{nonlinear wave equations}. 
Therefore, we state the next result verbatim from \cite{CK23}.  
\begin{Theorem}
\label{Thm2.2}
Let $(u_0, u_1) \in \dot{H}^{s_p} \times \dot{H}^{s_p -1}$, and $p$ satisfy the condition \eqref{p-condition}.
Let $I \in \mathbb{R}$  be an interval containing $0 \in I^{\circ}$. 
 Assume that  
\begin{equation}
\Vert (u_0, u_1) \Vert_{ \dot{H}^{s_p} \times \dot{H}^{s_p -1}} \leq A.
\end{equation}
Then, there exists $\delta_0 = \delta_0 (p, A) > 0$ such that if 
\begin{align*}
\Vert S(t) (u_0, u_1) \Vert_{S_p (I)} = \delta < \delta_0,
\end{align*}
the Cauchy problem \eqref{nonlinear wave equations} admits a unique solution $u(t,x)$  on $I \times \mathbb{R}^4$ 
as in Definition~\ref{def1}. Moreover, we have
\begin{align*}
\Vert u \Vert_{S_p (I)} & < 2 \delta \\
\Vert u \Vert_{\dot{S} (I)} &< \infty \\
\Vert u \Vert_{X(I)}&  \leq C_0 \delta^{\theta} A^{1 - \theta}
\end{align*}
for some $\theta \in (0,1)$. 
Furthermore, if $(u_{0,k}, u_{1,k}) \to (u_0, u_1)$  as $k \to \infty$ in $\dot{H}^{s_p} \times \dot{H}^{s_p -1}$, 
then 
\begin{align*}
(u_k, \partial_t u_k) \to (u, \partial_t u) \quad \quad \mbox{in~} C(I; \dot{H}^{s_p} \times \dot{H}^{s_p -1}) 
\end{align*}
where $u_k$ is the solution corresponding to $(u_{0,k}, u_{1,k})$.
\end{Theorem}
 Below, we list some immediate consequences of Theorem~\ref{Thm2.2}. These results may be obtained by following the arguments based on Strichartz estimates in \cite{KM1, C03} in a standard manner.
 \begin{Remark}
\label{R2.3}
\textbf{Small data results.} 
\begin{itemize}
    \item[$1$] \textit{There exists $\tilde{\delta} >0$ such that if $\norm{(u_0, u_1)}_{\dot{H}^{s_p}\times \dot{H}^{s_p -1}} < \tilde{\delta}$, Theorem~\ref{Thm2.2} holds with $I = \mathbb{R}$.   }
    
    \item[$2$] \textit{ Let $(u_0, u_1) \in \dot{H}^{s_p} \times \dot{H}^{s_p -1}$.  There exists a maximal open interval $I$, where $0 \in I$ and} $$ I = I ((u_0, u_1)) = (T_-, T_+)$$ 
    \textit{
such that the solution $\vec{u} (t)$ is defined on $I \times \mathbb{R}^4$ as in Definition~\ref{def1}.}
\end{itemize}
\textbf{Standard blow-up criterion.}
\textit{Let $\vec{u} (t)$ be the solution of \eqref{nonlinear wave equations} on $(T_-, T_+)$. If $T_+ < \infty$, then we have 
$$\norm{u}_{S_p ([0, T_+))} = \infty .$$} 
\textbf{Scattering.}
\textit{If $T_+ = \infty$ and 
$\norm{u}_{S_p ([0, \infty)])} < \infty$, then $\vec{u} (t)$ scatters to a free wave as $t \to \infty$, i.e., there exists $(u_0^+, u_1^+) \in \dot{H}^{s_p} \times \dot{H}^{s_p -1}$ so that 
$$\lim_{t \to \infty} \norm{\vec{u} (t) - \vec{S} (t) (u_0^+, u_1^+)}_{\dot{H}^{s_p} \times \dot{H}^{s_p -1}} =0.$$}
\end{Remark}
Analogously, the scattering property and a finite time blow-up criterion may be stated for $T_-$ as well.

Next, we state a long time perturbation theorem in Sobolev spaces. This result applies to approximate solutions to \eqref{nonlinear wave equations}, and it will be needed in the concentration compactness method in the next section. 
\begin{Theorem}
\label{perturbation_thm_Sobolev}
Let  $(u_0, u_1) \in \dot{H}^{s_p} \times \dot{H}^{s_p -1}$ and $I \subset \mathbb{R}$ be an open interval containing $t_0$. 
Assume that $v$ solves the equation
\begin{align*}
\partial^2_{t} v
    - \Delta v 
     =    |v|^{p-1} v + e  
\end{align*} 
in the sense of the corresponding  integral equation, and it satisfies
\begin{equation}
\sup_{t \in I} \Vert v \Vert_{\dot{H}^{s_p} \times \dot{H}^{s_p -1}} 
+\Vert D^{s_p -1}   v \Vert_{L^2_t L^{8}_x (I \times \mathbb{R}^4)} \leq A
\label{EQ2.19}
\end{equation}
Additionally, assume that we have
\begin{align}
\Vert (u_0 - v(t_0), u_1 - \partial_t v (t_0) ) \Vert_{\dot{H}^{s_p} \times \dot{H}^{s_p -1}} \leq A'
\label{EQ2.4}
\end{align}
as well as the smallness condition 
\begin{align}
\Vert D^{s_p -1} e \Vert_{L^{1}_t L^{2}_x} + \Vert S(t- t_0) (u_0 - v(t_0), u_1 - \partial_t v (t_0) ) \Vert_{S_p (I)} 
\leq \epsilon .
\label{EQ2.5}
\end{align}
Then there exists $\epsilon_0 = \epsilon_0 (p,  A, A') >0$  such that if  $0 < \epsilon < \epsilon_0$,  there is a unique solution $u$ of \eqref{nonlinear wave equations}
 in $I$ with $(u(t_0), \partial_t u (t_0))= (u_0, u_1)$ 
satisfying
 \begin{align}
\Vert u-v  \Vert_{\dot{S}(I)} & \leq C(p, A, A')  , \label{EQ2.20}
\\
 \Vert u - v \Vert_{S_p (I)} & \leq  C(p, A, A') \cdot \epsilon^{\theta} \label{EQ2.21}
 \end{align}
 where $\theta \in (0,1)$. 
 \end{Theorem}
As discussed for the previous results in this section, the statement and the proof of Theorem~\ref{perturbation_thm_Sobolev} is the four dimensional analog of \cite[Theorem~2.17]{CK23}. We refer our readers to the detailed discussion in \cite[Section~2]{CK23}.
Therefore, we omit the proof of Theorem~\ref{perturbation_thm_Sobolev}. 

\begin{Remark}
\label{R2.19}
    A standard application of perturbation results such as Theorem~\ref{perturbation_thm_Sobolev} is that the flow associated to \eqref{wave equations} obeys the following continuity property: let $\vec{u} (t)$ be a solution of \eqref{wave equations} on $(T_- (\vec{u}), T_+ (\vec{u}))$ with initial data
$(u_0, u_1) \in \dot{H}^{s_p} \times \dot{H}^{s_p -1}$. 
Suppose that 
\begin{equation}
    (u_{0,n}, u_{1,n}) \to (u_0, u_1) \quad \mbox{in }~ \dot{H}^{s_p} \times \dot{H}^{s_p -1}.
\end{equation}
Denote by $\vec{u}_n (t)$ the solution to \eqref{wave equations} with data $(u_{0,n}, u_{1,n})$. Then, we have
\begin{equation}
    (T_- (\vec{u}), T_+ (\vec{u})) \subset 
    (\limsup_n T_- (\vec{u}_n), \liminf_n T_+ (\vec{u}_n)). 
\end{equation}
Furthermore, we have
\begin{equation}
    (u_{n} (t), \partial_t u_{n} (t)) \to (u (t), \partial_t u( t)) \quad  \mbox{in }~ \dot{H}^{s_p} \times \dot{H}^{s_p -1}.
\end{equation}
for each $t \in T_- (\vec{u}), T_+ (\vec{u}))$. See \cite[Remark~2.17]{KM1} for more details. 

\end{Remark}

\section{Concentration compactness procedure}
\label{sec:reduction}

In this section, we will
recall how to implement the method of the concentration-compactness following the approach in \cite{KM1, KM11, KM10, DKM14}. Using this method, we extract minimal blow-up solutions to \eqref{nonlinear wave equations} and study their additional properties, which 
essentially 
initiates the proof of Theorem~\ref{main thm_1} and reduce the proof 
to a rigidity statement given in Proposition~\ref{Proposition1}.   

\subsection{Critical solutions}

In analogy with \cite{KM1, KM10, KM11, DL15, CK23}, we start out with some definitions. 

Let $\vec{u} (t)$ be the unique solution to \eqref{nonlinear wave equations} with initial data $(u_0, u_1)$. Denote by $\Imax (\vec{u}) = (T_- (\vec{u}), T_+ (\vec{u}))$ the lifespan of the solution 
$\vec{u} (t)$. Throughout this section, we deal with radial initial data $(u_0, u_1) \in \dot{H}^{s_p} \times \dot{H}^{s_p -1}$.

For $A>0$, we define
\begin{align}
\mathcal{B}_{rad} (A) := \{(u_0, u_1) \in \dot{H}^{s_p} \times \dot{H}^{s_p -1} :~  (u_0, u_1)~ \mbox{ radial},~\sup_{t \in [0, T_+ (\vec{u})) } \Vert \vec{u}(t)\Vert_{\dot{H}^{s_p} \times \dot{H}^{s_p -1}} \leq A \}. 
\end{align}
We say that $\mathcal{SC}_{rad} (A)$ holds if for each $(u_0, u_1) \in \mathcal{B}_{rad} (A)$ we have $T_+ (\vec{u}) = \infty$ and $\Vert \vec{u} \Vert_{S_p ([0, \infty))} < \infty$. In addition, we say that $\mathcal{SC}_{rad} (A; (u_0, u_1))$ holds if $(u_0, u_1) \in \mathcal{B}_{rad} (A)$ implies  $T_+ (\vec{u}) = \infty$ and $\Vert \vec{u} \Vert_{S_p ([0, \infty))} < \infty$.

In this context, Theorem~\ref{main thm_1} is equivalent to the statement that $\mathcal{SC}_{rad} (A)$ holds for all $A>0$. Revisiting Theorem~\ref{Thm2.2}
and Remark~\ref{R2.3} in Section~\ref{sec:A review of the Cauchy problem} guarantees that $\mathcal{SC}_{rad} (\delta)$ holds for sufficiently small $\delta >0$. 
 Thus, under the assumption that Theorem~\ref{main thm_1} fails, there exists a critical value $A_C >0$ so that for $A < A_C$, $\mathcal{SC}_{rad} (A)$ holds, and for $A > A_C$, $\mathcal{SC}_{rad} (A)$ fails. 
The next result summarizes the main conclusions of the concentration-compactness procedure.

\begin{Proposition} 
\label{Proposition3}
Suppose that Theorem~\ref{main thm_1} is false. Then, there exists radial data $(u_{0,C}, u_{1,C})$ such that $\mathcal{SC}_{rad} (A_C, (u_{0,C}, u_{1,C}))$ fails. 
Additionally, there exists a continuous function $\lambda: [0, T_+ (\vec{u}_C)) \to (0, \infty)$  such  that
\begin{align}
    \inf_{t \in [0, T_+ (\vec{u}_C))} \lambda (t)  >0
    \label{eq.Lem20}
\end{align}
and the set 
\begin{equation}
\left\lbrace \left(\frac{1}{\lambda(t)^{\frac{2}{p-1}} } u_C \left( t, \frac{x}{\lambda (t)} \right), \frac{1}{\lambda(t)^{\frac{2}{p-1}+1 }} \partial_t u_C \left(t, \frac{x}{\lambda (t)}  \right) \right):~ t \in [0, T_+ (\vec{u}_C) )  \right\rbrace
\label{eq3.1}
\end{equation}
has compact closure in $\dot{H}^{s_p} \times \dot{H}^{s_p -1}$. 
\end{Proposition}

Solutions such as $\vec{u}_C (t)$ in Proposition~\ref{Proposition3} are called \textit{critical solutions}. More precisely, we say that
 $\vec{u}_C (t)$ is a critical solution if 
it satisfies  
\begin{equation*}
\sup_{t \in [0, T_+ (\vec{u}_C))} \norm{\vec{u}_C (t)}_{\dot{H}^{s_p} \times \dot{H}^{s_{p-1}}} = A_C \comma \norm{\vec{u}_C }_{S_p ([0, T_+ (\vec{u}_C ) )} = \infty,
\end{equation*}
and there exists a continuous function $\lambda: [0, T_+ (\vec{u}_C)) \to (0, \infty)$ so that the set given in \eqref{eq3.1} is pre-compact in $\dot{H}^{s_p} \times \dot{H}^{s_{p-1}}$.


Going back to Proposition~\ref{Proposition3}, the existence of a radial initial data $(u_{0, C}, u_{1,C})$ and a positive scaling function $\lambda (t)$ and pre-compactness of the set in \eqref{eq3.1} may be justified following the procedure in \cite[Propositions~3.3--3.4]{KM10}. In Section~3 of \cite{KM10}, authors develop a version of the concentration-compactness procedure that does not use any conservation laws. The theory of the  Cauchy problem for cubic NLS in three dimensions \cite[Section~2]{KM10}, combined with a profile decomposition theorem for linear solutions by Bahouri-Gerard \cite{BG99} yields a critical element, which enjoys a compactness property. 
 To carry out an analogous concentration-compactness procedure, we may rely on our  local well-posedness results from the previous section and 
a higher dimensional version of the profile decomposition theorem by Bulut \cite{B10} for initial data in $\dot{H}^{s_p} \times \dot{H}^{s_p-1}$. 
The continuity  property of
$\lambda (t)$ on $[0, T_+ (\vec{u}_C))$ follows from the continuity of $\vec{u}_C (t)$ on $[0, T_+ (\vec{u}_C))$ in $\dot{H}^{s_p} \times \dot{H}^{s_p -1}$. Details to obtain a critical solution with a continuous scaling function $\lambda (t)$ are outlined in Remark~5.4 of \cite{KM1}. Similarly, we may construct a critical element with a scaling function that satisfies \eqref{eq.Lem20}. We refer the readers to \cite[Lemma~3.10]{KM10} for details.

\subsection{The compactness property}
In view of the discussion in the previous subsection and Proposition~\ref{Proposition3}, we observe that Theorem~\ref{main thm_1} follows from 
the next result.
\cole
\begin{Theorem} 
\label{Theorem3.1}
Let $\vec{u} (t)$ be a radial solution of \eqref{wave equations}.
Assume that there exists a continuous function $\lambda: [0, T_+ (\vec{u}) ) \to (0, \infty)$ 
so that
\begin{equation}
K_+ := \left\lbrace \left(\frac{1}{\lambda(t)^{\frac{2}{p-1}} } u \left( t, \frac{x}{\lambda (t)} \right), \frac{1}{\lambda(t)^{\frac{2}{p-1}+1 }} \partial_t u \left(t,  \frac{x}{\lambda (t)} \right) \right):~t \in [0, T_+ (\vec{u}))  \right\rbrace
\end{equation}
has compact closure in $\dot{H}^{s_p} \times \dot{H}^{s_p -1}$ and we have
\begin{equation}
\inf_{t \in [0, T_+ (\vec{u}) )} \lambda (t) >0 . 
\end{equation}
Then, $\vec{u} \equiv (0,0)$.
\end{Theorem}
\colb

First, we introduce the following definition in order to study properties of solutions to \eqref{wave equations} as described in Theorem~\ref{Theorem3.1}.

\begin{definition}\label{compactness property}
Let $\vec{u} (t)$ be a solution of  \eqref{wave equations} defined on its maximal interval of existence $\Imax (\vec{u}) = (T_- (\vec{u}), T_+ (\vec{u}))$. 
We say that $\vec{u} (t)$ has the compactness property if there exists  $\lambda: \Imax (\vec{u}) \to (0, \infty)$ 
so that
 the set 
\begin{align}
K = \left\lbrace \left( \frac{1}{\lambda(t)^{\frac{2}{p-1}}} u \left( t, \frac{x}{\lambda (t)}  \right), \frac{1}{\lambda(t)^{\frac{2}{p-1}+1}} \partial_t u\left( t, \frac{x}{\lambda (t)}  \right) \right) 
: t \in \Imax (\vec{u}) \right\rbrace
\label{eq_comp_2}
\end{align}
has compact closure in ${\dot{H}^{s_p} \times \dot{H}^{s_p -1}}$. 
\end{definition}

We point out that Definition~\ref{compactness property} takes into account solutions whose scaling functions are defined on the entire lifespan $(T_-(\vec{u}), T_+ (\vec{u}))$ as opposed to $[0, T_+ (\vec{u}))$, which leads to a trajectory set different than $K_+$.  Nevertheless, we may still use solutions with compactness property to give a proof of Theorem~\ref{Theorem3.1}. The next lemma justifies this transition. 
\begin{Lemma} 
\label{Lem23}
Let $\vec{u} (t)$ be a solution of \eqref{wave equations} as in Theorem~\ref{Theorem3.1}. Let $\{ t_n \}_n$ be a sequence of times 
in $[0, T_+ (\vec{u}))$ such that $\lim_n t_n = T_+ (\vec{u})$. Assume that there exists $(v_0, v_1) \in \dot{H}^{s_p} \times \dot{H}^{s_p -1}$
such that 
\begin{equation}
\left( \frac{1}{\lambda (t_n)^{\frac{2}{p-1}}} u \left( t_n, \frac{x}{\lambda (t_n )} \right), \frac{1}{\lambda (t_n)^{\frac{2}{p-1} +1 }} \partial_t u \left( t_n,  \frac{x}{\lambda (t_n )} \right)    \right)
\to 
(v_0, v_1)
\quad 
\mbox{as~} n \to \infty
\label{eq3.2}
\end{equation}
in $\dot{H}^{s_p} \times \dot{H}^{s_p -1}$. Let $\vec{v} (t)$ be the solution of \eqref{wave equations} with initial data $(v_0, v_1)$ at $t=0$. 
Then, $\vec{v} \not\equiv (0,0)$  provided that $\vec{u}  \not\equiv (0,0)$. 
Additionally, $\vec{v} (t)$ has the compactness property. 
\end{Lemma}

Similar lemmas with detailed arguments may be found in \cite{KM11, DKM15.2, CK23}. See \cite[Lemma~3.11]{CK23} for an outlined proof with further references. 

Below, we include the main application of Lemma~\ref{Lem23}. Let $\vec{u} (t)$ be a solution of \eqref{wave equations} as in Theorem~\ref{Theorem3.1}. By Lemma~\ref{Lem23}, we obtain a solution $\vec{v} (t)$ that satisfies the compactness property. The next result shows that $\vec{v}(t)$ must be a global solution. Since the initial data  $(v_0, v_1)$  is obtained through the limiting process described in Lemma~\ref{Lem23}, we deduce that $\vec{u} (t)$ is global due to the continuity of the flow as stated in Remark~\ref{R2.19}.
In other words, the next proposition eliminates the case $T_+ (\vec{u}) <\infty$ for solutions $\vec{u} (t)$ as in Theorem~\ref{Theorem3.1}.

\cole
\begin{Proposition}[\text{\cite[Proposition~3.1]{DKM15}}]
\label{Prop3.5}
Let 
$\vec{u}(t)$ be a solution of \eqref{wave equations} with the compactness property. Then, $\vec{u} (t)$ is global. 
\end{Proposition}
\colb

The proof of Proposition~3.1 in \cite{DKM15} directly applies to our setting, proving Proposition~\ref{Prop3.5} above. An analogous result is also provided with a detailed outline of the original proof for the sake of completeness in \cite[Proposition~3.12]{CK23}.



We close this section with the following rigidity result, the proof of which will be covered in Section~\ref{sec: rigidity}. 
Note that the statement below incorporates the compactness property and the conclusion of Proposition~\ref{Prop3.5} into Theorem~\ref{Theorem3.1}. 

\cole
\begin{Proposition}
\label{Proposition1}
Let $\vec{u} (t)$ be a  radial solution of  \eqref{wave equations} with $\Imax (\vec{u}) = \mathbb{R}$, which has the compactness property on $\mathbb{R}$. Suppose that
we have
\begin{equation}
\inf_{t \in (- \infty, \infty) } \lambda (t) > 0. 
\label{eq3.3}
\end{equation}
Then,  $\vec{u}  \equiv (0,0)$. 
\end{Proposition}
\colb


Finally, we may proceed as in \cite[Lemma~6.3--6.6]{KM10} to show that Theorem~\ref{Theorem3.1} reduces to  Proposition~\ref{Proposition1}.




\section{Decay results for solutions with the compactness property}
\label{sec:Decay}

This section is devoted to study the decay properties of solutions satisfying the hypotheses of Proposition~\ref{Proposition1}. Namely, we will assume that $\vec{u} (t)$ is a radial solution of 
\eqref{wave equations} with 
$\Imax (\vec{u}) = \mathbb{R}$. Additionally, $\vec{u} (t)$ has  the compactness property on $\mathbb{R}$
and the corresponding scaling parameter $\lambda$  satisfies
\begin{equation}
\inf_{t \in \mathbb{R}} \lambda (t) >0.
\label{eq4.58}
\end{equation}
Below is the main result of this section showing that such solutions are uniformly bounded in time in a lower order homogeneous Sobolev space. 
\begin{Proposition}
\label{decay_proposition}
Let $\vec{u} (t)$ be a solution to \eqref{wave equations} as in Proposition~\ref{Proposition1}. 
Then, for all $t \in \mathbb{R}$ 
\begin{align}
\Vert\vec{u} (t) \Vert_{\dot{H}^{3/4} \times \dot{H}^{-1/4} (\mathbb{R}^d)} \leq C_p
\label{eq_decay1}
\end{align}
where the constant $C$ is uniform in time. 
\end{Proposition}
Combined with the compactness property in $\dot{H}^{s_p} \times \dot{H}^{s_p -1}$, estimate \eqref{eq_decay1} will yield pre-compactness in the energy space $\dot{H}^1 \times L^2$. This fact will allow us to apply the rigidity method in Sections~\ref{sec: rigidity section 1}--\ref{sec: rigidity section 2}.

The proof of Proposition~\ref{decay_proposition} is obtained via a double Duhamel technique. Analogous decay results were proven in \cite{DL15} for dimension five and in  \cite{CK23} for odd dimensions $d\geq 7$. In this application, the treatment is adjusted to compensate for the lack of strong Huygen's principle in even dimensions.

We begin by recalling a couple of preliminary results that will be applicable for solutions $\vec{u} (t)$ as in Proposition~\ref{decay_proposition}. 

\begin{Lemma}[Radial Sobolev inequality~\cite{TVZ07} ]
\label{Lem21}
Let $1 \leq p,q \leq \infty$, $0< s < 4$, and $\beta \in \mathbb{R}$ obey the conditions 
\begin{equation*}
\beta > - \frac{4}{q'} \comma 1 \leq \frac{1}{p} + \frac{1}{q} \leq 1+s 
\end{equation*}
and the scaling condition
\begin{equation*}
 4 - \beta -s =  \frac{4}{p'} + \frac{4}{q'}
\end{equation*}
with at most one of the equalities 
\begin{equation*}
p = 1, ~ p = \infty,~ q =1,~ q = \infty,~ \frac{1}{p} + \frac{1}{q} = 1+s
\end{equation*}
holding. Then, for any radial function $f \in \dot{W}^{s, p} (\mathbb{R}^4)$, we have
\begin{equation*}
\Vert |x|^{\beta} f   \Vert_{L^{q'} (\mathbb{R}^4)} \leq C \Vert  D^s f  \Vert_{L^{p} (\mathbb{R}^4)}
\period
\end{equation*}
\end{Lemma}
For solutions with the compactness property, the linear part of the evolution vanishes weakly in $\dot{H}^{s_p} \times \dot{H}^{s_p -1}$. Combined with the standard Duhamel formula, this fact leads to the following weak limit identities. Analogous results may be found in  \cite[Proposition~3.8]{S13} and \cite[Section~6]{TVZ08}.

\begin{Lemma}[Weak Limits]
\label{Lemma_weak_limits}
Let $\vec{u} (t)$ be a solution of the equation \eqref{wave equations} 
as in Proposition~\ref{Proposition1}. Then, for any $t_0 \in \mathbb{R}$ we have
\begin{align*}
u (t_0) &  =  - \lim_{T \to \infty } \int_{t_0}^{T} \frac{\sin ((t_0 - \tau)  \sqrt{- \Delta})}{\sqrt{- \Delta}} ~  |u|^{p-1} u ~d \tau
\qquad   \mbox{ weakly in }~\dot{H}^{s_p} (\mathbb{R}^4)
\\
u_{t} (t_0) & = - \lim_{T \to \infty} \int_{t_0}^{T} {\cos  ((t_0 - \tau)  \sqrt{- \Delta})} ~  |u|^{p-1} u ~d \tau
\qquad \mbox{ weakly in }~\dot{H}^{s_p -1} (\mathbb{R}^4)
\\ 
u (t_0) &  =   \lim_{T \to -\infty } \int_{T}^{t_0} \frac{\sin ((t_0 - \tau)  \sqrt{- \Delta})}{\sqrt{- \Delta}} ~  |u|^{p-1} u ~d \tau
\qquad   \mbox{ weakly in }~\dot{H}^{s_p} (\mathbb{R}^4)
\\
u_{t} (t_0) & =  \lim_{T \to -\infty } \int_{T}^{ t_0} {\cos  ((t_0 - \tau)  \sqrt{- \Delta})} ~ |u|^{p-1} u ~d \tau
\qquad \mbox{ weakly in }~\dot{H}^{s_p -1} (\mathbb{R}^4).
\end{align*}
\end{Lemma}

In a similar fashion, we may obtain the following uniform estimates for solutions with the compactness property. We refer the readers to \cite{DL15, CK23} for analogous results in odd dimensions. 

\begin{Lemma}[Uniformly Small Tails] 
\label{Lemma_tails}
Let $\vec{u} (t)$ be a solution of the equation \eqref{wave equations} as in Proposition~\ref{Proposition1}. 
Then for any $\eta >0$ there are  $0< c(\eta) < C(\eta) < \infty$ such that for all $t \in \mathbb{R}$
we have
\begin{align}
\begin{split}
& \int_{|x| \geq \frac{C(\eta)}{\lambda (t)}} \left|  D^{s_p} u (t,x)  \right|^2 dx
+
\int_{|\xi| \geq C(\eta) \lambda (t)} |\xi|^{2 s_p} |\hat{u} (t, \xi)|^2 d\xi 
 \leq \eta^2 \\
& \int_{|x| \leq \frac{c(\eta)}{\lambda (t)}} \left|  D^{s_p} u (t,x)  \right|^2 dx
+
\int_{|\xi| \leq c(\eta) \lambda (t)} |\xi|^{2 s_p} |\hat{u} ( t, \xi )|^2 d\xi 
 \leq \eta^2 \\
& \int_{|x| \geq \frac{C(\eta)}{\lambda (t)}} \left|  D^{s_p -1} u_t (t,x)  \right|^2 dx
+
\int_{|\xi| \geq C(\eta) \lambda (t)} |\xi|^{2 s_p -2} |\hat{u}_t (t, \xi)|^2 d\xi 
 \leq \eta^2 \\
& \int_{|x| \leq \frac{c(\eta)}{\lambda (t)}} \left|  D^{s_p -1} u_t (t, x)  \right|^2 dx
+
\int_{|\xi| \leq c(\eta) \lambda (t)} |\xi|^{2 s_p -2} |\hat{u}_t ( t, \xi )|^2 d\xi 
 \leq \eta^2 .
\end{split}
\label{eq4.uniform_tails}
\end{align}
\end{Lemma}

\begin{proof}[Proof of Proposition~\ref{decay_proposition}]
As in the previous applications of the double Duhamel technique in \cite[Section~4.2]{DL15} and \cite[Section~4]{CK23}, we define
\begin{align}
v = u + \frac{i}{\sqrt{- \Delta}} u_t 
\label{eq4.1}
\end{align}
and observe that
 \begin{align}
 \Vert \vec{u} (t) \Vert_{\dot{H}^{s_p} \times \dot{H}^{s_p -1}} \cong
 \Vert v(t) \Vert_{\dot{H}^{s_p}}
 \period
 \label{eq4.4}
 \end{align}
Also, $v(t)$ solves the equation
\begin{align}
v_t & = u_t + \frac{i}{\sqrt{- \Delta}} \left( \Delta u +  |u|^{p-1} u \right) \\
& = -i \sqrt{- \Delta} v + \frac{i}{\sqrt{- \Delta}}  |u|^{p-1} u \period
\label{eq4.2}
\end{align} 
By Duhamel's formula, we get
\begin{align}
v(t) = e^{-it \sqrt{- \Delta}} v(0) + i  \int_0^{t} \frac{e^{-i (t - \tau) \sqrt{- \Delta}}}{\sqrt{-\Delta}}
 |u|^{p-1} u (\tau)~d\tau
\period
\label{eq4.3}
\end{align} 
Moreover, by Lemma~\ref{Lemma_weak_limits}, we deduce that  for any $t_0 \in \mathbb{R}$
\begin{align}
\int_{t_0}^{T} \frac{e^{-i (t_0 - \tau) \sqrt{- \Delta}}}{\sqrt{-\Delta}}
u |u|^{p-1} (\tau)~d\tau
\rightharpoonup 
i v (t_0)
\comma \mbox{as } ~T \to \pm \infty
\label{eq4.3p}
\end{align}
weakly in $\dot{H}^{s_p}$.

Let $P_k$ denote the Littlewood-Paley projection corresponding to the dyadic number $2^k$. Recall the definition of $P_k$ in \eqref{eq1.7p}. Equivalently, we may express $P_k f$ as a convolution operator, that is 
in four dimensions,
 \begin{align}
 P_k f = 2^{4 k} \check{\phi} (2^k \cdot) \ast f
 \label{convolution_pk}
 \end{align}
 where $\check{\phi}$ belongs to the Schwartz class. 

 To prove \eqref{eq_decay1}, we proceed as in \cite{CK23} and look for a sequence $\beta = \{ \beta_k \}$ of positive numbers such that 
\begin{align}
\sup_{t \in \mathbb{R}} \Vert \left( P_k u (t), P_k u_t (t) \right) \Vert_{\dot{H}^{3/4} \times \dot{H}^{-1/4}} \les 2^{-\frac{3k}{4}} \beta_k
 \label{eq4.5}
\end{align}
 for all $k \in \mathbb{Z}$, and
 \begin{align}
 \Vert \{ 2^{-\frac{3k}{4}} \beta_k \}_k\Vert_{\ell^{2}} \les 1 \period
 \label{eq4.6}
 \end{align}

  \begin{Claim}
 \label{claim_1}
 A frequency envelope $\{\beta_k \}_k$ of positive numbers obeying \eqref{eq4.5}--\eqref{eq4.6} may be defined as below: we set
 \begin{align}
 \beta_k := 1 \inon{ for }~ k \geq 0
 \label{eq4.7}
 \end{align}
 and for $k <0$, we denote by
 \begin{align}
 \beta_k := \sum_{j} 2^{- |j-k|} a_j
 \label{eq4.8}
 \end{align}
 where 
 \begin{align}
 a_j = 2^{s_p j} \Vert P_j u \Vert_{L_t^{\infty} L^2} + 2^{( s_{p} -1) j} \Vert P_j u_t \Vert_{L_t^{\infty} L^2} 
 \inon{ for }~ j \in \mathbb{Z}.
 \label{eq4.9}
 \end{align}
 \end{Claim}
 
We immediately observe that the choice $\beta_k =1$ for $k \geq 0$ in \eqref{eq4.7} is sufficient due to the uniform in time control of the $\dot{H}^{s_p} \times \dot{H}^{s_p -1}$ norm of $\vec{u} (t)$ and the fact that $s_p \geq 3/2$. Toward the end of the proof, we will observe that we may simply use the estimate
\begin{equation*}
    a_j \cong \norm{P_j v}_{L_t^{\infty} \dot{H}^{s_p}} 
    \les 1
\end{equation*}
when $j >0$. Therefore, we will focus on estimating $a_j$ for $j \in \mathbb{Z}^{-}$ for the remainder of the proof. 

Let $j \in \mathbb{Z^-}$. Recalling \eqref{eq4.1} and Bernstein's inequalities, we observe that
\begin{align}
a_j \cong
 \Vert P_j  v \Vert_{L_t^{\infty} \dot{H}^{s_p}} \cong 
  2^{s_p j} \Vert P_j  v \Vert_{L_t^{\infty} L^2}
 \period
 \label{eq4.10}
 \end{align}

We will first estimate the term $\Vert P_j  v (0) \Vert_{\dot{H}^{s_p}}$. Throughout this case, the estimates will be uniform in $t$, yielding an estimate for \eqref{eq4.10}.


For $j \in \mathbb{Z}^-$, we let $M= 2^j$ denote an arbitrary dyadic number in $(0,1)$. 
Proceeding as in \cite[Section~4]{CK23}, by \eqref{eq4.3}--\eqref{eq4.3p}, we obtain the reduction
 \begin{align}
 \begin{split}
&  \left\langle P_M v(0), P_M v(0) \right\rangle_{\dot{H}^{s_p}}  \\
 & \qquad = 
   \left\langle P_M \int_0^{\infty} \frac{e^{i \tau \sqrt{-\Delta}}}{\sqrt{-\Delta}}  |u|^{p-1} u (\tau)~d\tau , 
   P_M \int^0_{-\infty} \frac{e^{i \tau \sqrt{-\Delta}}}{\sqrt{-\Delta}}  |u|^{p-1} u (\tau)~d\tau \right\rangle_{\dot{H}^{s_p}} \\
 & \qquad = 
  \left\langle \int_0^{\infty} e^{i \tau \sqrt{-\Delta}}  D^{s_p -1} P_M  (|u|^{p-1} u (\tau))~d\tau , 
   \int^0_{-\infty} e^{i \tau \sqrt{-\Delta}} D^{s_p -1}  P_M (|u|^{p-1} u (\tau))~d\tau \right\rangle_{L^2}
   \period
   \end{split}
   \label{eq4.12}
\end{align}

Let $\chi \in C_0^{\infty} (\mathbb{R}^4)$ be a radial non-increasing bump function, which satisfies
 \begin{align*}
\chi (x) = 
\begin{cases}
1 & \inon{if  }~  |x| \leq 1, \\
0 & \inon{if  }~  |x| \geq 2 \period
\end{cases}
\end{align*}
We then split the $L^2$ inner product in
  \eqref{eq4.12} as
\begin{align}
\begin{split}
\langle A+B, \tilde{A}+ \tilde{B} \rangle_{L^2} = 
\langle A, \tilde{A}+ \tilde{B} \rangle_{L^2} + 
\langle A+B,  \tilde{A} \rangle_{L^2}
- \langle A, \tilde{A} \rangle_{L^2} 
+ \langle B, \tilde{B} \rangle_{L^2} 
\end{split}
\label{EQ4.2}
\end{align} 
 where 
\begin{align}
 \begin{split}
 A &:= 
  \int_0^{\Lambda/M} {e^{i \tau \sqrt{-\Delta}}}  
  D^{s_p -1}
  P_M ( |u|^{p-1} u (\tau))~d\tau
  +  
  \int_{\Lambda/M}^{\infty} {e^{i \tau \sqrt{-\Delta}}} (1- \chi) \left(4x/{\tau} \right)
  D^{s_p -1}
  P_M ( |u|^{p-1} u (\tau))~d\tau \\
  & = A_1 + A_2
 \end{split}
 \label{eq4.13}
 \end{align}
 and 
\begin{equation}
B := \int_{\Lambda/M}^{\infty} 
e^{i \tau \sqrt{-\Delta}}
\chi \left(4 x/{\tau} \right) 
D^{s_p -1}
P_M ( |u|^{p-1} u (\tau))~d\tau \period \label{eq4.14} 
\end{equation} 
The terms $\tilde{A}$ and $\tilde{B}$ are defined analogously in the negative time direction. 
 The constant $\Lambda>0$ will be determined later.  

We note that the splitting in \eqref{EQ4.2} is slightly different than the one in \cite{CK23}. The change in the order of $\chi (4x / \tau)$ and $D^{s_p -1}$ will help us in the treatment of $\langle B, \tilde{B} \rangle$ in dimension four, where we do not have the strong Huygen's principle. 

The readers may also observe a minor change in estimates \eqref{eq4.17}--\eqref{eq4.19}. In \cite[Section~4]{CK23}, analogous decompositions were done around the frequency $M/2p$ due to a  mistake regarding the size of the support for the Littlewood-Paley multiplier associated to $P_M$. Below, we correct this miscalculation by selecting $M/8p$ for the high-low frequency decomposition cut-off, i.e., we use
\begin{equation*}
    u = P_{\leq M/8p} u + P_{> M/8p} u . 
\end{equation*}
See \eqref{eq4.18.2}--\eqref{eq4.18.3} for more details on the selection of $M/8p$. 

Now, we turn to estimating  $\norm{P_j v(0)}_{\dot{H}^{s_p}}$ as in \cite[Claim~4.7]{CK23}. We begin with the following claim. 
 \begin{Claim}
\label{claim_2}
Let $\eta >0$ be arbitrary. There is $N_0 >0 $ so that
\begin{equation}
\Vert A_1 \Vert_{L^2} \les \Lambda M^{s_p} \eta^{p-1}  \Vert P_{> M/8p} u \Vert_{L_{t}^{\infty} L^{2}} + \Lambda M^{s_p} N_0^{-s_p} \period
\label{eq4.28}
\end{equation}
\end{Claim} 
By the Plancherel and Bernstein's inequalities, we first observe that
\begin{align}
\begin{split}
\norm{A_1}_{L^2}
&  \les 
\frac{\Lambda}{M} \sup_{\tau \in \left[ 0, {\Lambda}/{M} \right]}
   \norm{  {e^{i \tau \sqrt{-\Delta}}} D^{s_p -1} P_M (  |u|^{p-1} u ) (\tau)}_{L^2} \\
& \qquad \les \Lambda M^{s_p -2} \norm{ P_M  (u |u|^{p-1}) }_{L_t^{\infty} L^2} \period  
\end{split}
\label{eq4.15}
\end{align}
Let $\eta >0$ be a small number to be determined later.
As in \cite{CK23}, Lemma~\ref{Lemma_tails} guarantees that there is a positive dyadic number $N_0 (\eta)$ such that for every $0 \leq K \leq N_0$ we have 
\begin{equation}
\Vert P_{\leq K} u \Vert_{\dot{H}^{s_p}}^2  \leq \sum_{N \leq N_0} N^{2 s_p} \Vert P_N u \Vert_{L^2}^2 \les   \eta^2
\label{EQ4.15p}
\end{equation}
which then leads to
\begin{equation}
\Vert P_{\leq K} u \Vert_{L^{2 (p-1)}} \les \eta
\label{eq4.16}
\end{equation}
 by the Sobolev embedding. 

 Analogously, we may
 estimate the term $\Vert P_M (u |u|^{p-1}) \Vert_{L_t^{\infty} L^{2}}$ in \eqref{eq4.15}
 as done in \cite{CK23}. We decompose
  \begin{align}
  \begin{split}
  \Vert P_M (u |u|^{p-1}) \Vert_{L^{2}}& = \Vert P_M ( u |P_{\leq N_0} u|^{p-1} + (u |u|^{p-1} - u |P_{\leq N_0} u |^{p-1}) ) \Vert_{L^2}
   \\  & \leq \Vert  P_M (u | P_{\leq N_0} u |^{p-1}) \Vert_{L^2}  +   \Vert P_M (u |u|^{p-1} - u |P_{\leq N_0} u |^{p-1})  \Vert_{L^2}
   \\ & = I +\tilde{I}.
	\end{split}
	\label{eq4.17p} 
  \end{align}
 We then  write  $I$
 as a product of two factors decomposed into high-low frequencies around ${M}/{8p}$ and $N_0$. 
 In other words, we get
 \begin{align}
 \begin{split}
I &
= \Vert P_M (( P_{\leq M/8p} u + P_{>M/8p} u ) | P_{\leq N_0} u |^{p-1} ) \Vert_{L^{2}}
\\ & \leq \Vert P_M ( P_{\leq M/8p} u  | P_{\leq N_0} u |^{p-1} ) \Vert_{L^{2}} 
+ \Vert P_M (  P_{> M/8p} u  | P_{\leq N_0} u |^{p-1} )  \Vert_{L^{2}}
\\ &  = I_1 + I_2.
\end{split}
\label{eq4.17}
\end{align} 

We begin with
 \begin{align}
 I_1 = \Vert P_M (( P_{\leq M/8p} u )  | P_{\leq N_0} u |^{p-1} ) \Vert_{L^{2}} \period
 \label{eq4.18}
 \end{align}
 Recalling \eqref{eq1.7p} and the fact that $p-1$ is an even integer, we deduce that 
 \begin{equation}
 \mbox{ supp } \left( \mathcal{F} \left( ( P_{\leq N_0 } u )^{p-1} \right)  \right) \subset \{ \xi: ~|\xi| \leq 2 (p-1) N_0\} .
 \label{eq4.18.2}
 \end{equation}
Note that 
\begin{equation}
    P_M \left( (P_{K} u)^p \right) =0
    \label{eq4.18.3}
\end{equation}
 if $K < M/4p$. 
 Therefore, in the case $N_0 \leq M/8p$, we get $I_1 =0$. 
 
 Without loss of generality, we assume that 
 $N_0 > M/8p$, and split $I_1$ in two parts
 \begin{align}
 \begin{split}
 & I_1 \les 
 \Vert 
 P_M (( P_{\leq M/8p} u )  | P_{ \leq {M}/{8p}} u |^{p-1} )
 \Vert_{L^2}
 \\ & \qquad + 
 \Vert
 P_M ( ( P_{\leq M/8p} u ) (   | P_{ \leq N_0} u |^{p-1} - | P_{\leq M/8p} u |^{p-1} ))
 \Vert_{L^2}
 \\ & \qquad = I_{11} + I_{12} \period
 \end{split}
 \label{eq4.19}
 \end{align}
 As we just discussed, the term $I_{11} =0$. 
 For $I_{12}$, using the fundamental theorem of calculus, we express the difference term  by 
 \begin{align}
\left|  | P_{ \leq N_0} u (x) |^{p-1} - | P_{\leq M/8p} u (x) |^{p-1} \right|
  \leq  \left|	P_{{M}/{8p} < \cdot < N_0} u (x) \right| 
  \int_0^1   \left| \lambda P_{\leq N_0} u (x) +  (1- \lambda) P_{\leq M/8p} u (x) \right|^{p-2} d\lambda
 .
 \end{align}
 Next, 
 we recall \eqref{convolution_pk}, and  apply Young's and H\"older's inequality
 \begin{align}
 \begin{split}
I_{12}  & \les 
 \Vert    M^4   \check{\phi} (M \cdot)    \Vert_{L^{2}}
 ~
 \Vert \left|  P_{\leq M/8p} u \right| ( | P_{\leq N_0} u |^{p-2}+ | P_{\leq M/8p} u|^{p-2} ) \Vert_{L^{2}}
 ~
  \Vert     P_{{M}/{8p} < \cdot < N_0} u    \Vert_{L^2}
\\ 
 & \les
 M^2 \eta^{p-1}  \Vert     P_{> {M}/{8p}} u    \Vert_{L^2} 
 \end{split}
 \label{eq4.20}
 \end{align}
 where we used \eqref{eq4.16} on the last line. 
By \eqref{eq4.19}, we then obtain
\begin{equation}
I_{1} \les M^2 \eta^{p-1} \Vert P_{> M/8p} u \Vert_{L^2}.
 \label{EQ4.52}
\end{equation}

 Similarly, we estimate
  \begin{align}
 \begin{split}
 I_2  & \les \Vert    M^4   \check{\phi} (M \cdot)    \Vert_{L^{2}}
 ~ 
 \Vert     P_{> M/8p} u    \Vert_{L^2}
 ~\Vert | P_{\leq  N_0}    u  |^{p-1}   \Vert_{L^{2}}
  \\
 & \les
 M^2 \eta^{p-1}  \Vert     P_{> {M}/{8 p}} u    \Vert_{L^2} .
 \end{split}
 \label{eq4.21}
 \end{align}
 Combining \eqref{EQ4.52} and \eqref{eq4.21}, we obtain
 \begin{align}
 I \les  M^2 \eta^{p-1}  \Vert     P_{> {M}/{8p}} u    \Vert_{L^2} 
 \period
 \label{eq4.21p}
 \end{align}
 
Going back to \eqref{eq4.17p}, we next estimate  $\tilde{I}$. In a similar manner, we control the difference as above:
 \begin{align*}
\left|  | u |^{p-1} - | P_{\leq N_0 } u |^{p-1} \right|
\leq 	\left| P_{>  N_0} u \right|
  \int_0^1 |  \lambda u + (1- \lambda )  P_{\leq N_0} u |^{p-2}  d\lambda
  .
 \end{align*}
 Next, applying Young's and H\"older's inequality followed by the Sobolev embedding and Bernstein's inequality 
 at the last step, we obtain
 \begin{align}
 \begin{split}
 \tilde{I}  & \les 
 \Vert    M^4 \check{\phi} (M \cdot)    \Vert_{L^{2}}
 ~
 \Vert u~ ( | u|^{p-2} +  |P_{\leq N_0} u|^{p-2} ) \Vert_{L^{2}}
 ~
  \Vert     P_{ > N_0} u    \Vert_{L^2}
\\ 
 & \les M^2 
 \Vert u \Vert_{L^{2 (p-1)}}^{p-1} 
 ~ 
 \Vert   P_{>  N_0} u  \Vert_{L^2}
 \\ 
 & \les M^2 
 N_0^{- s_p} 
 ~ 
 \Vert u \Vert_{L^{\infty}_t \dot{H}^{s_p}}^{p} . 
 \end{split}
 \label{eq4.22}
 \end{align}

 Then, we collect the estimates \eqref{eq4.17p}, \eqref{eq4.21p}, and \eqref{eq4.22}, and we obtain
 \begin{align}
 \| P_M (u |u|^{p-1}) \|_{L^2} \les M^2 \eta^{p-1} \| P_{> M/8p} u \|_{L^2} + 
 M^2 N_0^{- s_p} .
 \label{eq4.22p}
 \end{align}
 We note that the implicit constant included via the inequality sign $\les$ may depend on the norm $\norm{v}_{L^{\infty}_t \dot{H}^{s_p}}$.
By \eqref{eq4.15}, we arrive at the estimate
 \begin{align}
\begin{split}
 \norm{A_1}_{L^2} 
& \les \Lambda M^{s_p -2} \Vert P_M  |u|^{p-1} u \Vert_{L_t^{\infty} L^2} \\
&  \les 
\Lambda M^{s_p} \eta^{p-1} \Vert P_{> M/8p} u \Vert_{L_t^{\infty} L^{2}} 
+ \Lambda M^{s_p} N_0^{- s_p}
\period  
\end{split}
\label{eq4.27}
\end{align}
 This completes the proof of Claim~\ref{claim_2}. 
Going back to \eqref{eq4.13}, we next show the following estimate on $\norm{A_2}_{L^2}$.
\begin{Claim}
\label{claim_3}
\begin{align}
   \norm{A_2}_{L^2}
 \les
 \frac{M^{s_p}}{\Lambda^{1/4}} \norm{ P_{> M/8p} u}_{L^{\infty}_t  L^2 } . 
 \label{EQ4.45p}
 \end{align}
\end{Claim}

Revisiting \eqref{eq4.13}, we get
\begin{align}
\norm{A_2}_{L^2} \les 
\int_{\Lambda/M}^{\infty} \norm{ (1 - \chi) (4 x / \tau)
    P_M (|u|^{p-1} u (\tau))
    }_{L^2} \, d \tau
    \period
    \label{eq4.29}
\end{align}
We recall that 
\begin{equation*}
    \mbox{ supp} (1 - \chi) (4x/ \tau) 
    \subset 
    \left\lbrace x: ~ |x| \geq \frac{\tau}{4} \right\rbrace
    .
\end{equation*}
This observation, combined with Lemma~\ref{Lem21} for radial functions
\begin{equation*}
    \norm{|x|^{5/4}  f}_{L^2} \les \norm{D^{3/4} f}_{L^1}
\end{equation*}
 leads to an upper bound with integrability in $\tau$:
\begin{align}
\begin{split}
    \norm{(1 - \chi) (4 x / \tau)
    P_M (|u|^{p-1} u (\tau))}_{L^2}
    &\les 
    \frac{1}{|\tau|^{5/4}} \norm{|x|^{5/4} D^{s_p -1} 
    P_M (|u|^{p-1} u (\tau))    }_{L^2}
    \\ & \les
    \frac{1}{|\tau|^{5/4}} \norm{ D^{s_p - 1/4} 
    P_M (|u|^{p-1} u (\tau))    }_{L^1}    
    . 
    \end{split}
    \label{eq4.30}
\end{align}
Furthermore, using Bernstein's inequalities on the upper bound above, we obtain
\begin{align}
\begin{split}
    \norm{(1 - \chi) (4 x / \tau)
    P_M (|u|^{p-1} u (\tau))}_{L^2}
     \les
    \frac{M^{s_p - 1/4}}{|\tau|^{5/4}} \norm{
    P_M (|u|^{p-1} u (\tau))    }_{L^1}    
    \end{split}
    \label{eq4.31}
\end{align}
which implies that
\begin{equation*}
    \norm{A_2}_{L^2} \les
    \frac{M^{s_p}}{\Lambda^{1/4}}
    \norm{
    P_M (|u|^{p-1} u (\tau))    }_{L^{\infty}_t L^1}   . 
\end{equation*}
To close the estimate on $\norm{A_2}_{L^2}$, we must estimate $\norm{P_M (|u|^{p-1} u (\tau))}_{L^{\infty}_t L^1}$. As done in the proof of Claim~\ref{claim_2}, we
decompose $u$ around the frequency $M/8$, which leads to
\begin{align*}
\begin{split}
\norm{P_M ( |u|^{p-1} u)}_{L^{1}}
& \leq
 \norm{P_M (( P_{\leq M/8p} u + P_{> M/8p} u ) |u|^{p-1} ) }_{L^1}
\\  & \leq
 \norm{P_M( ( P_{\leq M/8p} u ) |P_{\leq M/8p} u|^{p-1} )}_{L^1 }
 \\ & \indeq \indeq  + 
\norm{P_M( ( P_{\leq M/8p} u ) ( |u|^{p-1} - |P_{\leq M/8p} u|^{p-1} ) )}_{L^{1} }
\\ & \indeq \indeq 
+ \norm{P_M ( (  P_{> M/8p} u ) |u|^{p-1} ) }_{L^{1} }.
\end{split}
\end{align*}
As in \eqref{eq4.18}--\eqref{eq4.19}, we deduce that the first term $ \norm{P_M \left( P_{\leq M/8p} u \right) |P_{\leq M/8p} u|^{p-1} }_{L^{1}}=0$. 
To treat the last two terms we simply apply Young's inequality and the Sobolev embedding, leading to the estimate
 \begin{align}
 \begin{split}
 \norm{P_M  \left(  |u|^{p-1} u (t)  \right) }_{L^{1} } 
&  \les \norm{M^4  \check{\phi} (M \cdot) }_{L^1} 
 \norm{ P_{> M/8p} u}_{L^2}
 \norm{u}^{p-1}_{L^{2(p-1)}} 
 \\ 
 & \les \norm{P_{> M/8p} u}_{L^2} 
 \end{split}
 \label{EQ4.44p}
 \end{align}
 as done in \eqref{eq4.22}. 
 As a result, we obtain the estimate
 \begin{align}
 \norm{A_2}_{L^2}
 \les
 \frac{M^{s_p}}{\Lambda^{1/4}} \norm{ P_{> M/8p} u}_{L^{\infty}_t  L^2 } . 
 \end{align}
 
  Next, we combine the estimates for $A_1$ and $A_2$ from Claim~\ref{claim_2} and Claim~\ref{claim_3}, and we set $\Lambda:= \eta^{-2}$ for $\eta \in (0,1)$ to be determined later. The estimate for $\norm{A}_{L^2}$ 
  adds up to the following upper bound:
 \begin{align}
 \begin{split}
 \norm{A}_{L^2 } & \les \norm{A_1}_{L^2} + \norm{A_2}_{L^2}
 \\ &  \les 
\eta^{p-3} M^{s_p} \norm{ P_{>M/8p} u}_{L^{\infty}_t L^2}
+ \eta^{-2} M^{s_p} N_0^{-s_p}
+ \eta^{1/2} M^{s_p} \norm{P_{>M/8p} u}_{L^{\infty}_t L^2}
\\ & \les 
\eta^{-2} M^{s_p} N_0^{-s_p}
+ \eta^{1/2} M^{s_p} \norm{P_{>M/8p} u}_{L^{\infty}_t L^2}
 \end{split}
 \label{EQ4.46p}
 \end{align}
 where we relied on the fact that $p-3 > 1/2$ at the last step. 
 Note that the upper bound above also controls the term $ \Vert {\tilde{A}} \Vert_{L^2}$.  

 Estimate \eqref{EQ4.46p} may also be utilized to estimate $\langle A, \tilde{A} + \tilde{B}
 \rangle_{L^2 }$ and $\langle A+B, \tilde{A} \rangle_{L^2}$. By Lemma~\ref{Lemma_weak_limits}, we have
  $e^{it \sqrt{ - \Delta} } v(t) \rightharpoonup 0$ weakly in 
 $\dot{H}^{s_p} $ as $t \to \pm \infty$. As in \cite[Section~4]{CK23}, we observe that
 \begin{align*}
     \langle A, \tilde{A} + \tilde{B}
 \rangle_{L^2 } 
 & = \left \langle D^{s_p} D^{-s_p} A, D^{s_p} P_M \int_{- \infty}^0 \frac{e^{it \sqrt{- \Delta}}}{\sqrt{- \Delta}} |u|^{p-1} u (t) \, dt  \right \rangle_{L^2}
 \\ & = 
 \langle
 D^{s_p} D^{-s_p} A, D^{s_p}  P_M v(0)
 \rangle_{L^2}. 
 \end{align*}   
By Cauchy-Schwarz, we may control $ \langle A, \tilde{A} + \tilde{B}\rangle_{L^2}$ as below: 
 \begin{align}
 \begin{split}
& 
 | \langle A, \tilde{A}+ \tilde{B} \rangle_{L^2 } |
 \les 
   \norm{A}_{L^2}
   \norm{P_M v }_{L_t^{\infty} \dot{H}^{s_p}}
 \\ & \quad \quad \les 
 \left(
(\eta^{p-3} + \eta^{1/2} ) M^{s_p} \norm{ P_{>M/8p} u}_{L^{\infty}_t L^2} + \eta^{-2} M^{s_p} N_0^{-s_p}
\right)
  M^{s_p}  \norm{P_M v }_{L_t^{\infty} L^2} .
 \end{split}
 \label{eq4.40}
\end{align}  
We note that the same upper bound also controls the term $ | \langle A+B , \tilde{A}  \rangle_{L^2 } |$. 

Finally, recalling \eqref{EQ4.2}--\eqref{eq4.14}, we estimate 
$\langle B, \tilde{B} \rangle_{L^2 }$. 
We denote by
\begin{equation}
    f_M (s) = P_M (|u|^{p-1} u (s)). 
    \label{eq5.29}
\end{equation} 
Explicitly, 
$\langle B, \tilde{B} \rangle_{L^2 }$ is given by
\begin{align*}
  \left \langle 
 \int_{\Lambda/M}^{\infty} e^{i \tau \sqrt{- \Delta}} \chi (4x/ \tau) D^{s_p -1} f_M (\tau)\, d\tau, 
 \int_{- \infty}^{ - \Lambda /M}
 e^{i t \sqrt{-\Delta}} \chi (4x/ t) D^{s_p -1} 
 f_M (t) \,  d t 
 \right \rangle_{L^2} .
\end{align*}
Similar to  \cite[Equation (4.19)]{CK23}, we observe that it is sufficient to estimate the real part of $\langle B, \tilde{B} \rangle_{L^2 }$ since the left hand side of the equation \eqref{EQ4.2} is given by
\begin{align*}
\langle A+B, \tilde{A} + \tilde{B} \rangle_{L^2}
=  
\langle P_M v(0), P_M v(0) \rangle_{\dot{H}^{s_p}} \geq 0.
\end{align*}
We remark that both $f_M (t)$  and $D^{s_p -1 } f_M (t)$ are real-valued because the corresponding multipliers for $P_M$ and $D^{s_p -1}$ are real valued functions. Furthermore, using the commutativity of $e^{i \tau \sqrt{- \Delta}}$ and $e^{i t \sqrt{- \Delta}}$, we may express $\mbox{Re} (\langle B, \tilde{B} \rangle_{L^2} ) $ by
\begin{align*}
    \int_{\Lambda /M}^{\infty} \int_{- \infty}^{- \Lambda/ M} 
    \langle  \chi (4x/ \tau) D^{s_p -1} f_M (\tau),
   \cos ((\tau - t) \sqrt{- \Delta}) \chi (4x/ t) D^{s_p -1} f_M (t)
    \rangle_{L^2} \, dt \, d\tau .
\end{align*}
Therefore, to estimate $\mbox{Re} (\langle B, \tilde{B} \rangle_{L^2})$, we study the decay rate in $t$ and $\tau$ of $\cos ((\tau - t) \sqrt{- \Delta}) \chi (4x/ t) D^{s_p -1} f_M (t)$ in dimension four. 
Note that this agrees with the representation of solutions for the linear wave equation
\begin{align}
    \begin{split}
         & \partial_s^2 h - \Delta h = 0
        \comma  (s,x) \in \mathbb{R} \times \mathbb{R}^4 \\
        & \evaluat{\vec{h}}_{s=0} = (g_0, 0) 
    \end{split}
    \label{eq5.18}
\end{align}
for time $s= \tau - t$. Letting 
\begin{equation}
    g_0 (x) = \chi (4x/ t) D^{s_p -1} f_M (t)
    \label{eq5.19}
\end{equation}
we have 
\begin{equation*}
    \mbox{supp} (g_0) \subset \{x:~ |x|\leq |t|/2 \}. 
\end{equation*}
On the other hand, we keep in mind that 
\begin{equation*}
    \mbox{supp} ( \chi (4x/ \tau) D^{s_p -1} f_M (\tau)) \subset \{x:~ |x|\leq |\tau|/2 \}. 
\end{equation*}
By \cite[Theorem~4.2]{SS98}, the solution to the Cauchy problem \eqref{eq5.18}--\eqref{eq5.19} is given by
\begin{equation*}
    h(x,s) = \partial_s \left( \frac{1}{s} \partial_s \right) \left( s^3 \int_{B_1 (0)} g_0 (x+ ys) \frac{1}{\sqrt{1 - |y|^2}} dy \right). 
\end{equation*}
We write
\begin{equation*}
    h(x,s) = v_0 (x,s) + v_1 (x,s)
\end{equation*}
where
\begin{equation*}
v_0 (x,s) = \partial_s \left( \frac{1}{s} \partial_s \right) \left( s^3 \int_{ \sqrt{1 - \epsilon_0}\leq |y|\leq 1} g_0 (x+ ys) \frac{1}{\sqrt{1 - |y|^2}} dy \right)
\end{equation*}
for some $\epsilon_0 >0$ and $v_1 (x,s) = h(x,s) - v_0 (x,s)$. This setting lets us express
\begin{align}
    \begin{split}
\mbox{Re} (\langle B, \tilde{B} \rangle ) 
=
    \int_{\Lambda /M}^{\infty} \int_{- \infty}^{- \Lambda/ M} 
    \langle  \chi (4x/ \tau) D^{s_p -1} f_M (\tau),
   v_0 (x,s) + v_1 (x,s)
    \rangle_{L^2} \, dt \, d\tau
    . 
        \label{eq5.20}
    \end{split}
\end{align}
First, we focus on the contribution of $v_0$ on the right hand side of \eqref{eq5.20}. Recalling $s = \tau + |t|$, when $|x+ ys| \leq {|t|}/{2}$ and $|y| \geq \sqrt{1 - \epsilon_0}$, we observe that
\begin{align*}
    |x| \geq |y|s - |x+ ys| 
        \geq \sqrt{1 - \epsilon_0} (\tau + |t|) - \frac{|t|}{2}   
        \geq \sqrt{1 - \epsilon_0}~ \tau 
        + (\sqrt{1 - \epsilon_0} - \frac{1}{2}) |t|
\end{align*}
which implies that
\begin{align*}
    |x| \geq \sqrt{1 - \epsilon_0} ~ \tau \geq \frac{\tau}{2}
\end{align*}
provided that $\epsilon_0 \in (0, 3/4)$. Fixing $\epsilon_0 \in (0, 3/4)$, we find that
\begin{align}
    \begin{split}
     \int_{\Lambda /M}^{\infty} \int_{- \infty}^{- \Lambda/ M} 
    \langle  \chi (4x/ \tau) D^{s_p -1} f_M (\tau),
   v_0 (x,s) 
    \rangle_{L^2} \, dt \, d\tau
    =0. 
        \label{eq5.21}
    \end{split}
\end{align}
Next, we consider $v_1 (x,s)$. Note that
\begin{equation*}
    v_1 (x,s) = 
    \partial_s \left( \frac{1}{s} \partial_s \right) \left( s^3 \int_{B_1 (0)} g_0 (x+ ys) \tilde{\alpha}(y) dy \right)
\end{equation*}
where
\begin{align*}
    \tilde{\alpha} (y) =
    \begin{cases}
        \frac{1}{\sqrt{1 - |y|^2}} & \inon{if }~ |y| \leq \sqrt{1 - \epsilon_0} \\
        0 & \inon{if }~ |y| > \sqrt{1 - \epsilon_0}. 
    \end{cases}
\end{align*}
We may instead use a smooth cut-off function $\alpha \in C_0^{\infty} (B_1 (0))$ such that 
\begin{equation*}
    \alpha (y) = \tilde{\alpha} (y) \quad  \mbox{if  }~ |y| \leq \sqrt{1 - \epsilon_0} . 
\end{equation*}
This way, the contribution of $v_1 (x,s)$ in \eqref{eq5.20} is given by
\begin{align}
    \begin{split}
    \tilde{v}_1 (x,s) = 
    \partial_s \left( \frac{1}{s} \partial_s \right) \left( s^3 \int_{B_1 (0)} g_0 (x+ ys) \alpha (y) dy \right). 
        \label{eq5.22}
    \end{split}
\end{align}
In other words, we have
\begin{align}
    \begin{split}
 &     \int_{\Lambda /M}^{\infty} \int_{- \infty}^{- \Lambda/ M} 
    \langle  \chi (4x/ \tau) D^{s_p -1} f_M (\tau),
   v_1 (x,s) 
    \rangle \, dt \, d\tau
  \\ & \quad \quad \quad  = 
   \int_{\Lambda /M}^{\infty} \int_{- \infty}^{- \Lambda/ M} 
    \langle  \chi (4x/ \tau) D^{s_p -1} f_M (\tau),
   \tilde{v}_1 (x,s) 
    \rangle \, dt \, d\tau .
        \label{eq5.23}
    \end{split}
    \end{align}
Therefore, we now focus on $\tilde{v}_1 (x,s)$ and carry out the derivatives in $s$ in \eqref{eq5.22}. We obtain
\begin{align}
    \begin{split}
    \tilde{v}_1 (x,s)  & =  3 \int_{B_1 (0)} g_0 (x+ ys) \alpha (y) \, dy
    +  5s \int_{B_1 (0)} \nabla g_0 (x+ ys) \cdot y  \alpha (y) \, dy
    \\ & \quad \quad  + 
    s^2 \int_{B_1 (0)} \langle \nabla^2 g_0 (x+ ys) y, y \rangle \alpha (y) \, dy
    \\ & = I + II + III. 
        \label{eq5.24}
    \end{split}
\end{align}
First, we treat $I$, and we observe that
\begin{align}
    I = 3 \int_{\mathbb{R}^4} g_0 (x+z) \frac{\alpha (z/s)}{s^4} \, dz
    = \alpha_s * g_0
    \label{eq5.25}
\end{align}
with $\alpha_s (z) = 3 \frac{\alpha (z/s)}{s^4}$. 

For $II$, we note that
\begin{equation*}
    s \nabla g_0 (x + ys) = \nabla_y \, g_0 (x+ys).  
\end{equation*}
Integrating by parts in $y$, we get
\begin{align*}
    II = -5 \int_{\mathbb{R}^4} g_0 (x+ys) \divv (y \alpha (y)) \, dy .
\end{align*}
Since $\alpha (y)$ and $y \alpha (y)$ are compactly supported in $B_1 (0)$, we may express
\begin{align*}
    II = \int_{\mathbb{R}^4} g_0 (x +ys) \beta (y) \, dy
\end{align*}
where $\beta \in C_0^{\infty} (B_1 (0))$. By changing the variables $z =ys$ as done in \eqref{eq5.25}, we obtain
\begin{align}
    II = \beta_s * g_0 
    \label{eq5.26}
\end{align}
with $\beta_s (z) = \frac{\beta(z/s)}{s^4}$. 

We note that $III$ may also be treated in the same manner. We integrate by parts in $y$ twice to take out the derivatives from $g_0$. Analogously, we get
\begin{align}
    III = \gamma_s * g_0
    \label{eq5.27}
\end{align}
where $\gamma_s (z) = \frac{\gamma (z/s)}{s^4}$, with $\gamma \in C_0^{\infty} (B_1 (0))$. 
Denote by $\delta_s$ the sum of smooth cut-off functions in \eqref{eq5.25}--\eqref{eq5.27}, i.e.,  
\begin{equation*}
    \delta_s = \alpha_s + \beta_s + \gamma_s 
\end{equation*}
Then, summing up the formulae for $I$, $II$, and $III$, we calculate
\begin{align}
    \norm{\tilde{v}_1 (\cdot, s)}_{L^{\infty}} \les 
    \norm{g_0}_{L^1} \norm{\delta_s}_{L^{\infty}}
    \label{eq5.28}
\end{align}
by Young's inequality. Recalling \eqref{eq5.19} and the fact that $s= |t|+ \tau$, we get 
\begin{equation*}
    \norm{\tilde{v}_1 (\cdot, s)}_{L^{\infty}}
\les \left( \frac{1}{|t|+ \tau} \right)^4 \norm{D^{s_p -1} f_M (t)}_{L^1} . 
\end{equation*}
We may now revisit \eqref{eq5.20}, \eqref{eq5.21}, and \eqref{eq5.23} to bound $ \mbox{Re} (\langle B, \tilde{B} \rangle)$ as follows:
\begin{align*}
0\leq 
   | \mbox{Re} (\langle B, \tilde{B} \rangle) |
    & \les 
     \int_{\Lambda /M}^{\infty} \int_{- \infty}^{- \Lambda/ M} 
    \norm{ \chi (4x/ \tau) D^{s_p -1} f_M (\tau)}_{L^1}
   \norm{\tilde{v}_1 (x,s) }_{L^{\infty}}
    \,  dt \, d\tau 
    \\ & \les 
\norm{D^{s_p -1} f_M }_{L_t^{\infty} L^1} 
    \norm{D^{s_p -1} f_M}_{L_t^{\infty} L^{1}} 
  \int_{\Lambda /M}^{\infty} \int_{- \infty}^{- \Lambda/ M}   \left( \frac{1}{|t|+ \tau} \right)^4 
  \, dt \, d\tau.
\end{align*}
Using Bernstein's inequality on the first two factors above and recalling \eqref{eq5.29} and \eqref{EQ4.44p}, we obtain  
\begin{equation}
    | \mbox{Re} (\langle B, \tilde{B} \rangle) |
     \les \frac{M^{2 s_p}}{\Lambda^2} \norm{P_M (|u|^{p-1} u ) }^2_{L_t^{\infty} L^{1}} 
     \les \eta^4 {M^{2 s_p}} \norm{P_{> M/8p} u }^2_{L_t^{\infty} L^{2}}. 
     \label{eq5.30}
\end{equation}
Note that we also relied on the definition $\Lambda= \eta^{-2}$ above. 

By adding up the estimates in \eqref{EQ4.46p}, \eqref{eq4.40}, and \eqref{eq5.30}, we finally obtain 
the estimate
\begin{align}
\begin{split}
 \langle  P_M v(0), P_M v(0) \rangle_{\dot{H}^{s_p}} 
& \les 
\left( \eta^{1/2} M^{s_p} \norm{ P_{>M/8p} u}_{L^{\infty}_t L^2} + \eta^{ -2} M^{s_p} N_0^{-s_p} \right)^2
 \\ & \quad  +
\left( \eta^{1/2} M^{s_p} \norm{ P_{>M/8p} u}_{L^{\infty}_t L^2} + \eta^{ -2} M^{s_p} N_0^{-s_p} 
\right) 
M^{s_p} \norm{ P_M v}_{L^{\infty}_t L^2}
 \period
 \end{split}
 \label{eq4.43}
 \end{align}    
The upper bound on the right hand side above, as in \cite[Section~4]{CK23}, gives us uniform control in time. As a result, we have established
 \begin{align}
 \begin{split}
 \langle  P_M v(t), P_M v(t) \rangle_{\dot{H}^{s_p}} 
& \les 
\left( \eta^{1/2} M^{s_p} \norm{ P_{>M/8p} u}_{L^{\infty}_t L^2} + \eta^{ -2} M^{s_p} N_0^{-s_p} \right)^2
 \\ & \quad  +
\left( \eta^{1/2} M^{s_p} \norm{ P_{>M/8p} u}_{L^{\infty}_t L^2} + \eta^{ -2} M^{s_p} N_0^{-s_p} 
\right) 
M^{s_p} \norm{ P_M v}_{L^{\infty}_t L^2}
 \period
 \end{split}
 \label{eq4.44}
 \end{align}
 We may close the proof of Proposition~\ref{decay_proposition} as in \cite{CK23}. For that end, we revisit the definitions of $a_j$ and $\beta_k$ in Claim~\ref{claim_1}. Recalling $ M = 2^j$ 
 for $j \in \mathbb{Z}^{-}$, the above estimate turns into
\begin{align*}
\begin{split}
a_j^2  & \les 
\left( \eta^{1/2} M^{s_p} \norm{ P_{>M/8p} u}_{L^{\infty}_t L^2} + \eta^{ -2} M^{s_p} N_0^{-s_p} \right)^2
 \\ & \quad  +
\left( \eta^{1/2} M^{s_p} \norm{ P_{>M/8p} u}_{L^{\infty}_t L^2} + \eta^{ -2} M^{s_p} N_0^{-s_p} 
\right) 
M^{s_p} \norm{ P_M v}_{L^{\infty}_t L^2}
\\ & \les 
\left( \eta^{1/2} \sum_{i > j- \ell_p} 2^{s_p (j-i)} a_i + 2^{s_p j} \eta^{-2 } N_0^{-s_p}  \right)^2 
\\ & \quad + 
a_j \left(  \eta^{1/2}   \sum_{i > j-\ell_p} 2^{s_p (j-i)} a_i  + 2^{s_p j} \eta^{-2 } N_0^{-s_p} \right)
\end{split}
\end{align*}
with $\ell_p = [\log (8p)] +1$ for $j <0$. It implies that
\begin{align*}
a_j \les 
 \eta^{1/2}   \sum_{i > j-\ell_p} 2^{s_p (j-i)} a_i  + 2^{s_p j} \eta^{-2 } N_0^{-s_p}
\end{align*} 
 for $j <0$.
 For $j >0$, we will simply use the estimate
 \begin{align}
 a_j \cong \norm{P_j v}_{L^{\infty}_t \dot{H}^{s_p}} \les 1 \period
\end{align}  
 By \eqref{eq4.7}-\eqref{eq4.8}, we obtain for $k <0$ 
  \begin{align}
  \begin{split}
  \beta_k &  \les \sum_{j>0} 2^{- |j-k|} + \eta^{1/2}  \left( \sum_{j <0} 2^{- |j-k|}  \left( \sum_{i > j- \ell_p} 2^{- s_p |j-i|} a_i \right) \right)
  \\ & \quad \quad \quad 
  + \eta^{-2 } N_0^{-s_p}  \sum_{j <0} 2^{- |j-k|} 2^{s_p j} 
  \\ & \les 
  \eta^{1/2} \beta_k + \sum_{j>0} 2^{- |j-k|} + \eta^{-2} N_0^{-s_p}   \sum_{j<0} 2^{- |j-k|} 2^{s_p j} 
  \period
  \end{split}
  \label{eq4.45}
  \end{align}
  Setting $\eta = 1/4C^2$, where $C$ is the implicit constant in \eqref{eq4.45}, we absorb the first term 
 on the last line above into the  left hand side, and we obtain 
 \begin{align*}
 \beta_k \les \sum_{j} 2^{-|j-k|} \min (1, 2^{s_p j})
 \end{align*}
 which yields
 \begin{align*}
 \beta_k \les 2^k \inon{ for}~ k <0 \period
 \end{align*}
 As we set $\beta_k =1$ for $k \geq 0$, we conclude that
 $\{ 2^{-3k/4} \beta_k\}_{k \in \mathbb{Z}} \in \ell^2$,
 which
 completes the proof of Proposition~\ref{decay_proposition}.
\end{proof}

The main consequence of Proposition~\ref{decay_proposition} is given in the next result, which lands the trajectory set of $\vec{u} (t)$ in the energy space $\dot{H}^1 \times L^2 (\mathbb{R}^4)$. Corollary~\ref{Cor7.1}
 is achieved as a result of 
 boundedness in $\dot{H}^{3/4} \times \dot{H}^{-1/4} (\mathbb{R}^4)$ and pre-compectness in $\dot{H}^{s_p} \times \dot{H}^{s_p -1} (\mathbb{R}^4)$. 
\begin{Corollary}\label{Cor7.1}
    Let $\vec{u} (t)$ be as in Proposition~\ref{Proposition1}. Then we have $\vec{u} (t) \in \dot{H}^1 \times L^2 (\mathbb{R}^4)$ for all $t \in \mathbb{R}$. Moreover, the trajectory
    \begin{equation}
        K_1 = \{\vec{u} (t): t \in \mathbb{R} \}
    \end{equation}
    is pre-compact in $\dot{H}^1 \times L^2 (\mathbb{R}^4)$. As a result, we have for all $R>0$
    \begin{equation*}
        \limsup_{t \to + \infty} \norm{\vec{u} (t)}_{\mathcal{H} (r \geq R + |t|)} 
        = \limsup_{t \to - \infty} \norm{\vec{u} (t)}_{\mathcal{H} (r \geq R + |t|)} 
        =0 .
    \end{equation*}
\end{Corollary}

An analogous version of this result is given in \cite[Corollary~7.1]{CK23}. Additionally, the arguments in the proof of Corollary~7.1 in \cite{CK23} stay valid in dimension four. For that reason, we refer the readers to \cite{CK23} for a proof.

\section{Preliminaries for the Rigidity Method}
\label{sec: rigidity}

In this section, we begin with the proof of Proposition~\ref{Proposition1}. The proof mainly follows the roadmap introduced in Sections~4--5 of \cite{DKMM22} and some preliminary results from \cite{CK23}. As we appeal to some of the results in \cite{DKMM22}, we will also rely on its notation. Below, we recall the function spaces commonly used in Sections~\ref{sec: rigidity}--\ref{sec: rigidity section 2}. 

\subsection{Notation and solutions outside a wave cone}
\label{subsec: notation}
The function spaces $L^p (\mathbb{R}^4)$ and $L^p (\mathbb{R}, L^q (\mathbb{R}^4))$ will simply be denoted by $L^p$ and $L^p L^q$, respectively. For the energy space, we will use 
$$\dot{H}^1 (\mathbb{R}^4) \times L^2 (\mathbb{R}^4) = \mathcal{H} $$
for shorthand.

Let $R>0$. Given a function space $A$ defined on $\mathbb{R}^4$, we denote by $A(R)$ the subspace of $A$ that consists of radial functions restricted to the domain $\{|x| \geq R \}$. Analogously, given a space-time function space $B$ on $\mathbb{R}_t \times \mathbb{R}_x^4$, we denote by $B(R)$ the subspace of functions radial in the space variable restricted to the outer cone $\{(t,x): |x|> R+|t|\}$. For instance, $L^p (R)$, $\dot{H}^1 (R)$, and $L^pL^q (R)$ will be the most commonly used spaces. The following norms will accompany these spaces:
\begin{align*}
     \Vert \varphi \Vert_{L^p (R)}^p = \int_R^{\infty} |\varphi(r)|^p r^3 \,dr
     \comma 
      \Vert \varphi \Vert_{\dot{H}^1}^2 = \int_R^{\infty} \left| \frac{d \varphi}{dr} (r) \right|^2 r^3 \,dr
\end{align*}
\begin{equation*}
      \Vert f \Vert_{L^p L^q (R)}^p = \int_{\mathbb{R}} \left(  \int_{R + |t|}^{\infty} |f (r)|^q r^3 \, dr \right)^{\frac{p}{q}} \, dt.
\end{equation*}


In Sections~\ref{sec: rigidity section 1}--\ref{sec: rigidity section 2}, we will be concerned with radial solutions of the nonlinear wave equation outside a wave cone. Let $\vec{u} (t)$ be as in Proposition~\ref{Proposition1}. By Lemma~\ref{lem7.4} and Remark~\ref{R5.6}, we will argue that
$\vec{u} (t)$ solves the equation below for every $R>0$:
\begin{align}
   \begin{split}
    & \partial^2_{t} u
    - \partial^2_{r} u - \frac{3}{r} \partial_r u
    = |u|^{p-1} u
      \comma r > R + |t|,
    \\
   &  \vec{u} (0)
    = (u_0, u_1) \in \mathcal{H} (R)  \comma r> R.
    \end{split}
   \label{equation4.1}
  \end{align}

\subsection{Singular stationary solutions}
\label{sec:stationary solutions}
Below, we state an analogous version of \cite[Proposition~6.1]{CK23} for dimension four. As in Proposition~6.1 in \cite{CK23}, the result below yields a family of smooth radial solutions to \eqref{SS_eq1}, the underlying  stationary equation away from the origin. Due to a singularity at $r=0$, these solutions fail to belong in the critical Sobolev space $\dot{H}^{s_p} (\mathbb{R}^4)$. This fact plays a key role in the energy-supercritical applications of the rigidity method, including Sections~\ref{sec: rigidity section 1}--\ref{sec: rigidity section 2} of the present article. 
In contrast,  the ground state solution 
\begin{equation*}
W(r) = \frac{1}{\left( 1+ \frac{r^2}{8} \right)}  
\end{equation*}
is a stationary radial solution for the energy-critical nonlinear wave equation and $W$ belongs to the energy space $\dot{H}^1 (\mathbb{R}^4)$.

\begin{Proposition}
\label{prop_singular_solns}
For any $l \in \mathbb{R}\backslash \{0\}$ there exists a radial $C^2$ solution $Z_l$ of 
\begin{equation}
\Delta Z_l + |Z_l|^{p-1} Z_l =0 \inon{ in }~  \mathbb{R}^4\backslash \{0\}
\label{SS_eq1}
\end{equation}
with the asymptotic behavior 
\begin{equation}
r^{2} Z_l (r) = l + O \left( r^{4 - 2 p}\right)  \inon{ as }~  r \to \infty
\period
\label{SS_eq2}
\end{equation}
Furthermore, $Z_l \notin L^{q_p} (\mathbb{R}^4)$, where $q_p := {2 (p-1)}$ is the critical Sobolev exponent corresponding to $\dot{H}^{s_p}$. This implies that $Z_l \notin \dot{H}^{s_p} (\mathbb{R}^4)$. 
\end{Proposition}

Proposition~6.1 in \cite{CK23} is stated
with the underlying assumption that $d \geq 7$ and $p$ is an odd integer. Nevertheless, its proof remains valid for $d > 4$ and $p \geq 3$ or $d=4$ and $p >3$.

\subsection{Channels of energy for the radial linear inhomogeneous wave equation }
An essential part of the rigidity method is to obtain exterior energy estimates for radial solutions of the linear problem. In dimension four, we may rely on Lemma~3.9 in \cite{DKMM22} (see also \cite{LSW24}). 
For the reader's convenience, we recall the setting and the statement of Lemma 3.9 from \cite[Section~3]{DKMM22} with minor changes to avoid confusion due to notational overlaps. 

We consider radial solutions of the linear inhomogeneous wave equation
\begin{align}
    \begin{split}
          & \partial_t^2 v - \Delta v =  
         f
        \comma  (t,x) \in \mathbb{R} \times \mathbb{R}^4 \\
        & \evaluat{\vec{v}}_{t=0} = (v_0, v_1) \in \mathcal{H}.
    \end{split}
    \label{eq5.11}
\end{align}
Let $R>0$. We note that $1/r^2$ is a solution of the free wave equation on the channel $\{|x|> |t| +R \}$, i.e., it solves the equation \eqref{eq5.11} with initial data $(1/r^2, 0)$ and $f=0$.  Let ${span}_R (1/r^2)$ denote the subspace of $\dot{H}^1 (R)$ spanned by $1/r^2$. 
\begin{Notation} Let $R>0$. We denote by $\pi_R f$ the orthogonal projection of the function $f$ on the one-dimensional subspace $\mbox{span}_R (1/r^2)$ under the $\dot{H}^1 (R)$ inner product. More specifically, we have the formula
\begin{equation}
    \pi_R f (r) = \frac{R^2 f(R)}{r^2} \comma r> R. 
    \label{eq5.14}
\end{equation}
   Additionally, $\pi_R^{\perp}$ denotes the orthogonal projection onto the complement of $\mbox{span}_R (1/r^2)$, and it is given by
   \begin{equation}
    \pi^{\perp}_R f (r) = f(r) - \frac{R^2 f(R)}{r^2} \comma r> R. 
       \label{eq5.15}
   \end{equation}
\end{Notation}
The following exterior energy estimate result is from \cite[Lemma~3.9]{DKMM22}. See also \cite[Proposition~1.11]{LSW24} for radial exterior energy estimates in even dimensions $d \geq 2$. 
\begin{Lemma}\label{Lemma3.9} 
Let $v_0 \in \dot{H}^1$ and $f \in L^1 ((0, \infty), L^2 (\mathbb{R}^4))$ have radial symmetry. Let $\vec{v} (t,r)$ be the solution of \eqref{eq5.11} with initial data $(v_0, 0)$. Then, 
\begin{equation}
    \norm{v_0}_{\dot{H}^1 (\mathbb{R}^4)}
    \leq 2 \lim_{t \to + \infty} \norm{\nabla v (t)}_{L^2 (|x|> t)} + 2 \norm{ f}_{L^1 ((0, \infty), \, L^2 (|x|>t)  )} 
    \label{eq5.12}
\end{equation}
and for $R>0$     
\begin{equation}
  \Vert{ \pi_R^{\perp}v_0} \Vert_{\dot{H}^1 ({R})}
    \leq \sqrt{\frac{20}{3}}
    \left( \lim_{t \to + \infty} \norm{\nabla v (t)}_{L^2 (|x|> t +R)} 
    + \norm{ f}_{L^1 ((0, \infty), \, L^2 (|x|>t + R)  )} 
    \right).
    \label{eq5.13}
\end{equation}
\end{Lemma}
The above estimate
establishes a lower bound of the asymptotic energy over channels $\{|x|> |t| + R\}$ of width $R$ for any $R>0$. We remark that this result is applicable for radial solutions of \eqref{eq5.11} when the data is of the form $(v_0, 0)$. 

\subsection{Cauchy problem outside a wave cone}
\label{subsec: Cauchy problem}

In this subsection, we will review the local wellposedness theory for radial solutions outside a wave cone, with finite energy. 

 First, we select a radial cut-off function $\chi \in C^{\infty} (\mathbb{R}^4)$ that agrees with 
\begin{align*}
    \chi (r) =
    \begin{cases}
        1 & \inon{if }~ r \geq 1/{2} \\
        0 & \inon{if }~ r \leq 1/{4} 
    \end{cases}
\end{align*}
For $r_0>0$, we denote by $\chi_{r_0} (r) = \chi (r/r_0)$.
Let $I \subset \mathbb{R}$ be an interval with $0 \in I$.
We then consider the following Cauchy problem on an interval $I \subset \mathbb{R}$ with $0 \in I$:
\begin{align}
    \begin{split}
        & \partial_t^2 h - \Delta h = 
        |V + \chi_{r_0} h|^{p-1} (V+  \chi_{r_0} h) - |V|^{p-1} V   
        \comma  \\
        & \evaluat{(h , \partial_t h)}_{t=0} = (h_0, h_1) \in \mathcal{H}.
    \end{split}
    \label{eq7.9}
\end{align}
We remark that when $V=0$ and $r_0 = R$, the nonlinearity for the Cauchy problem above satisfies
\begin{equation}
|\chi_{r_0}|^{p-1} \chi_{r_0} h =
   \left|\chi_{R+|t|}  h\right|^{p-1} \chi_{R+|t|} h = |h|^{p-1} h \comma
   \mbox{on }~
   r \geq R + |t|
   \label{eq7.30}
\end{equation} 
where $\chi_{R +|t|} (r) = \chi \left( r / (R+ |t|) \right)$.
 By the finite speed of propagation, we observe that solutions to \eqref{eq7.9} and \eqref{equation4.1} agree on the exterior cone $\{(t,r):~ r > R + |t| \}$
 provided that 
 \begin{equation}
     (h_0, h_1) = (u_0, u_1) \comma r > R. 
     \label{eq7.31}
 \end{equation}

We will also utilize \eqref{eq7.9} with 
\begin{equation}
    V(x, t) = \chi \left( \frac{|x|}{R+|t|}\right) Z_{\ell} 
\end{equation}
where $Z_{\ell}$ is as in Proposition~\ref{prop_singular_solns}. Below, we include a small data result addressing \eqref{eq7.9} in dimension four. We note that analogous statements were shown in \cite{DKM14, DL15, CK23} for odd dimensions $d \geq 3$. 
\begin{Notation}
Let  $q \in [1, \infty]$. We denote by $L^q_{I} := L^q \left(I, L^q (\mathbb{R}^4 ) \right)$ throughout Lemma~\ref{lem7.4} and its proof. This notation agrees with the one in  \cite[Lemma~7.4]{CK23}
\end{Notation}
\begin{Lemma}
\label{lem7.4}
    There exists $\delta_0 >0$ satisfying the following property: let 
    $V \in L^{\frac{5}{2} (p-1)}_I$ be a radial function in space satisfying 
    \begin{equation}
        \Vert D^{1/2} V \Vert_{L^{\frac{10}{3}}_I} \leq \delta_0 r_0^{\frac{p-3}{p-1}}~\mbox{ and }~ 
       \Vert V \Vert_{L^{\frac{5}{2} (p-1)}_I }
       \leq \delta_0 .
       \label{eq5.10}
    \end{equation} 
    Furthermore, let $(h_0, h_1) \in \mathcal{H}$ be radial functions with
    \begin{equation}
        \Vert (h_0, h_1) \Vert_{\mathcal{H}} \leq \delta_0 r_0^{\frac{p-3}{p-1}} .
        \label{eq5.6}
    \end{equation}
    Then the Cauchy problem \eqref{eq7.9} is well-posed on the interval $I$, and we have 
    \begin{equation}
        \sup_{t \in I} \Vert  h(t) - S(t) (h_0, h_1) \Vert_{\mathcal{H}} \leq \frac{1}{100} \Vert (h_0, h_1) \Vert_{\mathcal{H}}. 
        \label{eq5.4}
    \end{equation}
    Moreover, if $V=0$, we may have $I= \mathbb{R}$, and we obtain
    \begin{equation}
        \sup_{t \in I} \Vert  h(t) - S(t) (h_0, h_1) \Vert_{\mathcal{H}} \les \frac{1}{r_0^{p-3}} \Vert (h_0, h_1) \Vert^p_{\mathcal{H}}. 
        \label{eq5.5}
    \end{equation}
\end{Lemma}

\begin{proof}
  Even though Lemma~7.4 in \cite{CK23} addresses dimensions $d \geq 7$, its proof extends to cover Lemma~\ref{lem7.4}. More precisely, \eqref{eq5.4}--\eqref{eq5.5} may be obtained by simply setting $d=4$ in the proof of Lemma~7.4. For the convenience of the readers, we provide the corresponding norms and main estimates below. 

  To apply the fixed-point argument, we introduce the norm 
\begin{align}
\| h \|_{\mathcal{S}} := \norm{h}_{L^{5}_I} + {\| D^{1/2} h\|}_{L^{10/3}_I } + \sup_{t \in I} \norm{(h, h_t)}_{\mathcal{H}}
\label{eq5.1} 
\end{align}
and the space $B_\alpha$ 
\begin{align*}
B_{\alpha} := \{h \in L^{5}_I: h \mbox{ is radial},~ \| h \|_{\mathcal{S}}  \leq \alpha \}. 
\end{align*}
for $\alpha >0$. We note that the norms in \eqref{eq5.1} correspond to the admissible triples $(5,5,1)$, $(10/3, 10/3, 1/2)$, and $(\infty, 2, 0)$. By Duhamel's principle, we define 
$\Phi (v) (t)$ as in \cite[Lemma~7.4]{CK23}, and 
applying the Strichartz estimates with $(h_0, h_1) \in \mathcal{H}$, we find that
\begin{align}
\| \Phi (v) \|_{\mathcal{S}}  \leq C ( \norm{(h_0, h_1)}_{\mathcal{H}} + \|{ D^{1/2} F_V (v)}\|_{L^{10/7}_I} )
\label{eq5.2}
\end{align}
where $F_V (v) = G (V + \chi_{r_0} v) - G(V)$ and $G(h) = |h|^{p-1} h$. 
The rest of the proof proceeds as in \cite{CK23}. For instance, we estimate
\begin{align}
\begin{split}
& \| D^{1/2} F_V (v) \|_{L^{10/7}_I } 
=
\| D^{1/2} (G (V + \chi_{r_0} v) - G(V)) \|_{L^{10/7}_I} 
\\ & 
 \indeq \indeq \indeq \leq
C \left( \|  G' (V + \chi_{r_0} v) \|_{L^{5/2}_I}    +   \| G' (V) \|_{L^{5/2}_I } \right) \| D^{1/2} (\chi_{r_0} v) \|_{L^{10/3}_I }
\\ & 
 \indeq \indeq \indeq +
 C \left( \|  G'' (V + \chi_{r_0} v) \|_{L^{{\frac{5}{2} \frac{(p-1)}{(p-2)}}}_I}  
 +   \| G'' (V) \|_{L^{{\frac{5}{2} \frac{(p-1)}{(p-2)}}}_I } \right)  \| \chi_{r_0} v \|_{L^{\frac{5}{2} (p-1)}_I }
 \\ &  \indeq \indeq \indeq \indeq \indeq  \indeqtimes
 \left( \| D^{1/2}  (V + \chi_{r_0} v)\|_{L^{10/3}_I}  + \| D^{1/2} (V) \|_{L^{10/3}_I} \right).
\end{split}
\label{eq5.3}
\end{align}
Similarly, we may follow the same line of arguments and estimates as in \cite{CK23} to close the proof. 
\end{proof}

\begin{Remark}\label{R5.6} 
Let $\vec{u} (t)$ be as in Proposition~\ref{Proposition1}. By Corollary~\ref{Cor7.1}, Lemma~\ref{lem7.4}, and the finite speed of propagation, we can deduce that $\vec{u} (t)$ solves ~\eqref{equation4.1}~ for every $R>0$.  It is worthwhile to note that 
the proof of Lemma~\ref{lem7.4} also yields the following estimate: 
    \begin{align}
        \begin{split} 
  \sup_{t \in \mathbb{R}} \Vert  \vec{u} (t) - S(t) (u_0, u_1) \Vert_{\mathcal{H} (R + |t|) } 
 & + \Vert  u(t) - S(t) (u_0, u_1) \Vert_{L^2 L^8 (R)}  \\
&  + \Vert  u(t) - S(t) (u_0, u_1) \Vert_{L^3 L^6 (R)} 
  \les \frac{1}{R^{p-3}} \Vert \vec{u} (t) \Vert^p_{\mathcal{S} (R)} \\
& \qquad  \qquad \qquad \qquad \qquad \qquad 
\les \frac{1}{R^{p-3}} \Vert (u_0, u_1) \Vert^p_{\mathcal{H} (R) } 
            \label{eq5.7}
        \end{split}
    \end{align}
    To see this, we modify the $\norm{\cdot}_{\mathcal{S}}$ norm as follows:
    \begin{align*}
\| h \|_{\mathcal{S}} := \norm{h}_{L^2 L^8} +
\norm{h}_{L^3 L^6} +\norm{h}_{L^{5} L^5 } + {\| D^{1/2} h\|}_{L^{10/3} L^{10/3} } + \sup_{t \in \mathbb{R}} \norm{(h, h_t)}_{\mathcal{H}}. 
\end{align*}
The Strichartz estimates combined with the upper bound on estimate (7.20) in \cite[pg.~43]{CK23} lead to \eqref{eq5.7}. 

Furthermore, the Strichartz estimates for the linear problem yields the following estimate outside the cone:
\begin{align*}
\Vert  S(t) (u_0, u_1) \Vert_{L^2 L^8 (R)}  +
    \Vert  S(t) (u_0, u_1) \Vert_{L^3L^6 (R)}  + 
    \sup_{t \in \mathbb{R}} \Vert  \vec S(t) (u_0, u_1) \Vert_{\mathcal{H} (R + |t|)} 
    \les \Vert  (u_0, u_1) \Vert_{\mathcal{H} (R )} . 
\end{align*}
As a result, estimate \eqref{eq5.7}  implies that
\begin{align}
    \Vert  {u} (t) \Vert_{L^2 L^8 (R)} +
     \Vert  {u} (t) \Vert_{L^3L^6 (R)}  + 
    \sup_{t \in \mathbb{R}} \Vert  \vec  u (t) \Vert_{\mathcal{H} (R + |t|)} 
    \les \Vert  (u_0, u_1) \Vert_{\mathcal{H} (R )} . 
    \label{eq5.8}
\end{align}
\end{Remark}

The next result is the main rigidity theorem for our paper. 

\begin{Proposition}\label{rigidity proposition}
   Suppose that $u(t,r)$ is a radial solution of \eqref{wave equations}
with $\vec{u} (t) \in C_b (\mathbb{R}, \dot{H}^{s_p} \times \dot{H}^{s_{p}-1})$. Suppose that for all $R_0 >0$, $u(t,r)$ is a radial solution of  
\eqref{equation4.1}  for $r> R_0 + |t|$, with initial data in $\mathcal{H}(R_0)$ and $\norm{\vec{u} (t)}_{\mathcal{S}(R_0)} < \infty$. 
Assume that
       \begin{equation}
        \sum_{\pm} \lim_{t \to \pm \infty}     \int_{r > R_0 +|t|} |\nabla_{t,r} u (t,r)|^2  r^3 \, dr =0.
\label{eq5.16}
    \end{equation}
Then, there exists $R > 0$ such that 
    \begin{equation}
        (u_0 (r), u_1(r)) = (0,0) \comma \mbox{ for all }~ r > R. 
    \end{equation}
\end{Proposition}
The proof of Proposition~\ref{rigidity proposition} is divided into two key results, namely Proposition~\ref{Prop4.5} and Proposition~\ref{Prop5.1}, and it will be carried out through Sections~\ref{sec: rigidity section 1}--\ref{sec: rigidity section 2}. We close this section by showing that Proposition~\ref{rigidity proposition} implies Proposition~\ref{Proposition1}.

\begin{proof}[Proof of Proposition~\ref{Proposition1}] 
Let $\vec{u} (t)$ be as in Proposition~\ref{Proposition1}. As argued in Remark~\ref{R5.6} and Corollary~\ref{Cor7.1}, $\vec{u} (t)$ verifies the hypotheses of Proposition~\ref{rigidity proposition}. Therefore, there is $R>0$ so that
\begin{equation}
        (u_0 (r), u_1(r)) = (0,0) \comma \mbox{ for all }~ r > R. 
    \end{equation}
For a contradiction, suppose $(u_0, u_1)$ is not identically zero. Then, we deduce that
    \begin{equation}
    \rho= \inf \{\sigma>0: \quad  \int_{r>\sigma } |\nabla u_0|^2 + |u_1|^2 =0 \}  \in (0, R). 
\end{equation}
Next, we take a sufficiently small $\epsilon_0 >0$ such that
\begin{equation}
    0 < \norm{(u_0, u_1)}_{\mathcal{H}(\rho - \epsilon_0)} \ll 1
\end{equation}
and consider the solution $\vec{h}(t,r)$ to the Cauchy problem \eqref{eq7.9} with $V=0$, $r_0 = \rho - \epsilon_0$, and  with data
\begin{align*}
    (h_0 (r), h_1 (r)) = 
    \begin{cases}
         (u_0 (r), u_1 (r)) \comma & r > \rho - \epsi_0
        \\
         (u_0 (\rho - \epsilon), 0) \comma  & 0< r \leq  \rho- \epsilon_0. 
    \end{cases}
\end{align*}
By Lemma~\ref{lem7.4} and  the finite speed of propagation, the solution $\vec{h}(t,r)$ agrees with $\vec{u} (t,r)$ for $r\geq |t|+ \rho- \epsilon_0$, and by \eqref{eq5.8}, we have 
\begin{equation}
\sup_{t \in \mathbb{R}} \norm{\vec{u} (t)}_{\mathcal{H} (|t|+ \rho- \epsilon_0)} 
\les \norm{(u_0, u_1)}_{\mathcal{H}(\rho - \epsilon_0)} \ll 1. 
\label{eq5.17}
\end{equation} 
Combined with the radial Sobolev inequality, the previous estimate implies that there is a positive constant $C$ such that
for all $t \in \mathbb{R}$, 
\begin{equation}
    |u (t,r)|^{p-1} \leq \frac{C}{r^{ p-1}}
    \comma \mbox{ for all }~  r> |t|+ \rho- \epsilon_0
    \label{eq5.31}
\end{equation} 
where $p-1 >2$. 
 The setting in this proof, combined with \eqref{eq5.17}--\eqref{eq5.31}, allows us to  use \cite[Proposition~4.7]{DKM20}. 
 We remark that while \cite{DKM20} deals with odd dimensions, Proposition~4.7~is actually valid in both even and odd dimensions, with a proof that does not distinguish among them. Hence, by \cite[Proposition~4.7]{DKM20} we arrive  at 
 \begin{equation}
     \sum_{\pm}
     \lim_{t \to \pm \infty}
     \int_{r \geq |t|+ \rho - \epsilon_0 } |\nabla_{t,r} u(t,r)|^2 r^3 \, dr
     \geq 
     \frac{1}{8} \int_{r\geq \rho- \epsilon_0} 
     (|\nabla u_0|^2 + |u_1 |^2) r^3 \, dr >0
 \end{equation}
which contradicts \eqref{eq5.16} with $R_0 = \rho - \epsilon_0$. 
\end{proof}






\section{Rigidity Method Part I: Constant Sign Solutions} \label{sec: rigidity section 1}

\begin{Proposition}\label{Prop4.5}
Let $\vec{u} (t,r)$ be as in Proposition~\ref{rigidity proposition}. Then, 
there exist $t_0 \in \mathbb{R}$ and $R> 0$ so that
\begin{equation*}
    u (t_0, r) = 0 \quad \quad \mbox{for all }~ r > |t_0| +R.
\end{equation*}
\end{Proposition}
We start by proving:
\begin{Lemma}\label{Lemma4.7}
   Let $R>0$ and $\vec{u} (t)$ be a radial solution of \eqref{equation4.1} on $\{ r> R+ |t| \}$.  There exists $\varepsilon_0 >0$ with the following property: assume that
   \begin{equation}
       \int_{|x| > R} |\nabla u_0 (x)|^2 \, dx \leq \varepsilon_0^2
       \label{eq6.1}
   \end{equation}
   and 
   \begin{equation}
       \sum_{\pm} \lim_{t \to \pm \infty}  \int_{|x| > R + |t|} |\nabla_{t,x} u (t, x)|^2 \, dx =0.
       \label{eq6.2}
   \end{equation}
   Then, we have
   \begin{equation}
       |u_0 (\rho)| \simeq \frac{1}{\rho} \norm{u_0}_{L^4 (\rho)}
       \simeq \frac{1}{\rho} \norm{u_0}_{\dot{H}^1 (\rho)} 
       \comma \mbox{for all }~ \rho \geq R
       \label{eq6.3}
   \end{equation}
   \begin{equation}
       |u_0 (r)| \leq 2 \left( \frac{\rho}{r} \right)^{\frac{11}{6}} |u_0 (\rho)| 
       \comma \mbox{for all }~ r \geq \rho \geq R. 
       \label{eq6.4}
   \end{equation}
\end{Lemma}

\begin{proof}
    Let 
    \begin{align*}
        u_{+} (t) = \frac{u(t) + u(-t)}{2} \quad \mbox{ and } 
        u_{-} (t) = \frac{u(t) - u(-t)}{2} .
    \end{align*}
    Then, $u_+$ solves the following system of equations:
    \begin{align*}
        \begin{cases}
           & \partial_t^2 u_+ - \Delta u_+ = \frac{ F(u (t))}{2} + \frac{F( u (-t) )}{2} \\
           & \vec{u}_+ (0) = (u_0, 0)
        \end{cases}
    \end{align*}
    where $F(x) = |x|^{p-1} x$. Recall that $F$ is an odd function with $F' (x) = p |x|^{p-1}$. Letting $a = u(t)$ and $b = u(-t)$, we observe that
    \begin{align*}
        |F(a) + F(b)| & = \left||a|^{p-1} a - 
        \left(|-b|^{p-1} (-b) \right)  \right| \\
        & \leq   |a+b|   \int_0^1 \left| F' (s a - (1-s)b) \right| \, ds \\
        &  \les   |a+b| \left( |a|^{p-1} + |b|^{p-1} \right)
    \end{align*}
    leading to 
    \begin{align}
        \frac{1}{2} \left|  F(u (t))+  F(u (-t)) \right| \les |u_+ (t)| 
        \left( |u_+ (t)|^{p-1} + |u_- (t)|^{p-1} \right). 
        \label{eq6.5}
    \end{align}

Estimates \eqref{eq6.3}--\eqref{eq6.4} may be verified by following the same arguments as in Lemma~4.7 in \cite{DKMM22}. As discussed in Remark~\ref{R5.6}, we have 
\begin{equation}
  \Vert u \Vert_{L^3 L^6 (\rho)} 
  + 
  \sup_{t\in \mathbb{R}} \Vert \vec{u} (t) \Vert_{\mathcal{H} (\rho + |t|)}
  \les 
  \Vert(u_0, u_1) \Vert_{\mathcal{H} (R)}
  \les \varepsilon_0
  \label{eq6.7}
\end{equation}
for all $\rho \geq R$. 
Also, by the exterior energy estimates in Lemma~\ref{Lemma3.9}, we get
\begin{equation}
    \norm{\pi_{\rho}^{\perp} u_0}_{\dot{H}^1{(\rho)}}
\les \lim_{t \to \infty} \norm{\nabla u_+ (t)}_{L^2 (\rho + |t|)} 
+ \norm{F (u (t)) + F(u (-t))}_{L^1 L^2 (\rho)}
\label{eq6.6}
\end{equation}
By \eqref{eq6.2}, the first term on the right side of \eqref{eq6.6} is zero. For the second term, we use \eqref{eq6.5}, \eqref{eq6.7}, and the radial Sobolev inequality, and we obtain the upper bound
\begin{align}
    \begin{split}
 \frac{1}{\rho^{p-3}}
 \sup_{t \in \mathbb{R}} \norm{u (t)}_{\dot{H}^1{(\rho)}}^{p-3}
 \norm{u_+ (t)}_{L^3 L^6 (\rho)}
 \left( \norm{u_+ (t)}_{L^3 L^6 (\rho)}^2 + \norm{u_- (t)}_{L^3 L^6 (\rho)}^2 \right)
 \les \frac{\varepsilon_0^{p-1}}{\rho^{p-3}}\norm{u_0}_{\dot{H}^1 (\rho)} . 
    \end{split}
    \label{eq6.9}
\end{align}
Recalling $\pi_{\rho}^{\perp} u_0 (r) = u_0 (r) - \frac{\rho^2}{r^2} u_0 (\rho)$, we control the norm of $\pi_{\rho}^{\perp} u_0 (r)$ from below as follows:
\begin{equation}
    \norm{\pi_{\rho}^{\perp} u_0 (r)}_{\dot{H}^1 (\rho)} 
    \geq 
    \left| \norm{u_0 (r)}_{\dot{H}^1 (\rho)} -  
    \frac{\rho^2 |u_0 (\rho)|}{\sqrt{2} \rho}
    \right| .
    \label{eq6.8}
\end{equation}
As a result,  by \eqref{eq6.6}, \eqref{eq6.9}, and \eqref{eq6.8}, we obtain
\begin{equation}
    \norm{u_0}_{\dot{H}^1 (\rho)} \les \rho |u_0 (\rho)| + 
    \frac{\varepsilon_0^{p-1}}{R^{p-3}}
    \norm{u_0}_{\dot{H}^1 (\rho)} 
\end{equation}
Selecting $\varepsilon_0 = \varepsilon_0 (R, p) >0$ sufficiently small, we deduce that
\begin{equation}
    \norm{ u_0}_{\dot{H}^1 (\rho)} \leq C \rho |u_0 (\rho)|. 
    \label{eq6.10}
\end{equation}
By the Sobolev embedding, \eqref{eq6.6} and \eqref{eq6.9} also imply that 
\begin{equation}
 \norm{\pi_{\rho}^{\perp} u_0 (r)}_{L^4 (\rho)} =
    \norm{u_0 (r) - \frac{\rho^2}{r^2} u_0(\rho)}_{L^4 (\rho)} 
    \les \frac{\varepsilon_0^{p-1}}{R^{p-3}} \rho |u_0 (\rho)| .
    \label{eq6.11}
\end{equation}
Since we have
\begin{equation*}
    \norm{\frac{1}{r^2}}_{L^4 (\rho)} = \frac{1}{\sqrt{2} \rho}, \quad 
 \mbox{and } \quad
    \norm{\frac{1}{r^2}}_{\dot{H}^1 (\rho)} = 
    \frac{\sqrt{2}}{\rho}
\end{equation*}
estimate \eqref{eq6.11} leads to
\begin{equation} 
    \norm{u_0}_{L^4 (\rho)} \leq 
    \rho |u_0 (\rho)|
    \left( \frac{1}{\sqrt{2}} + \frac{C \varepsilon_0^{p-1}}{R^{p-3}} \right)
    \leq \frac{\rho |u_0 (\rho)|}{\sqrt{2} - \tilde{\varepsilon}_0}
\end{equation}
for some small $\tilde{\varepsilon}_0  (R, p) > 0$.

Moreover, estimate \eqref{eq6.11}  implies that
$ \rho |u_0 (\rho)| \les \norm{u_0}_{L^4 (\rho)}$.
Combined with the Sobolev embedding, this completes the proof of \eqref{eq6.3}. 

We skip the proof of \eqref{eq6.4} as it is identical to the proof of (4.21) in \cite[Lemma~4.7]{DKMM22}.
\end{proof}

\begin{Lemma}\label{Lemma4.8}
   Let $u$, $R$ and $\varepsilon_0$ be as in Lemma~\ref{Lemma4.7}. Assume that for all $t \in \mathbb{R}$, $u(t)$ is not identically $0$ on $[R+ |t|, \infty)$. Then, $u$ never vanishes and has a constant sign on $\{|x| > R + |t| \}$. 
   Furthermore, for all $t \in \mathbb{R}$, there exists $\ell (t) \neq 0$ such that
   \begin{equation}
       |r^2 u (t,r) - \ell (t)| \leq C r^4 |u (t,r)|^3 
       \label{eq6.12}
   \end{equation}
for all $r > R + |t|$. 
\end{Lemma}

\begin{proof}
    The proof of Lemma~4.8 in \cite{DKMM22} is largely applicable here. 
    In particular, the arguments in the beginning of Lemma~4.8 still hold true, yielding that for all $t$, $u(t)$ never vanishes on the set $\{r> R + |t| \}$. As a result, using the continuity of $u$, we deduce that $u(t,r)$ has a constant sign on the outer channel $\{r> R+ |t| \}$, and without loss of generality, we may assume that $u(t,r)>0$ for all $t \in \mathbb{R}$ and $r> R + |t|$. This implies that
    \begin{equation*}
        |u_- (t,r)| = \frac{1}{2} |u(t,r)- u (-t,r)| \leq \frac{u(t,r) + u (-t,r)}{2} = u_+ (t,r)
    \end{equation*}
   for all $t$ and $r> R+ |t|$. By following the estimates \eqref{eq6.6}-\eqref{eq6.9}, 
  we observe that for all $\rho \geq R$
  \begin{equation*}
      \norm{\pi_{\rho}^\perp u_0}_{\dot{H}^1_{\rho}} 
      \leq \frac{C}{\rho^{p-3}} \sup_{t \in \mathbb{R}} \norm{u (t)}^{p-3}_{\dot{H}^1_{\rho}} \norm{u_+}^3_{L^3 L^6 (\rho)}
      \leq \frac{C}{\rho^{p-3}} \norm{u_0}^{p}_{\dot{H}^1 (\rho)}
      .
  \end{equation*}
  Next, we use the definition of $\pi^{\perp}_{\rho}$ and Sobolev embedding, as done in the proof of Lemma~\ref{Lemma4.7}. For all $\rho \geq R$
  \begin{equation*}
      \left|\| u_0\|_{L^4 (\rho)}  - \frac{\rho}{\sqrt{2}} u_0 (\rho)    \right| 
      \leq C \norm{\pi_{\rho}^\perp u_0}_{\dot{H}^1_{\rho}} 
      \leq 
      \frac{C}{\rho^{p-3}} \norm{u_0}^{p}_{\dot{H}^1 (\rho)}
      \leq  \frac{ C \varepsilon_0^{p-3}}{R^p}
      \rho^3 |u_0 (\rho)|^3
  \end{equation*}
where we used \eqref{eq6.3} in the last step. We note that the estimate above agrees with the one in (4.36) from the proof of  \cite[Lemma~4.8]{DKMM22}. Hence, 
we may complete the proof by following the same  calculations done in \cite[Lemma~4.8]{DKMM22}. 
\end{proof}
\begin{Lemma}\label{Lemma4.10}
    Let $u$, $R$ and $\varepsilon_0$ be as in Lemmas~\ref{Lemma4.7}--\ref{Lemma4.8}. Then, 
    $\ell$ is independent of $t$ and 
    \begin{equation}
        |u (t,r)| + \frac{1}{r} \norm{u (t)}_{\mathcal{H} (R)} \les \frac{|\ell|}{r^2}
        \comma
        \mbox{for all }~ t\in \mathbb{R}~ \mbox{ and }~ r \geq R + |t|.
        \label{eq6.13}
    \end{equation}
\end{Lemma}
 For a proof, we refer our readers to \cite[Lemma~4.10, pg.~819-820]{DKMM22} since the proof of  Lemma~4.10 in \cite{DKMM22} applies to Lemma~\ref{Lemma4.10} without any changes.  

\begin{Lemma}\label{Prop4.6}
Let $R_0 >0$ and $u$ be a radial solution of \eqref{equation4.1} defined for $\{r> |t| + R_0 \}$ with $(u_0, u_1) \in \mathcal{H} (R_0)$. Assume that
\begin{equation}
     \sum_{\pm} \lim_{t \to \pm \infty}  \int_{|x| > |t|+ R_0} |\nabla_{t,x} u (t, x)|^2 \, dx =0.
       \label{eq6.14}
\end{equation}
   Then, one of the following holds:
   \begin{align}
       & \exists t_0 \in \mathbb{R},~
       \exists R > R_0,~
       \forall r> |t_0| + R
       \comma u(t_0, r) =0; \label{eq6.15}
       \\
       &  \exists l \in \mathbb{R}\setminus\{0\},~
       \exists R > R_0,~
       \forall r>  R
       \comma (u_0, u_1)(r) = (Z_{l} (r), 0). 
       \label{eq6.16}
   \end{align}
\end{Lemma}

\begin{proof}
    We begin assuming \eqref{eq6.15} fails. We select $\tilde{R} >0$ sufficiently large that 
    \begin{equation*}
        \int_{|x|> \tilde{R}} |\nabla u_0|^2 + |u_1|^2 \, dx
        \leq \varepsilon_0^2
        .
    \end{equation*}
    As a result, $u(t,r)$ satisfies the assumptions of Lemmas~\ref{Lemma4.7}--\ref{Lemma4.8} on $\{|x|> \tilde{R} + |t| \}$. In particular, there is $\ell \neq 0$ such that for all $t \in \mathbb{R}$ and for all $r \geq \tilde{R} + |t|$, we have
    \begin{equation}
        \left|u (t,r) - \frac{\ell}{r^2} \right| \leq C \frac{|\ell|^3}{r^4} 
        \label{eq6.17}
    \end{equation}
    for some $C > 0$. 
    Furthermore, we may assume that 
    $u$ is positive. By selecting a larger $\tilde{R}$, we may combine \eqref{eq6.17} with the asymptotic estimate \eqref{SS_eq2} in Proposition~\ref{prop_singular_solns}. Thus, we observe that for all $t$,
    \begin{equation}
           \left|u (t,r) - Z_{\ell} (r) \right| \leq  \frac{K}{r^4} 
           \comma 
           r > \tilde{R} + |t|
           \label{eq6.18}
    \end{equation}
        for some $K>0$. 

  Next, proceeding as in the proof of \cite[Proposition~4.6]{DKMM22}, we 
  define
  \begin{equation*}
      M (\rho) = \sup_{\substack{t_0 \in \mathbb{R} \\\sigma \geq \rho + |t_0|}}
      \sigma^2 |u (t_0, r) - Z_{\ell} (r)| 
      \comma \rho \geq \tilde{R}.
  \end{equation*}
\end{proof} 
By \eqref{eq6.18}, we directly observe that for $\rho \geq \tilde{R}$, 
\begin{equation}
    M (\rho) \leq \frac{K}{\rho^2}.
    \label{eq6.19}
\end{equation}
As in the proof of \cite[Proposition~4.6]{DKMM22}, we aim to show that there exists $R \geq \tilde{R}$ so that 
\begin{equation}
    M (2 \rho ) \geq \frac{1}{2} M (\rho)
    \label{eq6.20}
\end{equation}
for all $ \rho \geq R$. To that end, we fix $t_0 \in \mathbb{R}$ and $\rho \geq R$, and we define
 \begin{align*}
     h (t,r) = u(t, r) - Z_{\ell} (r)
      \comma 
     h_{t_0} (t,r) = h(2 t_0 - t, r)
     \comma h_+ (t,r) = \frac{h (t,r) + h_{t_0}(t, r)}{2} 
     .
 \end{align*}
 We note that
  \begin{align}
        \begin{cases}
           & \partial_t^2 h_+ - \Delta h_+ = \tilde{F} (t,r) \\
           & \vec{h}_+ (t_0) = (u (t_0) - Z_{\ell}, 0)
        \end{cases}
        \label{eq6.22}
    \end{align}
where the nonlinear term $\tilde{F} (t,r)$ on the right hand side is given by
\begin{equation*}
\tilde{F} (t,r)= \frac{F(u(t,r)) + F(u (2 t_0 -t,r))}{2} - F(Z_{\ell} (r))
\comma F(x) = |x|^{p-1} x. 
\end{equation*} 
We let $\sigma \geq \rho + |t_0|$.
As done in \eqref{eq6.5}, we may estimate
$ \Vert \tilde{F}\Vert_{L^1 L^2 (\sigma + |t - t_0|)}$ from above by
\begin{align}
\begin{split}
     \Vert \tilde{F}\Vert_{L^1 L^2 (\sigma + |t - t_0|)}
   & \les 
     \Vert |h|^p + |Z_{\ell}|^{p-1}| h| \Vert_{L^1 L^2 (\sigma + |t - t_0|)}
   \\ & \quad 
   + 
     \Vert |h_{t_0}|^p + |Z_{\ell}|^{p-1} |h_{t_0}| \Vert_{L^1 L^2 (\sigma + |t - t_0|)}
     .
\end{split}
   \label{eq6.21}
\end{align} 
We observe that when $r \geq \sigma + |t-t_0|$, we have
\begin{equation*}
    r \geq \max (\rho + |t|, \rho + |2 t_0 -t|)
\end{equation*}
which, further implies that 
\begin{equation}
    |h_+ (t,r)| + |h_{t_0} (t,r)|+ |h(t,r)|
    \les \frac{M (\rho)}{r^2}. 
\end{equation}
Therefore, the right hand side of \eqref{eq6.21} may be estimated by
\begin{align}
    \begin{split}
    &
        \int_{\mathbb{R}} 
        \left( 
        \int_{\sigma + |t- t_0|} \frac{|M (\rho)|^{2p}}{r^{4p}} r^3 dr
        \right)^{1/2} dt
        +
          \int_{\mathbb{R}} 
        \left( 
        \int_{\sigma + |t- t_0|} |Z_{\ell} (r)|^{2 (p-1)}   \frac{|M (\rho)|^{2}}{r^{4}} r^3 dr
        \right)^{1/2}  dt
    \\ & 
    \quad \les 
    \frac{|M (\rho)|^{p}}{\sigma^{2p -3}}
    + \sup_{r > \sigma} 
    \left| r^2 Z_{\ell} (r)\right|^{p-1} 
    \int_{\mathbb{R}} 
        \left( 
        \int_{\sigma + |t- t_0|} \frac{1}{r^{4 p -4}} \frac{|M (\rho)|^{2}}{r^{4}} r^3 dr
        \right)^{1/2} dt. 
    \end{split}
\end{align}
We recall that our choice of $\tilde{R}$ guarantees that 
$ \left| r^2 Z_{\ell} (r)\right| \les 1$
  for all $r \geq \tilde{R}$, and so
we obtain
\begin{align}
     \Vert \tilde{F}\Vert_{L^1 L^2 (\sigma + |t - t_0|)}
    \les 
    \frac{|M (\rho)|^{p}}{\sigma^{2p -3}}
    +
     \frac{|M (\rho)|}{\sigma^{2p -3}}
     . 
     \label{eq6.23}
\end{align}
On the other hand, applying Lemma~\ref{Lemma3.9} to the system \eqref{eq6.22}, we get
\begin{align*}
\norm{  h_+ (t_0, r) - \frac{\sigma^2}{r^2} h_+ (t_0, \sigma) }_{\dot{H}^1(\sigma)}
\les 
    \Vert  \pi_{\sigma}^{\perp} h_+ (t_0) \Vert_{\dot{H}^1(\sigma)}
    \les   \Vert \tilde{F}\Vert_{L^1 L^2 (\sigma + |t - t_0|)}. 
\end{align*}

We may close the proof by following the same estimates in \cite[Proposition~4.6]{DKMM22}. By radial Sobolev and \eqref{eq6.23}, the previous estimate leads to
\begin{equation*}
    \left|
    h_+ (t_0, 2 \sigma) - \frac{\sigma^2}{{(2 \sigma)}^2} h_+ (t_0, \sigma)
    \right| \leq C \frac{|M (\rho)|}{\sigma^{2p -2}}
    \left(  {|M (\rho)|^{p-1}}
    + 1 \right) 
\end{equation*}
Next, we utilize the fact that $2 \sigma \geq 2 \rho + 2 |t_0| \geq 2 \rho + |t_0|$ implying 
\begin{equation}
    |h_+ (t_0, 2 \sigma)| \leq \frac{M (2 \sigma)}{(2 \sigma)^2}. 
\end{equation}
Recalling the definition of $h_+ (t_0)$, we deduce that 
\begin{equation}
    \frac{1}{4} |h_+ (t_0,  \sigma)|
    =  \frac{1}{4} |u (t_0, \sigma) - Z_{\ell} (\sigma)|
    \leq C \frac{|M (\rho)|}{\sigma^{2p -2}}
    \left(  {|M (\rho)|^{p-1}}
    + 1 \right)  
    + \frac{M (2 \sigma)}{4 \sigma^2}
    .
\end{equation}
Note that the estimate above holds for all $t_0 \in \mathbb{R}$ and $\sigma \geq \rho + |t_0|$. Combined with \eqref{eq6.19}, we then have
\begin{equation}
    M (\rho) \leq M (2 \rho) + \frac{\tilde{C} |M (\rho)|}{\rho^{2p-4}}
\end{equation}
which proves \eqref{eq6.20} for sufficiently large $\rho$. Note that, as in \cite[Proposition~4.6]{DKMM22}, 
 \eqref{eq6.20} and \eqref{eq6.19}  imply that $M (\rho) =0$ for $\rho$ large, which in turn yields \eqref{eq6.16}.

\begin{proof}[Proof of Proposition~\ref{Prop4.5}]
Let $\vec{u} (t,r)$ be as in Proposition~\ref{rigidity proposition}. Then, $\vec{u} (t,r)$ satisfies the hypotheses of Lemma~\ref{Prop4.6} for all $R_0 >0$. Hence, it suffices show that \eqref{eq6.16} cannot be true. We argue by contradiction and assume that \eqref{eq6.16} is true. Denote by
\begin{align*}
   h (t,r) = u(t,r) - Z_{\ell} (r) 
   \comma
   (h_0, h_1) = (u_0  - Z_{\ell}, u_1 )
\end{align*}
for $t \in \mathbb{R}$ and $r>0$. By \eqref{eq6.16}, $(h_0, h_1)$ is compactly supported. Defining $\rho$ as follows
\begin{align}
    \rho = \inf \lbrace 
    \sigma >0: \int_{|x|> \sigma} |\nabla h_0|^2 + h_1^2 \, dx =0  \rbrace 
\end{align}
we note that $\rho < \infty$. On the other hand, $\rho$ cannot be zero because $Z_{\ell} \notin \dot{H}^{s_p}(\mathbb{R}^4)$ and $u_0 \in \dot{H}^{s_p} (\mathbb{R}^4)$, by hypothesis. Therefore, $0 < \rho < \infty$.  

We observe that for all $R>0$, ${h}(t,r)$ solves the equation
\begin{equation}
    \partial_t^2 h - \Delta h = |u|^{p-1} u - |Z_{\ell}|^{p-1} Z_{\ell} \comma r > R+ |t|
    \label{eq6.26}
\end{equation}
with initial data 
\begin{equation*}
    \vec{h}(0, r) = (u_0 (r) - Z_{\ell} (r), u_1(r)) \comma r >R. 
\end{equation*}
We denote the right hand side of \eqref{eq6.26} by $F(t,r)$.  By the fundamental theorem of calculus, we may express $F(t,r)$ as
\begin{equation}
    F(t,r) = \tilde{V} (t,r) h(t,r)
\end{equation}
with
\begin{equation}
    |\tilde{V} (t,r)| \les |h (t,r)|^{p-1} + |Z_{\ell} (r)|^{p-1}. 
\end{equation}
Moreover, by the radial Sobolev and \eqref{SS_eq2}, we get 
\begin{equation}
    |r h(t,r)| \les  \sup_{t \in \mathbb{R}} \norm{ \nabla u (t)}_{L^2(R)} + \frac{1}{r} 
\end{equation}
for all $t \in \mathbb{R}$ and $r > R+ |t|$. 

Next, we invoke \cite[Proposition~4.7]{DKM20} on an outer channel $\{r> \rho - \delta + |t| \}$, where $\delta >0$ is sufficiently small. More specifically, we observe that \eqref{eq6.26} is of the form
\begin{equation*}
    \partial_t^2 h - \Delta h = \tilde{V} h
    \comma r> \rho - \delta +|t|
\end{equation*}
where the potential $\tilde{V}$ obeys
\begin{equation*}
    |\tilde{V} (t,r)| \leq \left( |h(t,r)|^{p-1} + |Z_{\ell} (r)|^{p-1} \right)
    \leq \frac{C_{\rho}}{r^{p-1}}
    \comma r> \rho- \delta + |t| . 
\end{equation*}
As noted in the proof of Proposition~\ref{Proposition1}, the proof of Proposition~4.7 in \cite{DKM20} remains valid in four dimensions. 
As a result, the following lower bound holds either for all $t\geq 0$ or for all $t\leq 0$: 
\begin{equation}
    \int_{ \rho - \delta + |t|}^{\infty}
    \left|\partial_{t, r} h (t,r) \right|^2 r^3 \, dr \geq \frac{1}{8} 
    \int_{ \rho - \delta + |t|}^{\infty}
    \left(
    \left| \nabla h_0 (r) \right|^2 + 
    \left| h_1 (r) \right|^2  
    \right)r^3 \, dr 
    >0.
    \label{eq6.25.2}
\end{equation}
On the other hand, we have
\begin{equation*}
    \lim_{t \to \pm \infty} \int_{R +|t|}^{\infty} |\partial_r Z_{\ell} (r)|^2 \, r^3 dr =0 
\end{equation*}
for all $R>0$. By  \eqref{eq6.25.2}, this implies that 
\begin{equation*}
    \sum_{\pm} \lim_{t \to \pm \infty} 
    \int_{\rho- \delta + |t|}^{\infty} 
    \left|\partial_{t,r} u(t,r)\right|^2 r^3 \, dr >0
\end{equation*}
which contradicts \eqref{eq5.16}. 
 \end{proof}

\section{Rigidity Method Part II: Odd Solutions} \label{sec: rigidity section 2}

As in \cite[Section~5]{DKMM22}, the second part of the rigidity program focuses on the following
result. 

\begin{Proposition}\label{Prop5.1}
   Let $R>0$ and $\vec{u} (t)$ be a radial solution of \eqref{equation4.1} on $\{ r> R+ |t| \}$ with data $(0, u_1) \in \mathcal{H} (R)$. Assume that
   
   \begin{equation}
       \sum_{\pm} \lim_{t \to \pm \infty}  \int_{|x| > R + |t|} |\nabla_{t,x} u (t, x)|^2 \, dx =0.
       \label{eq7.1}
   \end{equation}
   Then, we have for all $t \in \mathbb{R}$ 
   \begin{equation}
     u(t,r) =0
    \comma r \geq R+ |t| .
       \label{eq7.2}
   \end{equation}
\end{Proposition}
Assuming Proposition~\ref{Prop5.1}, we have:
\begin{proof}[Proof of Proposition~\ref{rigidity proposition}] 
    Suppose $u(t,r)$ is as in Proposition~\ref{rigidity proposition}. By Proposition~\ref{Prop4.5}, there exists $t_0 \in \mathbb{R}$ and $R>0$ so that 
    \begin{equation*}
        u(t_0, r) =0 \quad \quad \mbox{for all }~ r > |t_0| +R.
    \end{equation*} 
   We may then apply Proposition~\ref{Prop5.1} to $u(t + t_0, r)$, which satisfies the hypotheses above with a larger $\tilde{R} >  |t_0| + R$. Then, the conclusion of Proposition~\ref{Prop5.1} implies that of Proposition~\ref{rigidity proposition}. 
\end{proof}

We now turn to the proof of Proposition~\ref{Prop5.1}. 
\begin{Lemma}\label{Lemma5.2}
 Let $\vec{u} (t)$ be as in Proposition~\ref{Prop5.1}. Then there exists a large constant $R' \geq R$ such that 
 \begin{equation}
     r u_1 \in \dot{H}^1 (R') \label{eq7.3}
 \end{equation}
  and 
  \begin{equation}
      \frac{1}{\rho} \norm{u_1}_{\dot{H}^1 (\rho)} \simeq |u_1 (\rho)|
      \comma \mbox{for all }~ \rho \geq R'.
      \label{eq7.4}
  \end{equation}
Furthermore, $(\partial_t u, \partial_t^2 u)$ is the restriction to $\{r> R' + |t| \}$ of a $C(\mathbb{R}, \dot{H}^1 \times L^2)$ function, and it satisfies
\begin{equation}
     \sum_{\pm} \lim_{t \to \pm \infty}  \int_{|x| > R' + |t|} |\nabla_{t,x} \partial_t u (t, x)|^2 \, dx =0.
       \label{eq7.5}
\end{equation}
\end{Lemma}

\begin{proof}[Proof of Lemma~\ref{Lemma5.2}]
    We remark that  \eqref{eq7.3} and \eqref{eq7.4}  follow verbatim from the proof of Steps 1--2 in \cite[Lemma~5.2]{DKMM22} . The proof of \eqref{eq7.5} is also similar to that of Step~3 in \cite[Lemma~5.2]{DKMM22}. However, since the proof of this step depends on the equation, namely \eqref{equation4.1} in this case, we justify the details behind the proof of \eqref{eq7.5} below. 

   As in \cite[Step~3]{DKMM22}, let $\epsilon_0 >0$ to be a small constant to be determined below. We may guarantee that
   \begin{equation}
       \norm{r \partial_r u_1}_{L^2 (R')} + \norm{u_1}_{H^1 (R')} \leq \epsilon_0
   \end{equation}
   by selecting a larger $R'$ if needed. 

Next, we update the initial data by extending $u_1$ as follows:  
\begin{align*}
\tilde{u}_1 (r) = \begin{cases}
     u_1 (R') &  \inon{if }~ r< R', \\
     u_1 (r)  &  \inon{if }~ r \geq R'. 
\end{cases}
\end{align*}
   By the radial Sobolev embedding, we get 
   \begin{equation*}
       |u_1 (R')| \les \frac{1}{(R')^2} 
       \norm{r \partial_r u_1}_{L^2 (R')}
       \les \frac{\epsilon_0}{(R')^2}
   \end{equation*}
   which implies that $\norm{\tilde{u}_1}_{H^1} \les \epsilon_0$. 
   We then consider the solution 
   $(\tilde{u}, \partial_t \tilde{u})$ to the Cauchy problem \eqref{eq7.9} with $V=0$ and $r_0= R'$, and data $(0, \tilde{u}_1)$. Namely, we have 
      \begin{align}
    \begin{split}
        & \partial_t^2 w - \partial_r^2 w - \frac{3}{r} \partial_r w = \chi_{R'} |w|^{p-1} w 
        \comma  \\
        & \evaluat{(w , \partial_t w)}_{t=0} = (0, \tilde{u}_1) .
    \end{split}
\label{eq7.6}   
\end{align}   
As discussed in Remark~\ref{R5.6}, the Strichartz estimates in \eqref{eq5.7} yield
\begin{equation}
    \norm{\tilde{u}}_{L^3 L^6} 
    + 
  \sup_{t \in \mathbb{R}}  \norm{(\tilde{u}, \partial_t \tilde{u} )}_{\mathcal{H}} 
  \les \norm{\tilde{u}_1}_{L^2} + \frac{\epsilon_0^p}{(R')^{p-3}}
  \les \epsilon_0. 
\end{equation}
Similarly, we may bound 
\begin{equation*}
      \norm{\partial_t \tilde{u}}_{L^3 L^6} 
      \les \norm{\tilde{u}_1}_{\dot{H}^1} 
      + \frac{\epsilon_0^p}{(R')^{p-3}}
      \les \epsilon_0. 
\end{equation*}
By Duhamel's principle, we express the solution  $(\tilde{u}, \partial_t \tilde{u})$ as
\begin{equation*}
    \tilde{u} (t) = (- \Delta)^{- 1/2} \sin (t \sqrt{- \Delta}) \tilde{u}_1 
    + \int_0^t \frac{\sin ((t-s) \sqrt{- \Delta } )}{\sqrt{- \Delta}} \chi_{R'} |\tilde{u}|^{p-1} \tilde{u} (s) \, ds
\end{equation*}
and 
\begin{equation*}
   \partial_t \tilde{u} (t) = \cos (t \sqrt{- \Delta}) \tilde{u}_1 + \int_0^t \cos ((t-s) \sqrt{- \Delta}) \chi_{R'} |\tilde{u}|^{p-1} \tilde{u} (s) \, ds. 
\end{equation*}
Since $\tilde{u}_1 \in H^{1}$, we observe that the solution to the following Cauchy problem
  \begin{align}
    \begin{split}
        & \partial_t^2 h - \partial_r^2 h - \frac{3}{r} \partial_r h= \chi_{R'} |\tilde{u}|^{p-1} h 
        \comma  \\
        & \evaluat{(h, \partial_t h)}_{t=0} = (\tilde{u}_1, 0) \in \dot{H}^{1} \times L^2 .
    \end{split}
\label{eq7.7}   
\end{align}   
is indeed given by
$(\partial_t \tilde{u}, \partial_t^2 \tilde{u})$. Note that
\begin{equation*}
    \norm{|\tilde{u}|^{p-1} \partial_t \tilde{u}}_{L^1 L^2 (R')}
    \les 
    \frac{1}{(R')^{p-3}} \norm{\tilde{u}}_{\dot{H}^1 (R')}^{p-3}
    \norm{\tilde{u}}^2_{L^3 L^6 (R')} 
    \norm{\partial_t \tilde{u}}_{L^3 L^6 (R')}
    < \infty.
\end{equation*}
We then close the proof as in    \cite[Lemma~5.2]{DKMM22}. By finite speed of propagation, we observe that
\begin{equation*}
    u (t,r) = \tilde{u} (t,r)
\end{equation*}
for all $r > R' + |t|$. By \eqref{eq7.1}, this implies that 
  \begin{equation*}
       \sum_{\pm} \lim_{t \to \pm \infty}  \int_{|x| > R' + |t|} |\nabla_{t,x} \tilde{u} (t, x)|^2 \, dx =0.
   \end{equation*}
Invoking \cite[Proposition~A.1]{DKMM22}, we observe that $u(t,r)$ and $\tilde{u}(t,r)$ satisfy the zero limits in \cite[Proposition~A.1, (A.2)--(A.3)]{DKMM22} with $G=0$. As a result, by (A.4)--(A.5), we obtain \eqref{eq7.5}. 
\end{proof}

\begin{proof}[Proof of Proposition~\ref{Prop5.1}] We may prove Proposition~\ref{Prop5.1} in the spirit of  Subsection~5.3 in \cite{DKMM22}
and complete the second part of the rigidity method. 

\textit{Step 1. Approximate non-radiative solution.}  For a contradiction, we assume that $u_1$ is not identically zero for $r>R$. By Lemma~\ref{Lemma5.2}, $u_1 (\rho)$ is never zero when $\rho > R'$. Due to the continuity of $\partial_t u (t)$, we deduce that $u_1$ does not change sign. Without loss of generality, we may assume that for sufficiently large $R$
\begin{equation*}
    u_1 (r) >0 \comma \mbox{ for all }~ r>R.
\end{equation*}
By the radial Sobolev inequality and \eqref{eq7.3}, we also have
\begin{equation}
    \lim_{r \to \infty} r^2 u_1 (r) =0.
    \label{eq7.11}
\end{equation}
Next, we define 
\begin{equation}
    a (t,r) = \frac{t}{r^2 (\log r)^{1/2}},
    \quad  \quad 
    a_{\lambda} (t,r) = \lambda a (\lambda t, \lambda r). 
\end{equation}
We fix $\rho  \gg 1$ and set $\lambda = \lambda (\rho)$ so that  
\begin{align}
\begin{split}
     u_1 (\rho) = \partial_t a_{\lambda} (0, \rho )
     = \frac{1}{\rho^2 (\log (\rho \lambda ))^{1/2}} .
    \end{split}
    \label{eq7.10}
\end{align}
   We note that, as in \cite{DKMM22}, there is a unique $\lambda = \lambda (\rho) > \frac{1}{\rho}$ that satisfies \eqref{eq7.10}. Furthermore, \eqref{eq7.11} implies that 
   \begin{equation*}
       \lim_{\rho \to \infty} \rho \lambda (\rho) = \infty.
   \end{equation*}
By selecting $\rho \gg 1$ sufficiently large, we may assume, without loss of generality, that $\log (\rho \lambda (\rho)) >1$. We then consider
\begin{equation}
    w (t,r) : = \partial_t a_{\lambda} (t,r) - \partial_t u(t,r)
    = \frac{1}{r^2 (\log (r \lambda))^{1/2}} 
    - \partial_t u (t,r). 
\end{equation}
As we have seen in the proof of Lemma~\ref{Lemma5.2}, $\partial_t u(t,r)$ solves \eqref{eq7.7} and obeys the non-radiative property \eqref{eq7.5}. 
Also, we have
\begin{equation*}
    \Delta \left( \frac{1}{r^2 (\log (\lambda r))^{1/2}} \right) 
    =  
    \frac{1}{r^4 (\log (\lambda r))^{3/2}} + \frac{3}{4 r^4 (\log (\lambda r))^{5/2}}       .      
\end{equation*}
As a result, we observe that
$w (t,r)$
is a solution to the equation
\begin{align}
    \partial_{tt} w - \partial_r^2 w - \frac{3}{r} \partial_r w = - p |u|^{p-1} \partial_t u
    - \frac{1}{r^4 (\log (r \lambda))^{3/2}} - \frac{3}{4 r^4 (\log (r \lambda))^{5/2}}
    \label{eq7.12}
\end{align}
on the outer set $\{ r> \rho + |t| \}$, supplemented with the initial data
\begin{equation}
    w (0,r) = \left( \frac{1}{r^2 \log (\lambda r)^{1/2}} - u_1 (r), 0 \right).
    \label{eq7.13}
\end{equation}
Also, by \eqref{eq7.10}, we find that $w(0, \rho) = 0$, which implies that $\pi_{\rho}^{\perp} w (0,r) = w (0, r)$. 

\textit{Step 2. Channels of energy estimates.} 
 Noting that the data in \eqref{eq7.13} is of the form $(w_0, 0)$, we apply Lemma~\ref{Lemma3.9} to \eqref{eq7.12}--\eqref{eq7.13}. We denote by $f_w$  the right-hand side of \eqref{eq7.12},  and we get
 \begin{align}
 \norm{w (0)}_{\dot{H}^1 (\rho)}
= \norm{\pi_{\rho}^{\perp} w(0)}_{\dot{H}^1 (\rho)} 
\les \lim_{t \to \infty} \norm{\nabla w (t)}_{L^2 (\rho + |t|)} 
+ \norm{f_w}_{L^1 L^2 (\rho )}
\label{eq7.14}
 \end{align}
 following the notation introduced in Subsection~\ref{subsec: notation}. 
 The first term on the right-hand side above is simply given by
 \begin{align*}
    \lim_{|t| \to \infty} \norm{\nabla w (t)}_{L^2 ( \rho + |t|)} 
    = \lim_{|t| \to \infty} \norm{\nabla (\partial_t a_{\lambda} - \partial_t u)}_{L^2 ( \rho + |t|)}. 
\end{align*}
We directly calculate that
\begin{equation*}
    \partial_r \partial_t a_{\lambda} (r) = 
    -\frac{2}{r^3 (\log (\lambda r))^{1/2}} - \frac{1}{2 r^3 (\log ( \lambda r))^{3/2}}
\end{equation*}
and 
\begin{equation*}
 \norm{\nabla (\partial_t a_{\lambda} )}_{L^2 (|t| + \rho)} 
 = \bigO \left( \frac{1}{|t| + \rho} \right)
 .
\end{equation*}
The above estimate, combined with \eqref{eq7.5}, shows that
\begin{align*}
    \lim_{|t| \to \infty} \norm{\nabla w (t)}_{L^2 ( \rho + |t|)} 
    =0.     
\end{align*}
By the exterior Strichartz estimates, we also have 
\begin{align}
    \begin{split}
        \norm{w}_{L^2 L^8 (\rho)} 
        + \norm{w}_{L^3 L^6 (\rho)}
        + \sup_{t \in \mathbb{R}} \norm{\vec{w} (t)}_{\mathcal{H} (\rho + |t|)}
        \les 
         \norm{w (0)}_{\dot{H}^1 (\rho)}
         + 
          \norm{f_w}_{L^1 L^2 (\rho )} . 
    \end{split}
    \label{eq7.15}
\end{align}
By \eqref{eq7.14}, we note that it is sufficient to estimate $\norm{f_w}_{ L^1 L^2 (\rho )}$ to control the right hand side of \eqref{eq7.15}. 
Recalling that $f_w$ denotes the right-hand side of \eqref{eq7.12}, we estimate
\begin{align}
    \begin{split}
        \norm{f_w}_{L^1 L^2 ( \rho )}
        & \leq 
        \norm{p |u|^{p-1} w }_{L^1 L^2 (\rho )}
        + 
          \norm{p |u|^{p-1} |\partial_t a_{\lambda}| }_{L^1 L^2 (\rho )}
          \\ & \quad \quad
          +   \norm{ \frac{1}{r^4 (\log (r \lambda))^{3/2}} }_{L^1 L^2 (\rho )}
          \\ & \quad \quad
          +   \norm{ \frac{3}{ 4 r^4 (\log (r \lambda))^{5/2}} }_{L^1 L^2 (\rho )} 
          \\ & = I_1 + I_2 + I_3 + I_4.
    \end{split}
\end{align}
We begin with $I_1$: 
\begin{align*}
    I_1 \les \int_0^{\infty} 
    \left( 
    \int_{|t|+ \rho}^{\infty} 
    |u|^{2p -6} |u|^4 |w|^2 r^3 \, dr
    \right)^{1/2} \, dt.
\end{align*}
Using the radial Sobolev inequality on the first factor and applying Holder's inequality to $|u|^4 |w|^2$, we bound $I_1$ with
\begin{align}
    I_1 \les \frac{1}{\rho^{p-3}} \norm{u}^{p-3}_{L^{\infty} \dot{H}^1 (\rho)}
    \norm{u}^2_{L^3 L^6 (\rho)} 
    \norm{w}_{L^3 L^6 (\rho)}
    \les \frac{1}{\rho^{p-3}}   \norm{w}_{L^3 L^6 (\rho)}
    .
    \label{eq7.42}
\end{align}
We note that the factors depending on the norms of $u$ are embedded in the implicit constant by the Strichartz estimates on $u$.

In order to treat $I_2$, we write
\begin{equation*}
    p |u|^{p-1} \partial_t a_{\lambda} 
    = p \left( |u|^{p-1} - |a_{\lambda}|^{p-1} \right)
    \partial_t a_{\lambda} 
    + p |a_{\lambda}|^{p-1} \partial_t a_{\lambda} 
\end{equation*}
and express the difference on the first term above by the fundamental theorem of calculus
\begin{equation*}
    F(u) - F( a_{\lambda}) = 
    (u - a_{\lambda}) \int_0^1 F' ( s u + (1-s) a_{\lambda} ) \, ds
\end{equation*}
with $F (x) = |x|^{p-1}$. This way, we may estimate $I_2$ through the two main terms below:
\begin{align}
\begin{split}
    I_2 & \les
    \norm{ ( |u|^{p-2} + |a_{\lambda}|^{p-2} ) |u - a_{\lambda}| \, |\partial_t a_{\lambda}|}_{L^1 L^2 (\rho )}
    + 
    \norm{ |a_{\lambda}|^{p-1} \, \partial_t a_{\lambda}   }_{L^1 L^2 (\rho )}
    \\ & 
    = I_{21} + I_{22}. 
    \label{eq7.37}
    \end{split}
\end{align}
As $I_{21}$ is naturally split into two terms, we get
\begin{align}
    \begin{split}
        I_{21} & \leq 
        \norm{  |u|^{p-2}  |u - a_{\lambda}| \, |\partial_t a_{\lambda}|}_{L^1 L^2 (\rho )}
    + 
    \norm{  |a_{\lambda}|^{p-2}  |u - a_{\lambda}| \, |\partial_t a_{\lambda}|}_{L^1 L^2 (\rho )}
    \\
    & = 
    I_{211} + I_{212}.
    \label{eq7.34}
    \end{split}
\end{align}
For $I_{211}$, recalling $w = \partial_t (a_{\lambda} - u)$, we first note that
\begin{align}
   |u - a_{\lambda}| \, |\partial_t a_{\lambda}| 
   = \frac{1}{r^2 (\log (\lambda r))^{1/2}} \left| 
   \int_0^t w ( \tau, r) \, d\tau \right| . 
   \label{eq7.32}
\end{align}
By Lemma~C.1 in \cite{DKMM22}, the operator
\begin{equation*}
    (A w) (t,r) = \frac{1}{r^2} \int_0^t w (\tau, r) \, d\tau 
\end{equation*}
maps into $L^2L^{8/3} (\rho)$ and satisfies
\begin{equation}
    \norm{ A w }_{L^2 L^{8/3} (\rho)} \les \norm{w}_{L^2 L^8 (\rho)}. 
\end{equation}
Using $\log (\lambda r) > 1$, we apply the estimate above in \eqref{eq7.32}, and we obtain
\begin{equation}
    \norm{ |u - a_{\lambda}| \, |\partial_t a_{\lambda}| }_{L^2 L^{8/3} (\rho)} \les \norm{ w}_{L^2 L^8 (\rho)} 
    \label{eq7.33}
\end{equation}
which leads to 
\begin{equation}
    I_{211} \les \norm{ |u|^{p-2}}_{L^2 L^8 (\rho)} 
    \norm{w}_{L^2 L^8 (\rho)}. 
\end{equation}
Next, by  the radial Sobolev inequality, we observe that
\begin{align*}
   \norm{ |u|^{p-2}}_{L^2 L^8 (\rho)}  & = 
    \left( \int_0^{\infty} 
    \left( 
    \int_{|t|+ \rho}^{\infty}
    |u|^{8 (p-3)} |u|^8  
    r^3 \, dr
    \right)^{1/4} \, dt \right)^{1/2}
    \\ & \leq
         \frac{C}{\rho^{p-3}} \norm{u}^{p-3}_{L^{\infty} \dot{H}^1 (\rho)}
    \norm{u}_{L^2 L^8 (\rho)} 
    \\ & \les \frac{C(u, p)}{ \rho^{p-3}}
\end{align*}
where we included the admissible norms of $u$ in $C(u,p)$. Then, we obtain
\begin{align}
     I_{211} 
  \les  \frac{1}{\rho^{p-3}}  \norm{w}_{L^2 L^8 (\rho)}. 
  \label{eq7.16}
\end{align}

Next, we turn to $I_{212}$ in \eqref{eq7.34}. Once again, recalling 
$\log (\lambda r) > 1$, we observe that for $t>0$ and $r> \rho + t$
\begin{align*}
    |a_{\lambda}|^{p-2} \, |\partial_t a_{\lambda}|
     = \left( \frac{t^{p-2}}{r^{p-2}} \right) \frac{1}{r^p (\log (\lambda r))^{\frac{p-1}{2}}}
     \leq \frac{1}{r^p} .
\end{align*}
Hence, it is sufficient to estimate the integral
\begin{equation}
    \norm{\frac{1}{r^p} (u - a_{\lambda}) }_{L^1L^2 (\rho)}
    = \int_0^{\infty} \left( 
    \int_{\rho + t}^{\infty} \frac{1}{r^{2p}} \left(\int_0^t w (\tau, r) \, d\tau \right)^2  r^3 \, dr
    \right)^{1/2} \, dt. 
    \label{eq7.35}
\end{equation}
We may control  the inner integral above by Minkowski's inequality followed by Holder's inequality. We obtain
\begin{align*}
\begin{split}
   \left( 
    \int_{\rho + t}^{\infty} \frac{1}{r^{2p}} \left(\int_0^t w (\tau, r) \, d\tau \right)^2  r^3 \, dr
    \right)^{1/2}
   &  \leq \int_0^t \left(  
    \int_{\rho + t}^{\infty} \frac{w^2 (\tau, r)}{r^{2p}} r^3 \, dr
    \right)^{1/2} \, d\tau
    \\& 
    \leq 
 \int_0^t  
    \norm{{w(\tau, \cdot)}}_{L^6 (\rho + t)}
    \norm{{r^{-p}} }_{L^3 (\rho + t)} \, d\tau. 
    \end{split}
\end{align*}
As we have
\begin{equation*}
    \int_{\rho + t }^{\infty} \frac{1}{r^{3p}} r^3 dr \leq \frac{1}{ (\rho +t)^{3p -4}} \comma
\end{equation*}
we note that the integral on the right hand side of \eqref{eq7.35} 
may be estimated by
\begin{equation}
    \int_0^{\infty} 
    \int_0^t
    \frac{1}{(\rho + t)^{p - \frac{4}{3}}} 
    \left( \int_{\rho +t}^{\infty} |w (\tau, r)|^6\,  r^3 dr \right)^{1/6}  \, d\tau  \, dt. 
    \label{eq7.36}
\end{equation}
Next, we apply Fubini followed by Holder's inequality to bound the above integral as follows:
\begin{align}
    \begin{split}
   &    \int_0^{\infty} 
    \int_{\tau}^{\infty}   
    \frac{1}{(\rho + t)^{p - \frac{4}{3}}} 
    \left( \int_{\rho + \tau }^{\infty} |w (\tau, r)|^6\,  r^3 dr \right)^{1/6}
    \, d\tau  \, dt
    \\ & \indeq \les 
\int_0^{\infty}
    \left( \int_{\rho + \tau }^{\infty} |w (\tau, r)|^6\,  r^3 dr \right)^{1/6}
     \frac{1}{(\rho + \tau)^{p - \frac{7}{3}}} 
     \, d \tau
     \\ & \indeq \les 
     \norm{w}_{L^3 L^6 (\rho)} 
     \left( \int_0^{\infty} \frac{1}{(\rho + \tau)^{\frac{3}{2} (p- \frac{7}{3}) } } \right)^{2/3} . 
    \end{split}
\end{align}
The upper bound on the last line completes the estimate for the integral in \eqref{eq7.36}, which in turn leads to
\begin{equation}
  I_{212} \les  \frac{1}{\rho^{p-3}} \norm{w}_{L^3 L^6 (\rho)}. 
  \label{eq7.39}
\end{equation}
Next, we treat $I_{22}$ in \eqref{eq7.37}. Expanding the norm for $I_{22}$, we get 
\begin{align*}
    I_{22} = \int_0^{\infty} t^{p-1} \left( 
    \int_{\rho + t}^{\infty} 
    \frac{1}{r^{4p}} \frac{1}{(\log (\lambda r))^p} \, r^3  dr
    \right)^{1/2} dt. 
\end{align*}
The following integration by parts identity 
for $k_1 > 1$, $k_2 > -1$  
\begin{align}
    \int_{\sigma}^{\infty} \frac{1}{r^{k_1} \log^{k_2} (\lambda r)} \, dr 
    =    \frac{C_{k_1}}{r^{k_1 -1}} 
    \frac{1}{\left( \log (\lambda r) \right)^{k_2}}\biggr\rvert_{r= \sigma}^{\infty}
    +
    C_{k_1, k_2} \int_{\sigma}^{\infty}
    \frac{1}{r^{k_1} (\log (\lambda r))^{k_2 +1}} \, dr
    \label{eq7.40}
\end{align}
will be useful in our estimates below. Utilizing this identity, we directly estimate
\begin{align}
    I_{22} &\les 
    \int_0^{\infty} \frac{t^{p-1}}{(\rho +t)^{2p -2}} \frac{1}{(\log (\lambda \rho))^{p/2}} \, dt
    + 
    \int_0^{\infty} t^{p-1} 
    \left(  \int_{\rho +t }^{\infty} 
    \frac{1}{r^{4p -3}} \frac{1}{(\log (\lambda r))^{p+1}} \, dr
    \right)^{1/2} dt
    \\
    & \les 
    \frac{1}{\rho^{p-2}} \frac{1}{(\log (\lambda \rho))^{p/2}} 
    + 
    \frac{1}{\rho^{p-2}} \frac{1}{(\log (\lambda \rho))^{(p+1)/2}} .
    \label{eq7.38}
\end{align}
Collecting the estimates in \eqref{eq7.16}, \eqref{eq7.39}, and \eqref{eq7.38}, we control $I_2$ by
\begin{align}
\begin{split}
     I_2   \les \frac{1}{\rho^{p-3}} (  \norm{w}_{L^2 L^8 (\rho)}   + \norm{w}_{L^3 L^6 (\rho)} )
    + \frac{1}{\rho^{p-2}} \frac{1}{(\log (\lambda \rho))^{p/2}} .
\end{split}
    \label{eq7.41}
\end{align}

For $I_3$, the identity in \eqref{eq7.40} leads to an upper bound of the following form:
\begin{align*}  
\begin{split}
    I_3 & = \int_0^{\infty}
    \left( 
    \int_{\rho + |t|}^{\infty}
    \frac{1}{r^8 \log^3 (r \lambda)} r^3\, dr
    \right)^{1/2} \, dt
    \\ &
    \les 
     \int_0^{\infty}
    \left( 
    \frac{1}{{(\rho + |t|)}^4 \log^3 ((\rho + |t|)  \lambda)} 
    \right)^{1/2} \, dt
    +
    \int_0^{\infty}
    \left( 
    \int_{\rho + |t|}^{\infty}
    \frac{1}{r^5 \log^4 (r \lambda)}\, dr
    \right)^{1/2} \, dt
    .
    \end{split}
\end{align*}
We may then calculate that 
\begin{equation}
    I_3 \les \frac{1}{\rho (\log (\rho \lambda))^{3/2}}. 
    \label{eq7.17}
\end{equation}
Similarly, we get
\begin{align}
    I_4 = \int_0^{\infty}
    \left( 
    \int_{\rho + |t|}^{\infty}
    \frac{1}{r^8 \log^5 (r \lambda)} r^3\, dr
    \right)^{1/2} \, dt
    \les 
    \frac{1}{\rho (\log (\rho \lambda))^{5/2}}.
\label{eq7.18}
\end{align}
Finally, combining \eqref{eq7.42}, \eqref{eq7.41}, and \eqref{eq7.17}--\eqref{eq7.18} with  $p>3$ and $\log (\lambda \rho) >1$, we  obtain
\begin{align*}
     \begin{split}
        \norm{f_w}_{L^1 L^2 ( \rho )}
        & \les
        \frac{1}{\rho^{p-3}} 
  \norm{w}_{L^3 L^6 (\rho)}
  + \frac{1}{\rho^{p-3}} 
  \norm{w}_{L^2 L^8 (\rho)}
  + \frac{1}{\rho (\log (\rho \lambda))^{3/2}}
    \end{split}
\end{align*}
which lets us use the Strichartz estimates  \eqref{eq7.15} as follows. We estimate
\begin{align}
    \begin{split}
 & \norm{w}_{L^2 L^8 (\rho)} 
        + \norm{w}_{L^3 L^6 (\rho)}
        + \sup_{t \in \mathbb{R}} \norm{\vec{w} (t)}_{\mathcal{H} (\rho + |t|)}
   \\  & \quad \quad \quad \quad  
   \les 
        \frac{1}{\rho^{p-3}} 
  \norm{w}_{L^3 L^6 (\rho)} 
  + \frac{1}{\rho^{p-3}} 
  \norm{w}_{L^2 L^8 (\rho)} 
  + \frac{1}{\rho (\log (\rho \lambda))^{3/2}}
    \end{split}
    \label{eq7.19}
\end{align}
Letting $\rho \gg 1$, we ensure that
\begin{align}
    \norm{w}_{L^2 L^8 (\rho)} 
        + \norm{w}_{L^3 L^6 (\rho)}
        + \norm{{w} (0)}_{\dot{H}^1 (\rho )}
        \les 
        \frac{1}{\rho (\log (\rho \lambda))^{3/2}}. 
\end{align}
By \eqref{eq7.13}, we arrive at 
\begin{equation}
    \norm{u_1 - \frac{1}{r^2 (\log (r \lambda))^{1/2}}}_{\dot{H}^1 (\rho)} \les \frac{1}{\rho (\log (\rho \lambda))^{3/2}} .
    \label{eq7.20}
\end{equation}
 
\textit{Step 3.} As in  \cite[\textit{Step~3}]{DKMM22}, we aim to show that the following differential inequality
\begin{equation}
    \left| 
    \norm{u_1}^4_{L^4 (\rho)} - \frac{\rho^4 (u_1 (\rho))^4}{4}
    \right| \les \rho^8 (u_1 (\rho))^6
    \label{eq7.23}
\end{equation}
holds for sufficiently large $\rho$. 

First, the Sobolev inequality applied to \eqref{eq7.20} implies
\begin{equation}
    \left| 
    \norm{u_1}_{L^4 (\rho)} - 
    \norm{\frac{1}{r^2 (\log (r \lambda))^{1/2}}}_{L^4 (\rho)}
    \right| \les 
    \frac{1}{\rho (\log (\rho \lambda))^{3/2}}
. \label{eq7.21}
\end{equation}
After an integration by parts, we may directly calculate
\begin{align*}
    \norm{\frac{1}{r^2 (\log (r \lambda))^{1/2} }}^4_{L^4 (\rho)}
    & = 
    \int_{\rho}^{\infty} 
    \frac{1}{r^8 \log^2 (r \lambda)} r^3 \, dr
    \\ & 
    = 
    \frac{1}{4 \rho^4 \log^2 (\rho \lambda)}
    + \mathcal{O} \left( \frac{1}{\rho^4 \log^3 (\rho \lambda)} \right)
.
\end{align*}
Combining this with the estimate in \eqref{eq7.21}, we obtain
\begin{equation}
    \left|
     \norm{u_1}^4_{L^4 (\rho)} - 
     \frac{1}{4 \rho^4 \log^2 (\rho \lambda)}
    \right| \les \frac{1}{\rho^4 \log^3 (\rho \lambda)}
    .
    \label{eq7.22}
\end{equation}
Lastly, we close this step by using \eqref{eq7.10}. More specifically, inserting 
$u_1 (\rho) = \frac{1}{\rho^2 (\log (\rho \lambda))^{1/2}}$ into \eqref{eq7.22}, we obtain the estimate in \eqref{eq7.23}.

\textit{Step~4. Closing the proof of Proposition~\ref{Prop5.1}.} In this step, as in the last step of \cite[Proposition~5.1]{DKMM22}, we show that
\begin{equation}
    u_1 (\rho) \geq \frac{1}{C \rho^2 (\log \rho)^{1/2}}
    \quad \quad 
    \mbox{ as }~ \rho \to \infty. 
    \label{eq7.24}
\end{equation}
 Since the lower bound above does not belong to $L^2 (r^3 dr)$, we deduce that $u_1 (r) \notin L^2 (R)$ by \eqref{eq7.24}, contradicting the hypothesis of Proposition~\ref{Prop5.1} on $u_1$. Going back to the beginning of the proof, we deduce that $u_1$ is identically zero for $r> R$. 

To verify the lower bound in \eqref{eq7.24}, we argue as follows. Denote by
\begin{equation}
    f (\rho) = \norm{u_1}^4_{L^4 (\rho)}. 
    \label{eq7.26}
\end{equation}
By \eqref{eq7.22}, we have
\begin{equation}
    \rho^4 f(\rho) \les 
     \rho^8 u_1^4 (\rho) + \rho^{12} u_1^6 (\rho) .
    \label{eq7.25}
\end{equation}
We note that the right-hand side above approaches zero as $\rho \to \infty$ by \eqref{eq7.11}. Also, we observe that
\begin{equation}
    f' (\rho) = \frac{d}{d\rho} 
    \left( 
    \int_{\rho}^{\infty} u_1^4 (r) r^3 \, dr
    \right)
    = - u_1^4 (\rho) \rho^3. 
    \label{eq7.27}
\end{equation}
Combined with the estimate in \eqref{eq7.22}, we then obtain the following differential estimate for $f$:
\begin{equation}
    \left|f(\rho) + \frac{\rho}{4} f' (\rho) \right| \les f(\rho)^{3/2} \rho^2 .
    \label{eq7.28}
\end{equation}
Next, we define $g (\rho) = \rho^4 f(\rho)$ and rewrite the previous estimate for $g$. We get
\begin{equation}
    \left|
    \frac{g' (\rho)}{g(\rho)^{3/2}}
    \right| \les \frac{1}{\rho}
\end{equation}
which, upon integrating from $\rho_0$ to $\rho$, leads to
\begin{equation}
    \frac{1}{(g(\rho))^{1/2}} \leq C \log (\rho) 
\end{equation}
for large $\rho$, where $C$ is a positive constant depending on $g(\rho_0)$. In terms of $f$, the previous inequality implies that
\begin{equation}
    \rho^4 f(\rho) \geq \frac{1}{C \log^2 (\rho)}. 
\end{equation}
Recalling  \eqref{eq7.25}, for large $\rho$, we obtain 
\begin{equation}
    \rho^8 u_1^4 (\rho) 
    \geq 
    \frac{1}{{C}_1 \log^2 (\rho)} - C_2 \rho^{12} u_1^6 (\rho)
    \label{eq7.29}
\end{equation}
where $C_1, C_2$ are positive constants.
  Noting that $\rho^{12} u_1^6 (\rho)$  approaches zero much faster than the leading term above, we arrive at the asymptotic lower bound inequality in \eqref{eq7.24}. 
\end{proof}

\section*{Acknowledgments} 
CEK was partially supported by NSF grant DMS-2153794.


\begin{thebibliography}{CK}
\bibitem{BG99}
H.~Bahouri and P.~Gerard. 
\newblock High frequency approximation of solutions
to critical nonlinear wave equations. 
{\em Amer.~J.~Math.},
121: 131--175, 1999. 

\bibitem{BS98}
H.~ Bahouri and J.~Shatah. 
\newblock Decay estimates for the critical semilinear wave equation.
{ \em Ann. Inst. H. Poincare Anal. Non Lineaire}, 
15 (6): 783–789, 1998.

\bibitem{B10}
A.~Bulut
\newblock Maximizers for the Strichartz inequalities for the wave equation.
\newblock{\em Differential Integral Equations}, 
23: 1035--1072, 2010.

\bibitem{B12}
A.~Bulut
\newblock Global well-posedness and scattering for the defocusing energy-supercritical 
cubic nonlinear wave equation. 
\newblock{\em J.~Funct.~Anal.},
263(6): 1609--1660, 2012.

\bibitem{B14}
A.~Bulut
\newblock {The radial defocusing energy-supercritical cubic nonlinear
              wave equation in {$\Bbb R^{1+5}$}}. 
\newblock {\em Nonlinearity}, 
27(8): {1859--1877}, 2014. 


\bibitem{B15}
A.~Bulut
\newblock The defocusing energy-supercritical cubic nonlinear wave
              equation in dimension five. 
\newblock {\em Trans. Amer. Math. Soc.}, 
367(9): {6017--6061}, 2015. 

\bibitem{BCLPZ13}
A.~Bulut, M.~Czubak, D.~Li, N.~Pavlovic, and X.~Zhang.
\newblock Stability and Unconditional Uniqueness of Solutions 
for Energy Critical Wave Equations in High Dimensions.
\newblock {\em Comm. Partial Differential Equations},
38(4):
{575--607}, 2013.

\bibitem{BD21}
A.~Bulut and B.~Dodson.
\newblock Global well-posedness for the logarithmically
              energy-supercritical nonlinear wave equation with partial
              symmetry. 
\newblock {\em Int. Math. Res. Not. IMRN},
(8): {5943--5967}, 2021. 

\bibitem{CK23}
G.~Camliyurt and C.E.~Kenig. 
\newblock {Scattering for radial bounded solutions of focusing
              supercritical wave equations in odd dimensions}. 
\newblock {\em Nonlinear Anal.},
(236): Paper No: 113352, 2023.

\bibitem{C03}
T.~Cazenave. 
\newblock{Semilinear Schr\"{o}dinger equations}. 
{Courant Lecture Notes in Mathematics}, {10},
 {New York University, Courant Institute of Mathematical
              Sciences, New York; American Mathematical Society, Providence,
              RI},
{2003}. 



\bibitem{CKSTT08}
J.~Colliander, M.~Keel, G.~Staffilani, H.~Takaoka, and T.~Tao.
\newblock{Global well-posedness and scattering for the energy-critical nonlinear Schr\"{o}dinger equation in $\mathbb{R}^3$}.
\newblock{ \em Ann.~Math.},
167(3): 767--865, 2008.

\bibitem{Co18}
C.~Collot.
\newblock{Type II blow up manifolds for the energy supercritical semilinear wave equation}. 
\newblock{ \em Mem.~Amer.~Math.~Soc.}, 
252(1205): v+163,  2018.

\bibitem{CDKM24}
C.~Collot, T.~Duyckaerts, C.~Kenig, and F.~Merle.
\newblock{Soliton resolution for the radial quadratic wave equation in space dimension 6}. 
\newblock{ \em Vietnam~J.~Math.}, 52 (3): 735--773, 2024.
              


\bibitem{CS96}
G.M.~Constantine and T.H.~Savits.
\newblock{A multi-variate Fa\`{a} di Bruno formula with applications}.
\newblock{ \em Trans.~Amer.~Math.},
348(2): 503--520, 1996. 

\bibitem{DD21}
W.~Dai and T.~Duyckaerts.
\newblock{Self-similar solutions of focusing semi-linear wave equations
              in {$\Bbb R^N$}}.
\newblock{ \em J. Evol. Equ.},
21(4): 4703--4750, 2021.              

\bibitem{DL15}
B.~Dodson and A.~Lawrie. 
\newblock Scattering for radial, semi-linear, super-critical wave
              equations with bounded critical norm.
\newblock {\em Arch.~Ration.~Mech.~Anal.}, 218(3):1459--1529, 2015.

\bibitem{DL15.2}
B.~Dodson and A.~Lawrie. 
\newblock Scattering for the radial 3d cubic wave equation. 
\newblock{\em Anal.~PDE}, 8(2): 467--497, 2015. 

\bibitem{DLMM20}
B.~Dodson,  A.~Lawrie, D.~Mendelson, and J.~Murphy.
\newblock Scattering for defocusing energy subcritical nonlinear wave
              equations.
\newblock{\em Anal.~PDE}, 13(7): 1995--2090, 2020.               

\bibitem{DHKS14}
R.~Donninger, M.~Huang, J.~Krieger, and W.~Schlag.
\newblock Exotic blowup solutions for the $u^5$ focusing wave equation  
in $\mathbb{R}^3$. 
\newblock {\em Michigan Math.~ J.}, 63(3):{451--501}, 2014.

\bibitem{DS16}
R.~Donninger and B.~Sch\"{o}rkhuber.
\newblock {On blowup in supercritical wave equations}.
\newblock {\em Comm. Math. Phys.}, 346(3):{907--943}, 2016.





\bibitem{DKM11}
T.~Duyckaerts, C.~Kenig, and F.~Merle.
\newblock {Universality of blow-up profile for small radial type {II}
              blow-up solutions of the energy-critical wave equation}.
\newblock {\em J. Eur. Math. Soc. (JEMS)}, 13(3):{533--599}, 2011. 

\bibitem{DKM12}
T.~Duyckaerts, C.~Kenig, and F.~Merle.
\newblock {Universality of blow-up profile for small radial type {II}
              blow-up solutions of the energy-critical wave equation: the nonradial case}.
\newblock {\em J. Eur. Math. Soc. (JEMS)}, 14(5):{1389--1454}, 2012. 

\bibitem{DKM12.2}
T.~Duyckaerts, C.~Kenig, and F.~Merle.
\newblock Profiles of bounded radial solutions of the focusing,
              energy-critical wave equation.
\newblock {\em Geom. Funct. Anal.}, 22(3): {639--698}, 2012. 

\bibitem{DKM13}
T.~Duyckaerts, C.~Kenig, and F.~Merle.
\newblock Classification of radial solutions of the focusing,
              energy-critical wave equation. 
\newblock {\em Camb. J. Math.}, 1(1):{75--144}, 2013.

\bibitem{DKM14}
T.~Duyckaerts, C.~Kenig, and F.~Merle.
\newblock Scattering for radial, bounded solutions of focusing
              supercritical wave equations.
\newblock {\em Int.~Math.~Res.~Not.~IMRN},  1:224--258, 2014.


\bibitem{DKM15.2}
T.~Duyckaerts, C.~Kenig, and F.~Merle.
\newblock {Profiles for bounded solutions of dispersive equations, with
              applications to energy-critical wave and {S}chr\"{o}dinger
              equations}. 
\newblock {\em Commun. Pure Appl. Anal.}, 14(4): {1275--1326}, 2015.

\bibitem{DKM15}
T.~Duyckaerts, C.~Kenig, and F.~Merle.
\newblock Global existence for solutions of the focusing wave equation with the compactness property.
\newblock {\em Ann.~ Inst.~ H.~Poincar\'{e} Anal.~Non~Lin\'{e}aire},  33(6):1675--1690, 2016.

\bibitem{DKM20}
T.~Duyckaerts, C.~Kenig, and F.~Merle.
\newblock{Exterior energy bounds for the critical wave equation close to the ground state}.
\newblock{\em Comm.~ Math.~Phys.}, 379 (3): 1113--1175, 2020.

\bibitem{DKMM22}
T.~Duyckaerts, C.~Kenig, Y.~Martel, and F.~Merle.
\newblock{Soliton resolution for critical co-rotational wave maps and radial cubic wave equation}.
\newblock{\em Comm.~ Math.~Phys.}, 391(2): 779--871, 2022.

\bibitem{DKM23}
T.~Duyckaerts, C.~Kenig, and F.~Merle.
\newblock{Soliton resolution for the radial critical wave equation in all odd space dimensions}.
\newblock{\em Acta~Math.}, 230 (1): 1--92, 2023.
       
\bibitem{DR17}
T.~Duyckaerts and T.~Roy.
\newblock {Blow-up of the critical {S}obolev norm for nonscattering
              radial solutions of supercritical wave equations on {$\Bbb
              R^3$}}. 
\newblock {\em Bull. Soc. Math. France}, 145(3): {503--573}, 2017. 

\bibitem{DY18}
T.~Duyckaerts and J.~Yang.
\newblock Blow-up of a critical {S}obolev norm for energy-subcritical and energy-supercritical wave equations. 
\newblock {\em Anal.~PDE}, 11(4):{983--1028}, 2018. 


\bibitem{GV95}
J.~Ginibre and G.~Velo.
\newblock Generalized Strichartz inequalities for the wave equation.
\newblock{\em J.~Funct.~Anal.}, 133(1):50--68, 1995.

\bibitem{G90}
M.G.~Grillakis. 
\newblock Regularity and asymptotic behaviour of the wave equation with a critical
nonlinearity. 
\newblock {\em  Ann. of Math. (2)}, 
132(3): 485–509, 1990.

\bibitem{G92}
M.G.~Grillakis. 
\newblock Regularity for the wave equation with a critical nonlinearity. 
\newblock { \em Comm. Pure Appl. Math.},
45(6): 749–774, 1992.


\bibitem{JL23}
J.~Jendrej and A.~Lawrie.
\newblock{Soliton resolution for the energy-critical nonlinear wave equation in the radial case}.
\newblock{ \em Ann.~PDE}, 9 (2): Paper No: 18, 2023.  


\bibitem{K94}
L.~Kapitanski.
\newblock  Global and unique weak solutions of nonlinear wave equations. 
\newblock { \em Math. Res.
Lett.}, 1(2): 211–223, 1994.

\bibitem{K15}
C.~Kenig.
 \newblock {Lectures on the Energy Critical Nonliner Wave Equation}.
 { CBMS Regional Conference Series in Mathematics}, 122.  
\newblock { \em Published for the Conference Board of the Mathematical
              Sciences, Washington, DC; by the American Mathematical
              Society, Providence, RI}, 2015.             
              
    
\bibitem{KLS} 
C.~Kenig, A.~Lawrie, and W.~Schlag.
\newblock Relaxation of wave maps exterior to a ball to harmonic maps 
for all data.
\newblock {\em Geom.~Func.~Anal.}, 
24 (2), 610--647, 2014. 

\bibitem{KLLS15} 
C.~Kenig, A.~Lawrie, B.~Liu, and W.~Schlag.
\newblock Channels of energy for the linear radial wave equation. 
\newblock{\em Adv.~Math.},
285: 877--936, 2015.

\bibitem{KLLS} 
C.~Kenig, A.~Lawrie, B.~Liu, and W.~Schlag.
\newblock Stable soliton resolution for exterior wave maps in all
              equivariance classes.
\newblock {\em Adv. Math.}, 285:235--300, 2015.

\bibitem{KM11}
C.~Kenig and F.~Merle. 
\newblock Nondispersive radial solutions to energy 
supercritical non-linear wave equations, with applications.
\newblock{Amer.~J.~Math.}, 
133(4):1029--1065, 2011.

\bibitem{KM1}
C.~Kenig and F.~Merle. 
\newblock Global well-posedness, scattering and blow-up 
for the energy critical, focusing, non-linear Schr\"{o}dinger equation in the radial case. 
\newblock{Invent.~Math.},
166(3):645--675, 2006.

\bibitem{KM2}
C.~Kenig and F.~Merle. 
\newblock Global well-posedness, scattering and blow-up 
for the energy critical focusing non-linear wave equation. 
\newblock{Acta.~Math.},
201(2):147--212, 2008. 

\bibitem{KM10}
C.~Kenig and F.~Merle. 
\newblock Scattering for $\dot{H}^{1/2}$ bounded solutions to the 
cubic defocusing nls in $3$ dimensions.
\newblock{Trans.~Amer.~Math.~Soc.},
362(4):1937--1962, 2010. 



\bibitem{KV10}
 {R.~Killip and M.~Visan}.
\newblock{Energy-supercritical NLS: critical {$\dot{H}^s$}-bounds imply scattering}.
\newblock{\em Comm.~ Partial Differential Equations}, 
35(6): 945--987, 2010.  

\bibitem{KV10.2}
 {R.~Killip and M.~Visan}.
\newblock{The focusing energy-critical nonlinear Schr\"{o}dinger equation in dimensions five and higher}. 
\newblock{\em Amer.~J.~Math.},
132(2): 361--424, 2010. 
 

\bibitem{KV11}
 {R.~Killip and M.~Visan}.
\newblock{The radial defocusing energy-supercritical nonlinear wave
              equation in all space dimensions}.
\newblock {\em Proc. Amer. Math. Soc.},
139(5): {1805--1817}, 2011.

\bibitem{KV11.2}
 {R.~Killip and M.~Visan}.
\newblock{The defocusing energy-supercritical nonlinear wave equation in
              three space dimensions}.
\newblock {\em Trans. Amer. Math. Soc.},
363(7): {3893--3934}, 2011.

\bibitem{KV13}
 {R.~Killip and M.~Visan}.
 \newblock{Nonlinear Schr\"{o}dinger equations at critical regularity}. 
 \newblock{\em Evolution equations}, 
 volume~17 of {\em Clay ~Math.~Proc.},
 325--437. 
 Amer. Math. Soc., Providence, RI, 2013. 


\bibitem{KS14}
J.~Krieger and W.~Schlag.
\newblock {Full range of blow up exponents for the quintic wave equation
              in three dimensions}.
\newblock  { \em J. Math. Pures Appl. (9)},
101(6): {873--900}, 2014.

\bibitem{KST09}
{J.~Krieger, W.~ Schlag, and D.~Tataru}.
 \newblock    Slow blow-up solutions for the {$H^1(\Bbb R^3)$} critical
              focusing semilinear wave equation.
 \newblock {\em Duke Math. J.}, 
 147(1): 1--53, 2009. 

\bibitem{L74}
H.A.~Levine.
\newblock Instability and nonexistence of global solutions to nonlinear wave equations
of the form $P u_{tt} = - Au + \mathcal{F} (u)$. 
 \newblock { \em Trans. Amer. Math. Soc.}, 192: 1--21, 1974.

 \bibitem{LSW24}
 {L.~Li, R.~Shen, and L.~Wei}. 
 \newblock  Explicit formula of radiation fields of free waves with applications on channel of energy. 
 \newblock{ \em Anal. PDE}, 17(2): 723--748, 2024. 
  
 
 \bibitem{MRR15}
 F.~Merle, P.~Raphael, and I.~Rodnianski.
 \newblock Type II blow up for the energy supercritical NLS. 
 \newblock {\em Camb.~J.~Math.},
 3(4): 439--617, 2015.
 
 \bibitem{MRRS22}
  F.~Merle, P.~Raphael, I.~Rodnianski, and J.~Szeftel.
  \newblock On blow up for the energy super critical defocusing nonlinear
              {S}chr\"{o}dinger equations.
 \newblock {\em Invent. ~ Math.},
 227(1): 247--413, 2022.
 
 \bibitem{MZ05}
 F.~Merle and H.~Zaag.
 \newblock Determination of the blow-up rate for a critical semilinear wave equation.
 \newblock {\em Math.~Ann.},
 331(2): 395--416, 2005. 
 
  \bibitem{MZ15}
 F.~Merle and H.~Zaag.
 \newblock On the stability of the notion of non-characteristic point and blow-up profile 
 for semilinear wave equations. 
 \newblock {\em Comm.~Math.~Phys.},
 333(3): 1529--1562, 2015. 
 
 
\bibitem{SS93}
J.~Shatah, M.~Struwe. 
\newblock Regularity results for nonlinear wave equations.
\newblock {\em Ann. of Math.}, 
138(3): {503--518}, 1993.

\bibitem{SS94}
J.~Shatah, M.~Struwe. 
\newblock{Well-posedness in the energy space for semilinear wave equations with critical growth}.
\newblock {\em  Int. Math. Res. Notices}, 7:303--309, 1994.



\bibitem{SS98}
J.~Shatah, M.~Struwe. 
\newblock{Geometric Wave Equations}.
\newblock{ \em  Courant Lecture Notes in Mathematics, 2. New York University
Courant Institute of Mathematical Sciences, New York}, 1998.


 
 \bibitem{S13}  
R.~Shen.
\newblock   On the energy subcritical, nonlinear wave equation in {$\mathbb{R}^3$} with radial data.
\newblock {\em Anal. PDE}, 6(8):1929--1987, 2013.

 \bibitem{S20}  
R.~Shen.
\newblock Bounded solutions to an energy subcritical non-linear wave equation on $\mathbb{R}^3$.
\newblock {\em J.~Differential Equations}, 
269(4):3943--3986, 2020.  


\bibitem{S89}
M.~Struwe
\newblock Globally regular solutions to the $u^5$ Klein--Gordon equation. 
\newblock {\em Ann. Scuola
Norm. Sup. Pisa Cl. Sci.}, 
15(3):495--513, 1989.

\bibitem{R17}
C.~Rodriguez
\newblock Scattering for radial energy-subcritical wave equations in dimensions 4 and 5.
\newblock {\em Comm.~ Partial Differential Equations}, 
42 (6):852--894, 2017.


\bibitem{T10}
R. J.~Taggard. 
\newblock Inhomogeneous {S}trichartz estimates.
\newblock  { \em Forum Math.}, 22(5):825--853, 2010. 

\bibitem{Tao06}
T.~Tao. 
\newblock Nonlinear dispersive equations.  {CBMS Regional Conference Series in Mathematics}, 106.
\newblock { \em Published for the Conference Board of the Mathematical
              Sciences, Washington, DC; by the American Mathematical
              Society, Providence, RI}, 2006. 

\bibitem{TVZ08} 
T.~Tao, M.~Visan, and X.~Zhang.
\newblock Minimal-mass blowup solutions of the mass-critical {NLS}.
\newblock {\em Forum Math.}, 20(5): 881--919, 2008.


\bibitem{TVZ07}
 T.~Tao, M.~Visan, and X.~Zhang.
\newblock Global well-posedness and scattering for the defocusing
              mass-critical nonlinear {S}chr\"{o}dinger equation for radial data
              in high dimensions.
\newblock {\em Duke Math.~ J.},  140(1): 165--202, 2007.







\end{thebibliography}
\end{document}